\documentclass[a4paper,12pt]{amsart}

\usepackage{hyperref}

\usepackage[top=3cm,left=2.5cm,right=2.5cm,bottom=3cm]{geometry}
\usepackage{relsize}
\usepackage{lineno}

\usepackage{amsmath,amsthm,pb-diagram,amssymb,comment}
\usepackage{amsfonts,graphicx,color}
\usepackage{enumitem, fancyhdr, dsfont, multicol}
\usepackage{thmtools,cleveref}
\usepackage[normalem]{ulem}
\usepackage{stmaryrd}
\usepackage{textcomp}
\usepackage{mathtools}
\usepackage{fontawesome}
\usepackage{textcomp}
\usepackage[all]{xy}

\usepackage{tikz,float}
\usepackage{mleftright}
\usepackage{mathrsfs}
\usepackage{tikz-cd}
\usepackage{dsfont}
\usepackage{stmaryrd}

\usepackage{caption}
\usepackage{subcaption}

\hypersetup{
    colorlinks=true,
    linkcolor=magenta,
    filecolor=magenta,
    urlcolor=magenta,
    citecolor=magenta,
}
\definecolor{caribbeangreen}{rgb}{0.0, 0.8, 0.6}

\setlist{topsep=0ex,itemsep=1ex}

    \DeclareMathOperator{\dom}{dom}

    \DeclareMathOperator{\ran}{ran}

    \newcommand{\Awf}{\mathcal{A}}
    \newcommand{\calER}{\mathcal{ER}}
    \newcommand{\Bwf}{\mathcal{B}}
    \newcommand{\Cwf}{\mathcal{C}}

    \newcommand{\Fwf}{\mathcal{F}}

    \newcommand{\Nwf}{\mathcal{N}}

    \newcommand{\Swf}{\mathcal{S}}

    \newcommand{\bfrak}{\mathfrak{b}}
    
    \newcommand{\dfrak}{\mathfrak{d}}

    \newcommand{\menos}{\smallsetminus}
    
    \DeclareMathOperator{\pts}{\mathcal{P}}
    
    \newcommand{\frestr}{{\upharpoonright}}
    \DeclareMathOperator{\add}{add}
    \DeclareMathOperator{\cov}{cov}
    \DeclareMathOperator{\non}{non}
    \DeclareMathOperator{\cof}{cof}

    \newcommand{\Ncal}{\mathcal{N}}

    \newcommand{\Z}{\mathbb{Z}}
    
    \newcommand{\R}{\mathbb{R}}

    \newcommand{\imp}{\Rightarrow}

    \newcommand{\la}{\langle}
    \newcommand{\ra}{\rangle}

\newcommand{\llex}{<_{\mathrm{lex}}}

\DeclareMathOperator{\sqw}{\mathrm{seq}_{<\omega}}
\newcommand{\Fr}{\mathrm{Fr}}

\newcommand{\leqT}{\leq_{\mathrm{T}}}
\newcommand{\eqT}{=_{\mathrm{T}}}

\newcommand{\inter}{\mathrm{int}}

\newcommand{\baire}{{}^\omega\omega}

\newcommand{\cantor}{{}^\omega2}

\newcommand{\IP}{\mathcal{IP}}
\newcommand{\TP}{\mathcal{TP}}
\newcommand{\BP}{\mathcal{BP}}
\newcommand{\GP}{\mathcal{GP}}

\newcommand{\vol}{\mathrm{Vol}}

\newcommand{\alt}{\mathrm{ht}}

\newcommand{\B}{\mathrm{B}}

\DeclareMathOperator{\E}{E}

\newcommand{\gb}{\mathfrak{b}}

\newcommand{\gd}{\mathfrak{d}}

\newcommand{\calT}{\mathcal{T}}
\newcommand{\calC}{\mathcal{C}}

\newcommand{\calA}{\mathcal{A}}

\newcommand{\calB}{\mathcal{B}}
\newcommand{\calP}{\mathcal{P}}
\newcommand{\calS}{\mathcal{S}}
\newcommand{\calN}{\mathcal{N}}

\newcommand{\calL}{\mathcal{L}}

\newcommand{\calI}{\mathcal{I}}

\newcommand{\calF}{\mathcal{F}}

\newcommand{\calY}{\mathcal{Y}}

\newcommand{\varp}{\varepsilon}

\newcommand{\rest}{{\restriction}}

\newenvironment{PROOF}[2][\proofname.]
   {\begin{proof}[#1]\renewcommand{\qedsymbol}{\textsquare$_{\mathrm{#2}}$}}
   {\end{proof}}

\newcommand{\concat}[2]{#1{}^{\frown}#2}

\renewcommand{\menos}{\smallsetminus}
\newcommand{\conj}{\wedge}

\definecolor{bluet}{rgb}{0.0, 0.4, 0.6}

\hypersetup{
    colorlinks=true,
    linkcolor=bluet,
    filecolor=bluet,
    urlcolor=bluet,
    citecolor=bluet,
}

\DeclareMathOperator{\cl}{c {\ell}}

    \newcommand{\set}[2]{\{#1 \colon #2\}}
    \newcommand{\bigset}[2]{\left\{#1 \colon #2\right\}}
    \newcommand{\seq}[2]{\la #1 \colon #2\ra}

    \DeclareMathOperator{\hgt}{\mathrm{ht}}
    \DeclareMathOperator{\Lv}{\mathsf{Lv}}
    
    \DeclareMathOperator{\suc}{\mathsf{succ}}
    
    \DeclareMathOperator{\spt}{\mathsf{spt}}
    \DeclareMathOperator{\rt}{\mathsf{rt}}
    
    \DeclareMathOperator{\Lb}{\mathsf{Lb}}

    \newcommand{\Cv}{\mathsf{Cv}}

\definecolor{sub0}{RGB}{29,32,137}
\definecolor{sub1}{RGB}{1,71,157}
\definecolor{sub2}{RGB}{1,104,183}
\definecolor{sub3}{RGB}{0,160,234}
\definecolor{sug}{RGB}{0,154,68}
\definecolor{suy}{RGB}{208,219,1}

\definecolor{redun}{rgb}{0.65, 0.11, 0.19}
\definecolor{greenun}{rgb}{0.58, 0.71, 0.23}
\definecolor{dodger}{rgb}{0.0,0.5,1.0}
\definecolor{carrotorange}{rgb}{0.93, 0.57, 0.13}

\newcommand{\red}[1]{{\color{red}#1}}

\newcommand{\Andres}[1]{{\color{orange}Andres says: #1}}
\newcommand{\Diego}[1]{{\color{blue}Diego says: #1}}

\title[Probability Trees]{Probability Trees}

\author{Diego A.~Mej\'ia}
\address{Graduate School of System Informatics, Kobe University. 1-1 Rokkodai-cho, Nada-ku, Kobe, Hyogo 657-8501 Japan}
\email{\href{mailto:damejiag@people.kobe-u.ac.jp}{damejiag@people.kobe-u.ac.jp}}
\urladdr{\url{http://www.researchgate.com/profile/Diego\_Mejia2}}

\author{Andr\'es F.\ Uribe-Zapata}
\address{TU Wien, Faculty of Mathematics and Geoinformation, Institute of Discrete Mathematics and Geometry, Wiedner Hauptstrasse 8--10, A--1040 Vienna, Austria }
\email{\href{mailto:andres.zapata@tuwien.ac.at}{andres.zapata@tuwien.ac.at}}
\urladdr{\url{https://sites.google.com/view/andres-uribe-afuz}}

\thanks{The first author was supported by the Grant-in-Aid for Scientific Research (C) 23K03198, Japan Society for the Promotion of Science. The second author was supported by the Austrian Science Fund (FWF): project number P33895.} 

\date{\today}



\begin{document}

\makeatletter
\def\@roman#1{\romannumeral #1}
\makeatother

\newcounter{enuAlph}
\renewcommand{\theenuAlph}{\Alph{enuAlph}}

\theoremstyle{plain}
  \newtheorem{theorem}{Theorem}[section]
  \newtheorem{corollary}[theorem]{Corollary}
  \newtheorem{lemma}[theorem]{Lemma}
  \newtheorem{mainlemma}[theorem]{Main Lemma}
  \newtheorem{mainproblem}[theorem]{Main Problem}
  \newtheorem{construction}[theorem]{Construction}
  \newtheorem{prop}[theorem]{Proposition}
  \newtheorem{clm}[theorem]{Claim}
  \newtheorem{fact}[theorem]{Fact}
  \newtheorem{exercise}[theorem]{Exercise}
  \newtheorem{question}[theorem]{Question}
  \newtheorem{problem}[theorem]{Problem}
  \newtheorem{cruciallem}[theorem]{Crucial Lemma}
  \newtheorem{conjecture}[theorem]{Conjecture}
  \newtheorem{assumption}[theorem]{Assumption}
  \newtheorem{teorema}[enuAlph]{Theorem}

  \newtheorem*{corolario}{Corollary}
\theoremstyle{definition}
  \newtheorem{definition}[theorem]{Definition}
  \newtheorem{example}[theorem]{Example}
  \newtheorem{remark}[theorem]{Remark}
    \newtheorem{hremark}[theorem]{Historical Remark}
    \newtheorem{observation}[theorem]{Observation}
  \newtheorem{notation}[theorem]{Notation}
  \newtheorem{context}[theorem]{Context}

  \newtheorem*{defi}{Definition}
  \newtheorem*{acknowledgements}{Acknowledgements}

\numberwithin{equation}{theorem}
\renewcommand{\theequation}{\thetheorem.\arabic{equation}}

\def\sectionautorefname{Section}
\def\subsectionautorefname{Subsection}
%\maketitle

%\tableofcontents

\begin{abstract}
   In this article, we introduce a formal definition of the concept of \emph{probability tree} and conduct a detailed and comprehensive study of its fundamental structural properties. In particular, we define what we term an \emph{inductive probability measure} and prove that such trees can be identified with these measures. Furthermore, we prove that probability trees are completely determined by probability measures on the Borel $\sigma$-algebra of the tree's body. 

   We then explore applications of probability trees in several areas of mathematics, including probability theory, measure theory, and set theory. In the first, we show that the cumulative distribution of finitely many dependent and non-identically distributed Bernouli tests is bounded by the cumulative distribution of some binomial distribution. In the second, we establish a close relationship between probability trees and the real line, showing that Borel, measurable sets, and their measures can be preserved, as well as other combinatorial properties. Finally, in set theory, we establish that the null ideal associated with suitable probability trees is Tukey equivalent to the null ideal on $[0, 1]$. This leads to a new elementary proof of the fact that the null ideal of a free $\sigma$-finite Borel measure on a Polish space is Tukey equivalent with the null ideal of $\R$, which supports that the associated cardinal characteristics remain invariant across the spaces in which they are defined.
\end{abstract}

\maketitle

% \setcounter{secnumdepth}{1}
% \setcounter{tocdepth}{1}
% %\tableofcontents

{
\small
\hypersetup{linkcolor=black}
\tableofcontents
}

%%%%%%%%%%%%%%%%%%%%%%%%%%%%%%%%%%%%%%%%%%%%%%%%%%%%%%%%%%%%%%%%%%%%%%%%%%%%%%%%%%%%%%%%%%%%%%%%%%%%%%%%%%%%%%%%%%%%%%%%%%%%%%%%%%%%%%%%%%%%%%%%%%%%%%%%%%%%%%%%%%%%%%%%%%%%

\section{Introduction}\label{SecIntro}

In probability theory and statistics, the so-called \emph{probability trees} are graphic tools that allow the representation of problems and help in their understanding and subsequent study, particularly those problems involving some kind of dependency. These trees can be used to visually analyze all the possible results of an experiment and the probabilities associated with each of them, making it easier to analyze situations involving decisions and dependent events and calculate their respective probabilities.  

Intuitively, a probability tree starts with a single node —called the \emph{root} of the tree\footnote{However, in some contexts, such as in the study of genealogy, unrooted probability trees are also considered (see e.g. \cite{GRIFFUT}).} —which represents the first event in the experiment. In the next step, from the root, branches of the tree fan out to represent all possible outcomes of the first experiment, labeled with their respective probabilities whose sum must be equal to one since they cover all possible outcomes. Each branch ends in a new node, from which new nodes and new branches are generated as the experiment is carried out, which ultimately determines the tree structure, covering all possible events in the experiment under consideration. The probabilities of each branch can be multiplied along the paths to calculate cumulative probabilities —that is, the probability that two or more events occur simultaneously— allowing for the complete determination of any sequence of events in the experiment.

In some contexts beyond probability and statistics, probability trees are relevant in many different areas of the natural sciences, such as in the field of genetics, where they are used to help model the inheritance of biological traits as well as diseases. Specifically, they can be used to determine the probability that an individual inherits a certain gene from his or her parents, connecting information between recessive and dominant inheritance patterns (see~\cite{Friedman1980}). In an area of science directly related to genetics, namely genealogy, probability trees have other very significant applications. There they can be used to model and analyze situations of relationships between generations, as well as to calculate the probabilities of transmission of certain biological traits, similar to the case of genetics mentioned above. In this context, the trees in question help to represent family connections and possible routes of genetic inheritance, especially when dealing with dependent events, such as calculating the probability that a certain individual inherits a specific gene or disease from their ancestors. A well-known example of this is the so-called \emph{infinitely-many-sites models}, which are used to study genetic variability and mutation processes in populations over time (see~\cite{1989_Griffiths}). 

As we mentioned at the beginning of this introduction, probability trees are especially useful in the analysis of dependent events, that is, those in which the probability of an event is conditioned by the results of previous events. In these cases, the probability in the successive branches reflects such dependence, making these trees an essential tool for representing problems in the context of Bayesian Inference (see e.g.~\cite{BayesInf}). This makes this type of tree very useful in some sub-areas of computer science such as artificial intelligence and machine learning, where they are used in particular to study probabilistic and learning algorithms, such as decision trees, which allow to classify data and make predictions of future outcomes based on probabilities derived from previously known data. In particular, probability trees serve as a model for a certain type of probabilistic process, known as \emph{Causal Generative Processes}, which are essential in artificial intelligence as they allow modeling data generation based on causal relationships (see~\cite{Causal1} and~\cite{causal2}). Furthermore, probability trees are widely used to carry out simulations, including the well-known \emph{Monte Carlo methods}, where they are used in the representation and subsequent analysis of random processes (see~\cite{MonteC}).

While we have seen that probability trees arise in highly practical and diverse contexts such as genetics, genealogy, artificial intelligence, machine learning, and Bayesian inference, the origin of this paper stems from a far more abstract and unexpected field: Forcing Theory, in Set Theory. Recently, based on previous work by Saharon Shelah~\cite{Sh00} and Jakob Kellner, Saharon Shelah, and Anda T\u{a}nasie~\cite{KST}, Miguel A.\ Cardona, along with the authors in~\cite{CMU}, introduced a general theory of iterated forcing using finitely additive measures\footnote{A preliminary version of the general theory of iterated forcing using finitely additive measures was presented in the master's thesis of the second author (see~\cite{uribethesis}), where the first questions related to probability trees arose. In this thesis, an entire chapter was dedicated to the formalization and study of these trees (see \cite[Ch.~3]{uribethesis}), which served as the starting point for the development of this article.}. This theory is founded on what the authors referred to as \emph{$\mu$-$\calY$-linkedness}, a property of partially ordered sets, where $\mu$ is an infinite cardinal and $\calY$ is a class of ordered pairs. A central focus of the theory was to prove that \emph{random forcing} —the poset $\calB(\cantor) \menos \calN(\cantor)$ ordered by inclusion— satisfies this property for certain suitable parameters $\mu_{0}$ and $\calY_{0}$. Later, the authors in \cite{MUR23} extended this result, a task that was not only highly technical but also required the development and application of concepts related to probability trees (as in the original work of Shelah~\cite{Sh00}). 
The core of this proof was to prove the existence of two objects: a finite set $u$ and an element of random forcing $r^{\oplus}$ satisfying certain special conditions. Without delving into overly technical and unnecessarily complex details, the strategy was based on defining a probability tree whose nodes represented partial approximations of the set $u$. Then, instead of directly attempting to find the objects $u$ and $r^{\oplus}$, following the approach of the \emph{probabilistic method} (see~\cite{PMethod} and~\cite[pp.~$100$-$101$]{uribethesis}), the probability of their existence was calculated using the tree structure. Finally, it was shown that this probability is positive, ensuring that such objects satisfying the required conditions can indeed be found.  To achieve this, it was necessary not only to formalize the notion of a probability tree but also to develop and analyze the structural properties of such trees. 

In the references cited so far, there is no concrete definition of the notion of \emph{probability tree}, as it is usually tailored to the needs of each particular case. Moreover, there is currently no detailed study of the structure of these trees. For this reason, the main objective of this article is clear and specific: to formalize the notion of \emph{probability tree} and analyze its structure rigorously. Our formalization of this concept —presented in \autoref{pp46}— is quite intuitive and, in general terms, outlined in the second paragraph of this introduction.
%: formally we define a probability tree as a pair $\langle T, \bar{\mu} \rangle$, where $T$ is a tree —in the set-theoretic sense of the term (see \autoref{3})— that satisfies certain properties, and $\bar{\mu} = \langle \mu_t \colon  t \in T \menos \max(T) \rangle$ is a sequence where each $\mu_t$ is a probability measure defined on $\calP(\suc_T(t))$. Here, $\max(T)$ denotes the set of maximal nodes of $T$, and $\suc_T(t)$ represents the set of immediate successors of $t$ in $T$. 

The formal definition and structural study of probability trees are carried out mainly in~\autoref{2.3}. Our starting point is to prove that every probability tree induces a measure on the tree that generates probability space structure at each of its fronts and levels (see~\autoref{b035} and~\ref{b035-3}). The analysis of the converse of this result motivated the introduction of the notion of \emph{inductive probability measure}, which enabled a deeper exploration of the structure of theses trees. This led to the prove that probability trees can be identified with inductive probability measures and are completely determined by probability measures on the Borel $\sigma$-algebra of $[T]$ —the body of $T$, which is the set of maximal branches (see~\autoref{3}). The structural study developed in this article focuses primarily on the definition of four collections associated with these trees: \emph{tree probability sequences} ($\TP$), \emph{inductive probability measures} ($\IP$), \emph{Borel probability measures} ($\BP$), and \emph{general probability sequences} ($\GP$), as well as the connections between them, which are ultimately reduced to the commutativity of the diagram presented in~\autoref{f50}.

Once the structure of probability trees was analyzed, we laid the groundwork for exploring their applications in different areas. This structural study not only allowed us to formalize and better understand their properties but also to establish connections with other branches of mathematics. In particular, we apply probability trees to three distinct areas: probability theory, measure theory, and set theory, particularly in the combinatorics of real numbers and invariant cardinals. Below, we provide a brief description of each of these applications.

The first application —in Probability Theory—is related to a generalization of a well-known result: by adding a finite number of independent and identically distributed random variables with a Bernoulli distribution, we obtain a random variable with a binomial distribution. However, we face the situation where these random variables with  Bernoulli distribution are dependent and even not identically distributed. 
To address this problem, we use probability trees to achieve the following result, which corresponds to~\autoref{p70}. 

\begin{teorema}\label{thmA}
    Let $p\in[0,1]$, $n$ be a natural number, and $Y$ be the random variable representing the number of successes of $n$-many dependent Bernoulli distributed random variables, where the probability of success of each variable also depends on the previous events. If $p$ is a lower bound of the probability of success of each Bernoulli-distributed random variable, then the cumulative distribution of $Y$ is below the cumulative distribution of the binomial distribution with parameters $n$ and $p$.  
\end{teorema}

This situation arose in the proof of \cite[Main Lemma~7.17]{CMU} (see also \cite[Main Lemma 4.3.17]{uribethesis}), where \autoref{thmA} proved sufficient for the purposes required in that context. 

The second application —in Measure Theory— establishes a connection between probability trees and the real line. In particular, we will prove the following theorem, which corresponds to~\autoref{5.2}. 

\begin{teorema}
   Every probability tree $\la T,\bar\mu\ra$ defines a canonical probability measure $\lambda^{\bar\mu}$ on the Borel $\sigma$-algebra $\calB([T])$ of $[T]$ which has a connection with the Lebesgue measure of the unit interval.
\end{teorema}

The measure $\lambda^{\bar{\mu}}$ is defined through a function $g_{\bar{\mu}}$ defined on  $[0, 1]$, expect on a countable subset, constructed from the probability tree $\langle T, \bar{\mu} \rangle$. The connection between $\langle T, \bar{\mu} \rangle$ and $[0, 1]$ is established through $g_{\bar{\mu}}$ as, when restricted to an appropriate set, turns out to be a homeomorphism, thus preserving the Borel sets and, consequently, the measurable sets (see~\autoref{b0100} and~\autoref{b101}). Additionally, this function preserves the measure between the measurable spaces $\langle [T], \calB([T]), \lambda^{\bar{\mu}} \rangle$ and $\langle [0, 1], \calB([0, 1]), \Lb \rangle$, where $\Lb$ denotes the Lebesgue measure. On the other hand, the properties that define $\lambda^{\bar{\mu}}$ will also allow us to demonstrate that probability trees can characterize the probability measures defined on $\calB([T])$ (see~\autoref{ppp14.T}). 

The third application —in Set Theory— is related to \emph{cardinal invariants}. \emph{Cardinal invariants}, also called   \emph{cardinal characteristics}, are cardinal numbers that capture combinatorial properties of infinite spaces. Examples of this kind of cardinals arise considering ideals. Recall that for a non-empty set $X$, $\calI \subseteq \calP(X)$ is an \emph{ideal on $X$}, if it is closed under finite unions, $\subseteq$-downwards closed, $\emptyset\in\calI$, and $X \notin \calI$. In this context, we define the \emph{cardinals invariants associated with $\calI$} as follows: 

\begin{tabular}{rr@{~}l}
      \qquad\quad   
       & $\add(\calI)$ & $\coloneqq \min \{ \vert \mathcal{F} \vert \colon \calF \subseteq \calI \conj \bigcup \calF \notin \calI \}$ is the \emph{additivity of $\calI$}, \index{$\add(\calI)$}\\[1ex]

       & $\cov(\calI)$ & $\coloneqq \min \{ \vert \mathcal{F} \vert \colon \mathcal{F} \subseteq \mathcal{I} \conj \bigcup \mathcal{F} = X \}$ is the \emph{covering of $\calI$},\\[1ex]

       & $\non(\calI)$ & $\coloneqq \min \{ \vert Y \vert \colon Y \subseteq X \conj Y \notin \calI \}$ is the \emph{uniformity of $\calI$}, \\[1ex]

       & $\cof(\calI)$ & $\coloneqq \min \{ \vert \calF \vert \colon \calF \subseteq \calI \conj \forall A\in \mathcal{I} \exists B\in \mathcal{F} \ (A \subseteq B) \}$ is the \emph{cofinality of $\calI$}.\\
\end{tabular}

Apparently, there is no unanimous reason why they are called \emph{invariants}, but it is known that they possess an invariance property: in many cases, the associated cardinal characteristics do not depend on the space on which the ideal is defined, as long as the space satisfies certain properties (see~\autoref{a243}). The connection between probability trees and the real line that we establish in~\autoref{5.2} allows us to extend this invariance property to the null ideal of $\lambda^{\bar{\mu}}$. If $\langle X, \calA, \mu \rangle$ is a measure space, $\Ncal(\mu)$ denotes the ideal of all $\mu$-measure zero subsets of $X$. When the measure space is understood, we just write $\calN(X)$, e.g.\ $\calN(\R)$ and $\calN([0,1])$ with respect to the Lebesgue measure.

\begin{teorema}
   If $\langle T, \bar{\mu} \rangle$ is a probability tree such that $\lambda^{\bar{\mu}}$ is free, the cardinal invariants associated with $\calN(\lambda^{\bar{\mu}})$ and $\calN([0, 1])$ are the same. 
\end{teorema}

Moreover, these identities follow by the \emph{Tukey-equivalence} (\autoref{c46}) between structures associated with these ideals.

This theorem leads to a new elementary proof of the more general known fact that the invariance of the cardinal invariants associated with the null ideal applies for any free $\sigma$-finite measure $\mu$ on the Borel $\sigma$-algebra of a Polish space. Details are developed in \autoref{7}.

In this work, we adopt the set-theoretic treatment of natural numbers.

\begin{notation}
    We adopt the set-theoretic treatment of natural numbers (starting from $0$): $0 = \emptyset$, and each natural number is the set of its predecessors. Formally, if $n$ is a natural number, then $n = \{ 0, 1, 2, \dots, n-1 \}$. The entire set of natural numbers is denoted by $\omega$. Furthermore, as in the context of \emph{ordinal numbers} in set theory, we typically write ``$n < \omega$'' instead of ``$n \in \omega$'' and ``$\alpha\leq \omega$'' for ``$\alpha<\omega$ or $\alpha=\omega$''.\footnote{$\omega$ is the first infinite ordinal number, represented as the limit of natural numbers.}
\end{notation}

%%%%%%%%%%%%%%%%%%%%%%%%%%%%%%%%%%%%%%%%%%%%%%%%%%%%%%%%%%%%%%%%%%%%%%%%%%%%%%%%%%%%%%%%%%%%%%%%%%%%%%%%%%%%%%%%%%%%%%%%%%%%%%%%%%%%%%%%%%%%%%%%%%%%%%%%%%%%%%%%%%%%%%%%%%%%

\section{Elementary notions of measure and probability theory}

\subsection{Review of measure theory}

We review, in this section, basic notions related to measure theory and probability spaces. 

Let $X$ be a non-empty set and $\calC \subseteq \calP(X)$. The \emph{$\sigma$-algebra generated by $\Cwf$}, that is, the smallest $\sigma$-algebra of sets over $X$ that contains $\Cwf$, is denoted by $\sigma_X(\Cwf)$. If $Z\subseteq X$, then  $\calC|_Z \coloneqq \set{C\cap Z}{C\in\Cwf}$, which is a ($\sigma$-)algebra of sets over $Z$ whenever $\Cwf$ is a ($\sigma$-)algebra of sets over $X$.  Recall that, if $f\colon X\to Y$ is a function between non-empty sets, $\Awf$ a $\sigma$-algebra over $X$ and $\Awf'$ is a $\sigma$-algebra over $Y$, then $f$ is \emph{$\Awf$-$\Awf'$-measurable} if $f^{-1}[B]\in\Awf$ for all $B\in \Awf'$. Most of the time, we work with the Borel $\sigma$-algebra of a topological space: given a topological space $X$, $\Bwf(X)$ denotes the $\sigma$-algebra generated by the collection of open subsets of $X$, which is known as the \emph{Borel $\sigma$-algebra of $X$}. Any set $B\in\Bwf(X)$ is called \emph{Borel in $X$}, and a function $f\colon X\to Y$ between topological spaces is called a \emph{Borel map} if it is $\Bwf(X)$-$\Bwf(Y)$-measurable. 

Regarding functions. Let $f\colon X\to Y$ be a function between non-empty sets. For $\Cwf\subseteq\pts(X)$ and $\Cwf'\subseteq\pts(Y)$, define $f^{\leftarrow}(\Cwf')   \coloneqq \set{f^{-1}[B]}{B\in\Cwf'} $ and $f^\to(\Cwf)  \coloneqq \set{B\subseteq Y}{f^{-1}[B]\in\Cwf}.$ Functions can be used to transfer $\sigma$-algebras.

\begin{fact}\label{b013}
   Let $f\colon X\to Y$ be a function between non-empty sets. 
   \begin{enumerate}[label=\normalfont(\alph*)]
     \item\label{b013a} If $\Awf'$ is a ($\sigma$-)algebra over $Y$ then $f^{\leftarrow}(\Awf')$ is a ($\sigma$-)algebra over $X$.
     
     \item\label{b013b} If $\Awf$ is a ($\sigma$-)algebra over $X$, then $f^\to(\Awf)$ is a ($\sigma$-)algebra over $Y$.
     
     \item\label{b013c} If $\Cwf'\subseteq\pts(Y)$ then $\sigma_X(f^{\leftarrow}(\Cwf')) = f^{\leftarrow}(\sigma_Y(\Cwf'))$.
     
     %\item If $\Cwf\subseteq \pts(X)$ then $\sigma_Y(f^\to(\Cwf)) = f^\to(\sigma_X(\Cwf))$.
   \end{enumerate}
\end{fact}

As a consequence, we have the following results. 

\begin{corollary}\label{b015}
  Let $X$ be a non-empty set, $\Cwf\subseteq\pts(X)$ and $Z\subseteq X$. Then we have that $\sigma_Z(\Cwf|_Z) = \sigma_X(\Cwf)|_Z$.
\end{corollary}

\begin{corollary}\label{b016}
   If $X$ is a topological space and $Z\subseteq X$ is a subspace, then we have that $\Bwf(Z) = \Bwf(X)|_Z$. 
\end{corollary}

\begin{corollary}\label{b017}
  Any continuous function between topological spaces is Borel.
\end{corollary}

We now review the notion of \emph{measure}. Let $X$ be a non-empty set and $\Cwf\subseteq\pts(X)$ with $\emptyset\in\Cwf$. Recall that a function $\mu\colon \Cwf\to [0,\infty]$ is a \emph{finitely additive measure (fam)}, if $\mu(\emptyset)=0$ and $\mu\left(\bigcup_{k<n}A_k\right) = \sum_{k<n}\mu(A_k)$ whenever $\set{A_k}{k<n}\subseteq\Cwf$ is a finite collection of pairwise disjoint sets whose union is in $\calC$. Also, a fam $\mu$ is a \emph{($\sigma$-additive) measure} if $\mu\left(\bigcup_{n<\omega}A_n\right) = \sum_{n<\omega}\mu(A_n)$ for any collection $\set{A_n}{n<\omega}\subseteq\Cwf$ of pairwise disjoint sets whose union is in $\calC$. 

Given a fam $\mu\colon \Cwf\to [0,\infty]$, we say that $\mu$ is \emph{finite} if $X\in \Cwf$ and $\mu(X)<\infty$. When there is some $\set{A_n}{n<\omega}\subseteq \Cwf$ such that $X=\bigcup_{n<\omega}A_n$ and $\mu(A_n)<\infty$ for all $n<\omega$, $\mu$ is called \emph{$\sigma$-finite}. If   $X\in\Cwf$ and $\mu(X)=1$, then  $\mu$ is a \emph{probability fam}. Finally, $\mu$ is \emph{free} if, for any $x\in X$, $\{x\}\in \Cwf$ and $\mu(\{x\}) = 0$.

Recall that a \emph{measure space} is a triple $\la X,\Awf,\mu\ra$ where $\Awf$ is a $\sigma$-algebra over $X$ and $\mu$ is a measure on $\Awf$.

\begin{example}\label{b022}
   Let $W$ be a countable set. For any function $f\colon W\to [0,\infty]$ there is a unique measure $\mu^f\colon \pts(W)\to [0,\infty]$ such that $\mu^f(\{k\}) = f(k)$ for all $k\in W$. Indeed, for $A\subseteq W$, $\mu^f(A)$ must be $\sum_{k\in A}f(k)$. 
   
   Conversely, for any measure $\mu\colon\pts(W)\to[0,\infty]$ there is a unique function $f\colon W\to [0,\infty]$ such that $\mu=\mu^f$ ($f(k)$ must be $\mu(\{k\})$).
\end{example}

As a consequence, if $W$ is a countable or finite set, to define a probability measure on $\calP(W)$ it is sufficient to define a function $f \colon W \to [0, 1]$ such that $\sum_{w \in W} f(w) = 1$. For example, we can use this to introduce the \emph{uniform measure} on finite sets. 

\begin{definition}\label{b024.0}
    Let $X$ be a non-empty finite set. The \emph{uniform measure on $X$}, denoted by $\mu_{X}$, is the measure on $\calP(X)$ determined by $\mu_{X}(\{x \}) \coloneqq \frac{1}{\vert X \vert}$.   
\end{definition}

\begin{definition}\label{b024}
  The \emph{standard measure on $\omega$} is the (unique) probability measure on $\pts(\omega)$ obtained from the function $\omega\to[0,1]$, $k\mapsto 2^{-(k+1)}$ as in \autoref{b022}. 
\end{definition}

Measure zero sets play an essential role in measure theory. Let $\Awf\subseteq\pts(X)$ be a $\sigma$-algebra over the set $X$ and let $\mu\colon \Awf\to [0,\infty]$ be a measure. Define
   \begin{align*}
       \Nwf(\mu) &  \coloneqq  \set{N\subseteq X}{\exists\, A\in\Awf\colon N\subseteq A \text{ and }\mu(A) = 0},\\
       \Awf^\mu &  \coloneqq  \set{A\cup N}{A\in\Awf,\ N\in\Nwf(\mu)}.
   \end{align*}
   The sets in $\Nwf(\mu)$ are called \emph{$\mu$-null}, or just \emph{null}, or \emph{measure zero sets}.
   
   It is known that $\Awf^\mu$ is the $\sigma$-algebra generated by $\Awf\cup\Nwf$ and that there is a unique measure on $\Awf^\mu$ that extends $\mu$, namely, $A\cup N\mapsto \mu(A)$ for $A\in\Awf$ and $N\in\Nwf$. This measure is called the \emph{completion of $\mu$}, and we use the same letter $\mu$ to denote this completion. In terms of measure spaces, we say that $\la X,\Awf^\mu,\mu\ra$ is the \emph{completion of $\la X,\Awf,\mu\ra$}.
   
   Notice that $\Awf^\mu = \Awf$ iff $\Nwf(\mu)\subseteq \Awf$. In this situation, we say that $\mu$ is a \emph{complete measure} and $\la X,\Awf,\mu\ra$ is a \emph{complete measure space}. 

\begin{lemma}\label{b028}
    Let $f\colon X\to Y$ be a map between non-empty sets. If $\la X,\Awf,\mu\ra$ is a measure space, then $\la Y,f^{\to}(\Awf),\mu'\ra$ is a measure space where $\mu'(B) \coloneqq  \mu(f^{-1}[B])$ for $B\in f^{\to}(\Awf)$.
\end{lemma}

\begin{lemma}\label{b029}
  Let $\Awf$ be an algebra on a set $X$ and let $\mu$ be a fam on $\Awf$. Assume that $\mu(X)<\infty$ and $\{x\}\in\Awf$ for all $x\in X$. Then the set $\set{x\in X}{\mu(\{x\})>0}$ is countable.
\end{lemma}

Finally, recall:

\begin{theorem}[{\cite[\S13]{Halmos}}]\label{b070}
  Let $\Fwf$ be an algebra of sets over a set $X$ and assume that $\mu\colon \Fwf\to [0,\infty]$ is a $\sigma$-finite measure. Then, there is a unique measure on $\sigma_{X}(\Fwf)$ that extends $\mu$.
\end{theorem}

\section{Elementary notions of trees}\label{sec:tree}\label{3}

\subsection{Trees}\label{trees}

In this section, we introduce the set-theoretic notion of \emph{tree} (of height ${\leq}\omega$) and study some of its combinatorial and topological properties. 

We start by fixing some notation about functions and sequences. 

\begin{notation}\label{it:seq}
  We write $\seq{a_i}{i\in I}$ to denote a function such that $i\mapsto a_i$. We usually say that $\la a_i:\, i\in I\ra$ is a \emph{sequence} (indexed by $I$).
  
  We typically look at \emph{sequences of length ${\leq}\,\omega$}: for $n < \omega$, $\seq{a_i}{i<n}=\la a_0,\ldots,a_{n-1}\ra$ is a finite sequence of \emph{length} $n$, while $\seq{a_i}{i< \omega} = \la a_0,a_1,\ldots\ra$ is a sequence of \emph{length} $\omega$. The \emph{empty sequence} $\la\ \ra$ is the only sequence of length $0$. We use $|s|$ to denote the length of a sequence $s$ of length ${\leq}\,\omega$. When $s$ is a finite sequence and $t$ is a sequence of length ${\leq}\,\omega$, define \emph{the concatenation of $s$ and $t$} by $\concat{s}{t} \coloneqq \la s_0,\ldots, s_{|s|-1},t_0,t_1,\ldots\ra$. If $s$ and $t$ are sequences of length ${\leq}\,\omega$, then
  $s\subseteq t \text{ iff }|s|\leq |t| \text{ and }\forall \ i<|s| \ ( s_i=t_i),$ that is,  $s\subseteq t$ means that $t$ is a longer sequence extending $s$.
\end{notation}

Let $W$ be a set and $\alpha\leq\omega$. We define ${}^\alpha W$ as the set of sequences in $W$ of length $\alpha$; ${}^{< \alpha} W$ as the set of sequences in $W$ of length ${<}\,\alpha$; and  ${}^{\leq  \alpha} W$ as the set of sequences in $W$ of length ${\leq}\,\alpha$. Equivalently,  ${}^{<\alpha}W  \coloneqq \bigcup_{k<\alpha}{}^k W $ and ${}^{\leq \alpha}W   \coloneqq \bigcup_{k\leq \alpha}{}^k W$. 

Recall the notions of partial and linear orders: a \emph{partial order} is a pair $\la P,\leq \ra$ where $P$ is a non-empty set and $\leq$ is a  reflexive, transitive, and anti-symmetric relation. We say that $x$ and $y$ are \emph{comparable in $P$} if either $x\leq y$ or $y\leq x$. A \emph{chain in $P$}\index{Poset (also preorder)!chain} is a subset of $P$ such that any pair of elements are comparable. When $P$ is a chain in itself, we say that $\la P,\leq\ra$ is a \emph{linear order}. 

We are ready to introduce the notion of \emph{tree}. 

\begin{definition}\label{a10}
    A \emph{tree} (of height ${\leq}\,\omega$) is a partial order $\langle T, <_{T} \rangle$ containing a minimum element $\rt(T)$, called \emph{the root of $T$}, such that, for any $t\in T$, $t{\downarrow} \coloneqq \set{s\in T}{s<_{t}t}$ is a finite linear order (under the order of $T$). The members of $T$ are usually called the \emph{nodes of $T$}. 
\end{definition}

For example, when $W$ is a non-empty set, $\la{}^{<\omega}W,\subseteq\ra$ is a tree with  $\rt({}^{<\omega}W) = \la\ \ra$.

Now, we introduce some notation related to trees and their properties. 

\begin{definition}\label{a20}
Given a tree $\la T,\leq\ra$, we fix the following notation for $t\in T$ and $n<\omega$:

\begin{enumerate}
    \item $\hgt_T(t) \coloneqq  |t{\downarrow}|$ the \emph{height of $t$ in $T$}. 

    \item $\Lv_n(T) \coloneqq  \set{t\in T}{\hgt_T(t) = n}$ the \emph{$n$-th level of $T$}, so $\Lv_0(T) = \{\rt(T)\}$.

    \item $\hgt(T) \coloneqq  \min(\set{n<\omega}{\Lv_n(T) =\emptyset}\cup\{\omega\})$ the \emph{height of $T$}. 

     \item $\suc_T(t) \coloneqq  \set{t'\in T}{t\leq t'\text{ and }\hgt_T(t') = \hgt_T(t)+1}$ the \emph{set of immediate successors of $t$ (in $T$)}. 

     \item $\spt(T) \coloneqq  \set{t\in T}{|\suc_T(t)|\geq 2}$ the \emph{set of splitting nodes of $T$}. 

    % \item $T\wedge t \coloneqq  \set{s\in T}{s \text{ is comparable with $t$ in $T$}}$.    

    \item $T_{\geq t} \coloneqq \{ s \in T \colon t \subseteq s \}$ is the set of \emph{successors of $t$} in $T$.  

    \item $\max(T) \coloneqq \{ s \in T \colon \suc_{s}(T) = \emptyset \},$ that is, the set of \emph{maximal nodes} of $T$. 
\end{enumerate}
\end{definition}

\begin{example}\label{a21}
   Consider the tree $\la{}^{<\omega}W,\subseteq\ra$. Then, for $t\in {}^{<\omega}W$ and $n<\omega$:

    \begin{multicols}{2}
        \begin{enumerate}[label = \normalfont (\arabic*)]
            \item $\hgt(t)=|t\downarrow|  =|t|$, the length of $t$.
    
            \item $\Lv_n({}^{<\omega}W) = {}^n W$. 
    
            \item $\hgt({}^{<\omega}W) =\omega$. 
    
            \item $\suc(t) = \set{\concat{t}{\la \ell\ra}}{\ell\in W}$. 
        \end{enumerate}
    \end{multicols}
\end{example}

We introduce the following notions related to trees.

\begin{definition}\label{a30}
Let $\la T,\leq\ra$ be a tree.
\begin{enumerate}[label = \normalfont (\arabic*)]
    \item\label{a30.1} We say that $T'$ is a \emph{subtree of $T$} if $T'\subseteq T$ and, for any $t\in T'$ and $s\leq t$ in $T$, $s\in T'$.
    
    \item A tree $T$ is \emph{well-pruned} if, for any $t\in T$ and $\hgt(t)<n<\hgt(T)$, there is some $t'\in\Lv_n(T)$ above $t$.
    
    \item A tree $T$ is \emph{finitely branching} if $\suc_T(t)$ is finite for all $t\in T$.
    
    \item A tree $T$ is \emph{perfect} if, for every $t\in T$, there is some splitting node in $T$ above $t$.
\end{enumerate}
\end{definition}

Notice that, if $\alt(T) < \omega$ then $\max(T) = \Lv_{\alt(T)}(T)$ whenever $T$ is a well-pruned tree. 

On the other hand, it is not hard to check that, if $T$ is a tree and $T'$ is a subtree of $T$ then, for any $t\in T'$ and $n<\omega$, $\rt(T') = \rt(T)$,  $\hgt_{T'}(t)  = \hgt_T(t)$, $\Lv_n(T') = T'\cap\Lv_n(T)$,  $\hgt(T') \leq\hgt(T)$, and $\suc_{T'}(t) = T'\cap\suc_T(t)$. 

\begin{example}\label{a33}
    If $T$ is a subtree of ${}^{<\omega}W$ then $\rt(T) = \la\ \ra$ and $\Lv_n(T) = T\cap {}^n W$ for $n<\omega$. For any $A\subseteq W$ and $0<\alpha\leq\omega$, ${}^{<\alpha} A$ is a subtree of ${}^{<\omega}W$ of height $\alpha$. Note that ${}^{<\alpha}A$ is well-pruned, it is finitely branching iff $A$ is finite, and it is perfect iff $\alpha=\omega$ and $|A|\geq 2$.
\end{example}

\begin{theorem}\label{a40}
   Any tree of height ${\leq}\,\omega$ is isomorphic with a subtree of ${}^{<\omega}W$ for some non-empty set $W$.
\end{theorem}
\begin{PROOF}{\ref{a40}}
    Let $T$ be a tree. Put $W \coloneqq T$ (as a set). Define $f\colon T\to {}^{<\omega}W$ in the following way: for $t\in T$, enumerate $t \cup\{t\} = \set{t_i}{i\leq\hgt_T(t)}$ such that $i<j\imp t_i<t_j$ (so $t_{\hgt(t)} = t$), so set $f(t) \coloneqq \seq{t_{i+1}}{i<\hgt_T(t)}$. Then, for any $s,t\in T$, $s\leq t \text{ iff }f(s)\subseteq f(t)$ and $f[T]$ is a subtree of ${}^{<\omega}T$.
\end{PROOF}

\begin{notation}\label{seqtree}
Due to the previous theorem, from now on \textbf{all trees are trees of sequences}, i.e.\ subtrees of ${}^{<\omega}W$ for some non-empty set $W$, unless otherwise stated.
\end{notation}

From the next sections, the space of infinite branches of a tree will be very important to understand the combinatorics and topology of the reals.

\begin{definition}\label{a50}
    When $T$ is a subtree of ${}^{<\omega}W$, we define 
    \begin{align*}
     \lim T & \coloneqq \set{x\in {}^\omega W}{x\frestr n\in T \text{ for all $n<\omega$}};\\
     [T] & \coloneqq \lim T\cup \max T.
    \end{align*}
    Notice that $[T]=\max T$ when $T$ has finite height (because, in this case, $\lim T=\emptyset$). Since there is a bijection between $[T]$ and the maximal chains contained in $T$, we call $[T]$ the
    \emph{space of maximal branches of the tree} or \emph{the body of $T$}. 
    %This notation is mostly used when $T$ is well-pruned. 
    %When $T$ has infinite height, 
    Any $x\in \lim T$ is called an \emph{infinite branch of $T$}.
\end{definition}

It is clear that $[{}^{<\omega}W ] = \lim  {}^{<\omega}W = {}^\omega W$.

If $T$ is a subtree of ${}^{< \omega} W$, then $\lim T \neq \emptyset$ implies $\hgt(T) = \omega$. However, the converse is not always true, since there could be trees of height $\omega$ with $\lim T = \emptyset$. For example, for $n<\omega$ let $t_n  \coloneqq  \concat{\la n \ra}{\bar 0^n}$ where $\bar 0^n$ is the sequence of length $n$ composed with only $0$. Then $S \coloneqq \set{s\in {}^{<\omega}\omega}{\exists\, n<\omega\colon s\subseteq t_n}$ is a counterexample. The so-called \emph{K\"onig's Theorem} (see~\autoref{Konig}) gives us the sufficient conditions to have the equivalence. 

\begin{theorem}\label{Konig}
    Let $T$ be a subtree of ${}^{<\omega}W$. If\/ $T$ is finitely-branching, then the following statements are equivalent.
    
    \begin{multicols}{3}
        \begin{enumerate}[label=\normalfont(\roman*)]
          \item\label{k1} $\lim T\neq\emptyset$.
          \item\label{k2} $\hgt(T) = \omega$.
          \item\label{k3} $T$ is infinite. 
\end{enumerate}
    \end{multicols}
\end{theorem}

The set $\lim T$ is also related to the notion of well-foundedness.

\begin{definition}\label{a054}
    A tree $T$ is \emph{well-founded} if every non-empty subset of $T$ has a maximal element.
\end{definition}

\begin{lemma}\label{a055}
    Let $T$ be a tree of sequences. Then, the following statements are equivalent.
    \begin{multicols}{2}
    \begin{enumerate}[label = \normalfont (\roman*)]
        \item $T$ is well-founded. 
        \item $\lim T =\emptyset$.
    \end{enumerate}
    \end{multicols}
\end{lemma}
\begin{PROOF}{\ref{a055}}
    If $\lim T \neq \emptyset$ and contains some infinite branch $x$, then $\set{x\frestr n}{n<\omega}$ is a subset of $T$ without maximal elements. Conversely, if $T$ is not well-founded, i.e.\ there is some non-empty $A\subseteq T$ without maximal elements, then we can construct an increasing sequence $\seq{t_n}{n<\omega}$ in $A$. This increasing sequence determines a unique member of $\lim T$, so $\lim T \neq\emptyset$.
\end{PROOF}

\subsection{Tree-topology}\label{sec:limT}

In this Subsection, we assign a topology to $[T]$ for any tree $T$ and study some of its properties. We start by recalling some basic topological notions. 

Consider a topological space $\langle X, \tau \rangle$. 
%Recall that a collection $\Bwf\subseteq \pts(X)$ is a \emph{base of $\tau$} if $\Bwf\subseteq\tau$ and, for any $A\subseteq X$, if $A\in\tau$ then for any $x \in A$, there exists some $B \in \calB$ such that $x \in B \subseteq A$. 
A set $C\subseteq X$ is \emph{clopen in $\la X,\tau\ra$} if it is open and closed in $\la X,\tau\ra$. We say that $\langle X, \tau \rangle$ is \emph{zero-dimensional} if it has a base of clopen sets. For $A\subseteq X$, $\cl_X(A)$ denotes the \emph{closure of $A$}, and $\inter_X(A)$ is the \emph{(topological) interior of $A$}. The subindex is removed when clear from the context. Recall that a topological space $X$ is \emph{discrete} if every subset of $X$ is open.

If $\Swf\subseteq \pts(X)$, the smallest topology of $X$ containing $\Swf$, which is called the \emph{topology of $X$ generated by $\Swf$}, is denoted by $\tau_\Swf$.  

We define the topology of the branches of a tree.

\begin{definition}\label{def:toptree}
   Let $T$ be a tree. The \emph{tree-topology of\/ $[T]$} is the topology generated by $\set{[t]_T}{t\in T}$, where $[t]_T \coloneqq \set{x\in [T]}{t\subseteq x}$. We just write $[t]$ when $T$ is clear from the context. Denote $\Bwf_T\coloneqq \Bwf([T])$.
\end{definition}

Notice that, if $T$ is a tree of finite height, then the tree topology is the discrete topology. More generally, $\max T$ (in case it is non-empty) is a discrete subspace.

%From now on, fix a well-pruned subtree $T$ of ${}^{<\omega} W$ of height $\omega$. 
Assume that $T$ is a tree of height ${\leq}\, \omega$. 
We say that $s,t\in T$ are \emph{compatible (in $T$)} if either $s\subseteq t$ or $t\subseteq s$. Otherwise, they are \emph{incompatible}, which we represent by $s\perp_T t$, or just $s\perp T$. It is not hard to check the following. 

\begin{fact}\label{fct:clp}
    Let $T$ be a tree and $s,t\in T$. Then:
    \begin{enumerate}[label = \normalfont (\alph*)]
        \item If $s\subseteq t$ then $[t]\subseteq [s]$.
        \item $s\perp t$ iff $[s]\cap [t] = \emptyset$.
        \item\label{clp-c} $[t]\subseteq [s]$ iff either $s\subseteq t$, or $t\subseteq s$ and there are no splitting nodes between $t$, including it, and $s$, excluding it. (The latter implies $[t]=[s]$.)
    \end{enumerate}
\end{fact}

\begin{lemma}\label{fct:Tbase}
  The collection $\set{[t]_T}{t\in T}$ is a base of the topology of\/ $[T]$ and each $[t]_T$ is clopen in $[T]$. In particular, $[T]$ is a zero-dimensional space.
\end{lemma}

The countable product of discrete spaces can be expressed as a topological space of the form $[T]$. 

\begin{definition}\label{def:prodb}
   Let $b = \seq{b(n)}{n<\omega}$ be a sequence of non-empty sets. Define $\sqw b \coloneqq  \bigcup_{n<\omega}\prod_{k<n}b(k)$, which is a well-pruned tree of height $\omega$ and $\Lv_n(\sqw b)= \prod_{k<n}b(k)$ for all $n<\omega$. Notice that, 
   $$[\sqw b] = \prod b \coloneqq \prod_{n<\omega}b(n) = \set{x}{x \text{ is a function, }\dom x =\omega,\ \forall\, n<\omega\colon x(n)\in b(n)}.$$
   In the case that $b(n)=W$ for all $n<\omega$, $\sqw b = {}^{<\omega} W$ and $\prod b = {}^\omega W$.

   Two very important spaces are defined in this way: the \emph{Cantor Space} ${}^\omega 2 = {}^\omega \{0,1\},$ and the \emph{Baire space }  ${}^\omega \omega$.
\end{definition}

\begin{lemma}\label{lem:prodb}
  Let $b = \seq{b(n)}{n<\omega}$ be a sequence of discrete spaces. Then the tree-topology of\/ $\prod b$ is the same as the product topology.
\end{lemma} 

\begin{PROOF}{\ref{lem:prodb}}
    It is enough to show that $\set{[t]}{t\in\sqw b}$ is a base of the product topology. First, for $[t]\in\sqw b$, $[t]$ is open in the product topology because $[t]=\prod_{i<\omega}a_i$ where
    \[a_i \coloneqq 
      \begin{cases}
         \{t_i\} & \text{ if $i<|t|$,}\\
         b(i) & \text{ if $i\geq|t|$.}
      \end{cases}\]
    
    Now assume that $A$ is open in the product topology and $x\in A$. Then, there is some sequence $\la c_i\colon i<\omega\ra$ such that each $c_i\subseteq b(i)$, $\set{i<\omega}{c_i\neq b(i)}$ is finite, and $x\in \prod_{i<\omega}c_i \subseteq A$. Hence, there is some $m<\omega$ such that $\set{i<\omega}{c_i\neq b(i)}\subseteq m$. Therefore, $x\in[x\frestr m]\subseteq \prod_{i<\omega}c_i \subseteq A$.
\end{PROOF}

The Cantor space is compact as a consequence of the following result, which follows by \autoref{Konig}.

\begin{theorem}\label{Tcpt}
    Let $T$ be a tree. Then $[T]$ is compact iff $T$ is finitely branching. 
\end{theorem}

%%%%%%%%%%%%%%%%%%%%%%%%%%%%%%%%%%%%%%%%%%%%%%%%%%%%%%%%%%%%%%%%%%%%%%%%%%%%%%%%%%%%%%%%%%%%%%%%%%%%%%%%%%%%%%%%%%

\section{Probability trees}\label{2.3}

In this section, we formalize the notion of a \emph{probability tree} and explore some of its structural properties. This will allow us to prove that such trees can be identified with a specific class of sequences, which we call \emph{inductive probability measures}, as well as with probability measures on the Borel $\sigma$-algebra on the space of maximal branches of the tree (see~\autoref{pp46},~\autoref{pp160}, and~\autoref{ppp14.T}). Later, in~\autoref{4.2}, we will build on the concept of conditional probability to define a relative expected value within this framework.

\begin{notation}\label{b001}
    We denote by $\calT$ the collection of all countable trees of sequences.
\end{notation}

\begin{remark}\label{b001rem}
Although, for simplicity, we develop all the theory in this section for countable trees of sequences, thanks to \autoref{a40} it can be applied, in a natural way, to arbitrary countable trees, regardless of whether they are composed by sequences.
\end{remark}

\subsection{Trees as probability spaces}\ 

We show how to define probability spaces from a tree in the sense of \autoref{3}.  

\begin{definition}\label{b030}
   We say that $\la T,\bar\mu\ra$ is a \emph{probability tree} if $T\in\calT$  
   and $\bar\mu = \seq{\mu_{t}}{t\in T \menos \max(T)}$, where each $\mu_{t}$ is a probability measure on $\pts(\suc_T(t))$. 

   Furthermore, we define $\TP$ as the class of all sequences $\bar{\mu}$ such that $\langle T, \bar{\mu} \rangle$ is a probability tree for some $T \in \calT$. Notice that this $T$ is uniquely determined by (the domain of) $\bar\mu$, so it will be denoted by $T_{\bar{\mu}}$.
\end{definition}

The following are examples of probability trees. 

\begin{example}\label{b036}
    \
    \begin{enumerate}[label=\normalfont(\arabic*)]
    \item\label{b036-1} $T \coloneqq {}^{<\omega}2$ and, for any $t\in T$, $\mu_t$ is the uniform measure on $\suc_T(t)$, that is, $\mu_t(\{\concat{t}{\la 0\ra} \}) = \mu_t(\{\concat{t}{\la 1\ra} \}) \coloneqq \frac{1}{2}$ (see~\autoref{b024.0}). 
    
    \item\label{b036-2} $T \coloneqq {}^{<\omega}\omega$ and, for $t\in T$, $\mu_t$ resembles the standard measure on $\omega$, that is, $\mu_t(\{\concat{t}{\la\ell\ra}\}) \coloneqq 2^{-(\ell+1)}$ for $\ell<\omega$ (see \autoref{b024}).
    
    \item\label{b036-3} $T \coloneqq {}^{<\omega}\omega$ and, for $t\in T$ and $\ell<\omega$,
    \[\mu_t(\concat{t}{\la\ell\ra})  \coloneqq 
    \begin{cases}
       1 & \text{ if $\ell =5$,}\\
       0 & \text{ if $\ell\neq 5$.}
    \end{cases}
    \]
    \end{enumerate}
\end{example}

Order-isomorphisms preserve the probability tree structure. 

\begin{lemma}\label{pp1}
    Let $\langle T, \bar{\mu} \rangle$ be a probability tree, $S$ a partial order, and $\varphi \colon T \to S$ an order-isomorphism. Then $\langle S, \bar{\tau} \rangle$ is a probability tree, where $\bar{\tau} = \langle \tau_{s} \colon s \in S \menos \max(S) \rangle$, and for any $s \in S$ and $r \in \suc_{S}(s)$, $\tau_{s}(\{ r \})  \coloneqq  \mu_{\varphi^{-1}(s)}(\{ \varphi^{1-}(r) \})$. 
\end{lemma}

\begin{figure}[H]
        \centering
        \begin{tikzpicture}[scale=1]

            % nodos
            \node at (4, 4) {$\bullet$}; 
            \node at (2, 2) {$\bullet$}; 
            \node at (6, 2) {$\bullet$}; 
            \node at (1, 0) {$\bullet$}; 
            \node at (3, 0) {$\bullet$}; 
            \node at (5, 0) {$\bullet$}; 
            \node at (7, 0) {$\bullet$}; 

            %aristas 
            \draw (4, 4) -- (2, 2);
            \draw (4, 4) -- (6, 2);

            \draw (2, 2) -- (1, 0);
            \draw (2, 2) -- (3, 0);

            \draw (6, 2) -- (5, 0);
            \draw (6, 2) -- (7, 0);

             \node at (4, 4.45) {$\langle \, \rangle$}; 

            \node at (1.6, 2) {$w_{1}^{0}$};
            \node at (6.4, 2) {$w_{1}^{1}$};

            \node at (1, -0.4) {$w_{2}^{0}$};
            \node at (3, -0.4) {$w_{2}^{1}$};
            \node at (5, -0.4) {$w_{2}^{2}$};
            \node at (7, -0.4) {$w_{2}^{3}$};

            % probabilidades 

            \node[redun] at (2.6, 3.1) {$p_{1}^{0}$};
            \node[redun] at (5.4, 3.1) {$p_{1}^{1}$};

            \node[redun] at (1.1, 1.1) {$p_{2}^{0}$};
            \node[redun] at (2.9, 1.1) {$p_{2}^{1}$};
            \node[redun] at (5.1, 1.1) {$p_{2}^{2}$};
            \node[redun] at (6.9, 1.1) {$p_{2}^{3}$};
        \end{tikzpicture}
        
        \caption{Example of a probability tree.}
        \label{f36}
    \end{figure}
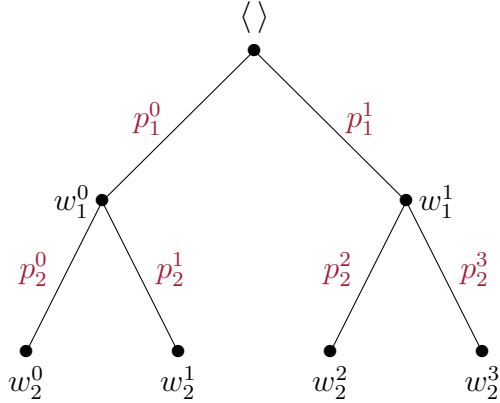

Let $T$ be the subtree of ${}^{< \omega} W$ whose set of nodes is $\{ w_{i}^{j} \colon i < 2 \conj j < 4 \}$ (see \autoref{f36}). 

If we define $\mu_{\langle \, \rangle}( \langle w_{0}^{i} \rangle ) \coloneqq p_{0}^{i},$ for $i \in \{0, 1 \}$; $\mu_{\langle w_{1}^{0} \rangle}(\langle w_{1}^{0}, \, w_{2}^{j} \rangle) \coloneqq p_{2}^{j},$ for $j \in \{0, 1\}$; and $\mu_{\langle w_{1}^{1} \rangle}(\langle w_{1}^{1}, \, w_{2}^{j} \rangle) \coloneqq p_{2}^{j},$ for $j \in \{2, 3\},$ then $T$ is a probability tree if, and only if: 
$$ p_{1}^{0} + p_{1}^{1} = 1, \ p_{2}^{0} + p_{2}^{1} = 1 \ \text{ and } \ p_{2}^{2} + p_{2}^{3} = 1. $$

Notice that, in that case, it satisfies the following: 
$$p_{1}^{0} p_{2}^{0} + p_{1}^{0}  p_{2}^{1} + p_{1}^{1}  p_{2}^{2} + p_{1}^{1}  p_{2}^{3} = p_{1}^{0}  (p_{2}^{0} + p_{2}^{1}) + p_{1}^{1}  (p_{2}^{2}+ p_{2}^{3}) = p_{1}^{0} + p_{1}^{1} = 1,$$ that is, if for any $t = \langle z_{j}^{k}, z_{l}^{m} \rangle \in \Lv_{2}(T)$ we define $\mu_{2}^{T}(t) \coloneqq p_{j}^{k} p_{l}^{m},$ then we have that $\langle \Lv_{2}(T),  \calP(\Lv_{2}(T)), \mu_{2}^{T} \rangle$ is a probability space. The same  happens trivially at $\Lv_{1}(T)$. 
    
As a consequence, %with a suitable definition of [$\mu_{n}^{2}$ for $n < 3$, 
it is possible to induce a probability space structure on each level of $T$, and even a measure on the whole tree. To formalize this idea, we introduce the following definition. 
    
\begin{definition}\label{b036.1}
    Let $ \langle T, \bar{\mu} \rangle$ be a probability tree. 
    Define the measure $\Xi^{\bar\mu}\colon \pts(T)\to[0,1]$ determined by
    \[\Xi^{\bar\mu}(\{t\})\coloneqq \prod_{k<n}\mu_{t\frestr k}( \{ t\frestr(k+1) \}) \text{ for $t\in T$.}\]
    For any $n < \hgt(T)$, $\Xi^{\bar\mu}_{n}\coloneqq \Xi^{\bar\mu}\frestr\pts(\Lv_n(T))$, which is a measure on $\pts(\Lv_n(T))$.
\end{definition}

In some cases, the probability of a successor's space in a probability tree can be determined using $\Xi^{\bar{\mu}}$. 

\begin{lemma}\label{pp13}
    Let $\langle T, \bar{\mu} \rangle$ be a probability tree, $t \in T \menos \max(T)$ and $s \in \suc_{T}(t)$. Then $ \Xi^{\bar{\mu}}(\{ s \}) = \mu_{t}(\{s \}) \Xi^{\bar{\mu}}(\{ t \})$. In particular, if $\Xi^{\bar{\mu}}(\{  t \}) > 0$, then 
    $\mu_{t}(\{s \}) = \frac{\Xi^{\bar{\mu}}(\{ s \})}{\Xi^{\bar{\mu}}(\{ t \})}.$
\end{lemma}

\begin{PROOF}{\ref{pp13}}
      Let $t \in T \menos \max(T)$ and $s \in \suc_{T}(t)$. Then, 
      \begin{equation*}
        \begin{split}
            \Xi^{\bar{\mu}}(\{ s \}) & = \prod_{k < \vert t \vert +1} \mu_{s \rest (k+1)} (\{ s \rest ( k + 1) \}) = \mu_{s \rest \vert t \vert}(\{ s \rest ( \vert t \vert + 1) \})\prod_{k < \vert t \vert} \mu_{s \rest k} (\{ s \rest ( k +1) \}) \\ 
            &= \Xi^{\bar{\mu}}(\{ s \rest \vert t \vert \}) \mu_{s \rest \vert t \vert}(\{ s \rest (\vert t \vert + 1) \}) = \Xi^{\bar{\mu}}( \{ t \}) \mu_{t}(\{ s \}).
        \end{split}
    \end{equation*}
    Thus, $\Xi^{\bar{\mu}}(\{ s \}) = \Xi^{\bar{\mu}}( \{ t \}) \mu_{t}(\{ s \}).$ 
\end{PROOF}

The induced measure $\Xi^{\bar\mu}$ satisfies the following properties.

\begin{lemma}\label{pp131}
    If $\la T,\bar\mu\ra$ is a probability tree, then:
    \begin{enumerate}[label = \normalfont (\alph*)]
        \item\label{pp131a} $\Xi^{\bar\mu}(\{\la\, \ra\}) =1$.
        \item\label{pp131b} $\Xi^{\bar\mu}(\{t\})=\Xi^{\bar\mu}(\suc_T(t))$ for any $t \in T \menos \max(T)$.  
    \end{enumerate}
\end{lemma}
\begin{PROOF}{\ref{pp131}}
    \ref{pp131a}: $\displaystyle\Xi^{\bar\mu}( \{ \langle \, \rangle \} ) =  \prod_{i < 0} \mu_{t \rest i}( \{ t \rest  (i+1) \} ) = 1$, where the last equality holds because empty products are equal to 1.

    \ref{pp131b}: By~\autoref{pp13}, we have that: 
    \begin{equation*}
        \begin{split}
            \Xi^{\bar{\mu}}(\suc_{T}(t)) & =  \sum_{s \in \suc_{T}(t)} \Xi^{\bar{\mu}}(\{ t \}) \mu_{t}(\{ s \}) = \Xi^{\bar{\mu}}(\{ t \}) \sum_{s \in \suc_{T}(t)}  \mu_{t}(\{ s \}) = \Xi^{\bar{\mu}}(\{ t \}). \qquad 
        \end{split}
        \qedhere
    \end{equation*}
\end{PROOF}

Considering the probability trees from~\autoref{b036}, we can calculate the corresponding measures on its levels.  

\begin{example}\label{pp12-1}
    \ 
    \begin{enumerate}[label=\normalfont(\arabic*)]
    
    \item Let $t \in T$, where $T$ is as in~\autoref{b036}~\ref{b036-1}. Then, 
    $$\Xi^{\bar{\mu}}( \{ t \}) = \prod_{i < \vert t \vert} \mu_{t \rest i}(\{  t \rest (i + 1) \}) = \prod_{i < \vert t \vert} \frac{1}{2}  = 2^{-|t|}.$$
    
    \item For $t \in T$, where $T$ is as in~\autoref{b036}~\ref{b036-2}, $\displaystyle\Xi^{\bar{\mu}}(\{ t \}) = 2^{-|t|} 2^{-\sum_{i<|t|}t(i)}$.

    \item In the case of~\autoref{b036}~\ref{b036-3}, for $n<\omega$ and $t\in {}^n\omega$, $\Xi^{\bar{\mu}}(\{t\}) = 1$ if $t(i)=5$ for all $i<n$, otherwise $\Xi^{\bar{\mu}}(\{t\}) = 0$.
    \end{enumerate}    
\end{example}

Based on~\autoref{b036.1}, we can prove that a well-pruned probability tree induces a probability space structure at each level.

\begin{theorem}\label{b035}
    If $\la T,\bar\mu\ra$ is a well-pruned probability tree then, for any $n < \alt(T)$, $\Xi^{\bar{\mu}}_n$ is a probability measure on $\calP(\Lv_{n}(T))$.  
\end{theorem}

\begin{PROOF}{\ref{b035}}
    We proceed by induction on $n < \alt(T)$. 
    The case $n=0$ follows by \autoref{pp131}~\ref{pp131a} because $\Lv_0(T)=\{\la\,\ra\}$.     
    Now, assume that $n + 1 < \hgt(T)$ and that $\langle \Lv_{n}(T), \calP(\Lv_{n}(T)), \Xi_{n}^{\bar{\mu}} \rangle$ is a probability space. Since $\Lv_{n+1}(T)=\bigcup_{t\in\Lv_n(T)}\suc_T(t)$ is a disjoint union and $\Lv_n(T)\subseteq T\menos\max(T)$ (the latter because $T$ is well-pruned and $n+1<\hgt(T)$), by \autoref{pp131}~\ref{pp131b} we get:
    \[\Xi^{\bar\mu}(\Lv_{n+1}(T)) = \sum_{t\in\Lv_n(T)}\Xi^{\bar\mu}(\suc_T(t)) = \sum_{t\in T}\Xi^{\bar\mu}(\{t\}) =\Xi^{\bar\mu}( \Lv_{n}(T))=1.\]
    % Then, using that $\langle T, \bar{\mu} \rangle$ is a probability tree and~\autoref{pp13}, we get: 
    % \begin{equation*}
    %     \begin{split}
    %         \Xi_{n+1}^{\bar{\mu}}(\Lv_{n+1}(T)) & =  \sum_{t \in \Lv_{n+1}(T)}  \Xi_{n+1}^{\bar{\mu}}(\{ t \}) = \sum_{t \in \Lv_{n}(T)}  \left( \ \sum_{s \in \suc_{T}(t)}  \Xi_{n+1}^{\bar{\mu}} ( \{ s \}) \right)\\ 
    %         & = \sum_{t \in \Lv_{n}(T)} \left( \sum_{s \in \suc(t)}   \mu_{t}( \{ s \})   \Xi_{n}^{\bar{\mu}}(\{ t \}) \right) = \sum_{t \in \Lv_{n}(T)}  \Xi_{n}^{\bar{\mu}}( \{ t \} ) = 1.
    %     \end{split}
    % \end{equation*}
    Thus, $(\Lv_{n+1}(T), \calP(\Lv_{n+1}(T)), \Xi_{n+1}^{\bar{\mu}})$ is a probability space. 
\end{PROOF}

\autoref{b035} may not hold when $T$ is not well-pruned. We leave this discussion to the following subsection (\autoref{b035-4}).

% \begin{remark}
%     When $\la T,\bar\mu\ra$ is a probability tree, the assumption that $T$ is well-pruned is essential to have that $\Xi^{\bar\mu}_n$ is a probability measure for each $n<\hgt(T)$. However, if we allow $T$ not being well-pruned, we still have that each $\Xi^{\bar\mu}_n$ is a measure but not necessarily of probability, i.e.\ $\Xi^{\bar\mu}_n(\Lv_n(T))$ can happen: if $T$ is not well-pruned then the minimal $n_0<\omega$ such that $\max T\cap \Lv_{n_0}(T)\neq\emptyset$ satisfies $n_0+1<\hgt(T)$. One can show that, for $n<\hgt(T)$, $\Xi^{\bar\mu}_n$ is not a probability measure iff $n\geq n_0$. 

%     Even without assuming that a probability tree is well-pruned, we can determine a probability space on every \emph{front} of the tree.
% \end{remark}

\subsection{Inductive probability measures}\ 

%In \autoref{b035}, appeared the following condition: ``$\Xi_{\vert t \vert + 1}^{\bar{\mu}}(\suc_{T}(t)) = \Xi_{\vert t \vert}^{\bar{\mu}}(\{ t \})$ whenever $t \in T \menos \max(T)$''. This will be fundamental for the analysis of the structure of probability trees (see~\autoref{ppp14.T}). For this reason, we introduce the notion of \emph{inductive probability measures}.

The properties presented in \autoref{pp131} are fundamental for the analysis of the structure of probability trees. This motivates the following notion of \emph{inductive probability measure}.

\begin{definition}\label{pp46}
    \ 
    \begin{enumerate}
        \item Let $T\in \calT$. We say that $\Xi$ is an \emph{inductive probability measure on $T$} if it is a measure on $\pts(T)$ such that $\Xi(\{\la\,\ra\})= 1$ and $\Xi(\suc_{T}(t)) = \Xi(\{ t \})$ whenever $t\in T\menos \max(T)$.

        For any $n < \hgt(T)$, $\Xi_{n}\coloneqq \Xi\frestr\pts(\Lv_n(T))$, which is a measure on $\pts(\Lv_n(T))$.

        \item  Denote by $\IP$ the collection of inductive probability measures $\Xi$ in some $T\in\calT$. Notice that $T$ is uniquely determined by $\Xi$, so it will be denoted by $T_{\Xi}$.

        \item Define the function $\Pi \colon \TP \to \IP$ such that, for any $\bar{\mu} \in \TP$, $\Pi(\bar{\mu}) \coloneqq \Xi^{\bar{\mu}}$, which is well-defined by virtue of~\autoref{pp131}. Notice that $T_{\Xi^{\bar\mu}} = T_{\bar\mu}$.
    \end{enumerate}
\end{definition}

Notice that the proof of \autoref{b035} only uses that $\Xi^{\bar\mu}$ is an inductive probability measure. For this reason, the same proof yields:

\begin{theorem}\label{pp46.1}
    If $\Xi$ is an inductive probability measure on a well-pruned tree $T$, then $\Xi_n$ is a probability measure on $\pts(\Lv_n(T))$ for all $n<\alt(T)$.
\end{theorem}

As in the case of \autoref{b035}, the previous theorem may not hold when $T$ is not well-pruned. To understand this, we generalize this theorem by using the following notion. 

\begin{definition}\label{b035-1}
    Let $T$ be a tree of sequences. A set $A\subseteq T$ is \emph{a front of $T$} if it satisfies the following:
    \begin{enumerate}[label = \normalfont (\roman*)]
         \item Any pair of members of $A$ are incompatible in $T$.
        \item Every maximal branch of $T$ intersects $A$, i.e.\ for any $x\in[T]$ there is some $n<\omega$ such that $x\frestr n\in A$.
    \end{enumerate}
    For example, for any $n<\alt(T)$, $\Fr_n(T)\coloneqq \Lv_n(T)\cup \set{s\in\max(T)}{|s|<n}$ (which is a disjoint union) is a front of $T$. 
\end{definition}

\begin{fact}\label{b035-2}
    Let $T$ be a tree of sequences. 
    \begin{enumerate}[label = \normalfont (\alph*)]
        \item\label{b035-2a} $T$ is well-pruned iff, for any $n<\hgt(T)$, $\Lv_n(T)$ is a front of $T$. 
        \item\label{b035-2b} Assume that $T$ is not well-pruned and let $n_0<\hgt(T)$ be the minimum number such that $\Lv_{n_0}(T)\cap \max(T)\neq \emptyset$. Then $n_0+1<\hgt(T)$ and, for any $n<\hgt(T)$, $\Lv_n(T)$ is a front of $T$ iff $n\leq n_0$.
    \end{enumerate}
\end{fact}
\begin{PROOF}{\ref{b035-2}}
    Regardless of whether $T$ is well-pruned or not, we can find the maximum $\gamma\leq\hgt(T)$ (which can be $\omega$) satisfying that $T|_\gamma\coloneqq \set{t\in T}{|t|<\gamma}$ is well-pruned. Notice that $\gamma>0$ and $T$ is well-pruned iff $\gamma=\hgt(T)$. Also, if $T$ is not well-pruned, then $n_0+1=\gamma<\hgt(T)$. 

    Therefore, it is enough to show that, for any $n<\hgt(T)$, $\Lv_n(T)$ is a front of $T$ iff $n< \gamma$. If $x\in[T]$ then its length must be ${\geq}\, \gamma$ (otherwise $T|_\gamma$ would not be well-pruned), so $x\frestr n\in\Lv_n(T)$ for any $n<\gamma$. This shows the implication from right to left.

    For the converse, assume that $\gamma\leq n<\hgt(T)$. Then $\gamma+1\leq \hgt(T)$ and $T|_{\gamma+1}$ is not well-pruned, so $\Lv_{n_0}(T)\cap\max T\neq\emptyset$. Any $x$ in this set has length $n_0<n$, so it does not have nodes of length $n$ below it. Thus, $\Lv_n(T)$ is not a front. 
\end{PROOF}

\autoref{pp46.1} is generalized as follows. 

\begin{theorem}\label{b035-3}
    If $\Xi\in\IP$ and $A$ is a front of $T\coloneqq T_\Xi$, then $\Xi\frestr\pts(A)$ is a probability measure. Moreover, 
    \begin{equation}\label{b035-3eq}
        \Xi(\{t\}) = \Xi(\set{s\in A}{t\subseteq s}) \text{ for any $t\in T$ below some member of $A$}.
    \end{equation}
\end{theorem}
\begin{PROOF}{\ref{b035-3}}
    For $t\in T$, let $A_{\geq t}\coloneqq \set{s\in A}{t\subseteq s}$. Since $A_{\geq \la\,\ra}=A$, \eqref{b035-3eq} and \autoref{pp131}~\ref{pp131a} imply that $\Xi(A) = \Xi(\{\la\,\ra\})=1$. Hence, it is enough to show \eqref{b035-3eq}. 
    We provide two proofs of \eqref{b035-3eq}; the first is presented below, and the second is in~\autopageref{proofb035-3}.

    Let $S$ be the set of nodes in $T$ below some member of $A$. This is a subtree of $T$ and $\lim S=\emptyset$: if $x\in \lim S$ then $x\in\lim T$, so $x\frestr n\in A$ for some $n<\omega$ because $A$ is a front, but all members of $A$ are pairwise incompatible, so $x\frestr m\notin S$ for all $m\geq n$, contradicting that $x\in\lim S$.

    Therefore, $S$ is a well-founded tree by \autoref{a055}. It is enough to show that $S^\sim\coloneqq\set{t\in S}{\Xi(\{t\}) \neq \Xi(A_{\geq t})} = \emptyset$. Assume the contrary, so $S^\sim$ contains a maximal element $t$. If $t\in A$ then $A_{\geq t}=\{t\}$, so $\Xi(\{t\}) = \Xi(A_{\geq t})$, that is $t\notin S^\sim$. Thus $t\notin A$. This implies that $\suc_T(t)\subseteq S$ and, since $t$ is maximal in $S^\sim$, $\suc_T(t)\cap S^\sim =\emptyset$, so $\Xi(\{t'\}) = \Xi(A_{\geq t'})$ for all $t'\in\suc_T(t)$. On the other hand, $A_{\geq t}=\bigcup_{t'\in\suc_T(t)}A_{\geq t'}$ is a disjoint union, so
    \[\Xi(A_{\geq t}) = \sum_{t'\in\suc_T(t)}\Xi(A_{\geq t'}) = \sum_{t'\in\suc_T(t)}\Xi(\{ t'\}) = \Xi(\suc_T(t))=\Xi(\{t\}),\]
    where the last equality holds by \autoref{pp131}~\ref{pp131b}. But this contradicts that $t\in S^{\sim}$.
\end{PROOF}

Notice that \autoref{pp46.1} (and \autoref{b035}) follow from \autoref{b035-2} and \autoref{b035-3}.

\begin{example}\label{b035-4}
    Assume that $\Xi\in\IP$ and that $T\coloneqq T_\Xi$ is not well-pruned. Let $n_0$ be as in \autoref{b035-2}~\ref{b035-2b}. For any $n_0<n<\alt(T)$, $\Fr_n(T)$ is a front of $T$, so $\Xi\frestr\pts(L_n)$ is a probability measure. Therefore, a necessary condition for $\Xi(\Lv_n(T))=1$ (i.e.\ $\Xi_n$ is a probability measure) is that $\Xi(\{t\})=0$ for all $t\in\max(T)$ with $|t|<n$. The latter condition does not always hold.
\end{example}

Given an inductive probability measure in $T$, we can compute the probability of $T_{\geq t}$ in each level of $T$ in terms of the probability of $t$. 

\begin{corollary}\label{ppp12}
    If $T \in \calT$ and $\Xi$ is an inductive probability measure in $T$, then for any $t \in T$ and $\vert t \vert\leq n <\hgt(T)$, if $T_{\geq t} \cap \Lv_{n}(T)$ is a front of $T_{\geq t}$ (which holds when $T$ is well-pruned), then $\Xi_{n} \left( T_{\geq t} \cap \Lv_{n}(T) \right) = \Xi(\{ t \})$.
\end{corollary}
\begin{PROOF}{\ref{ppp12}}
    Apply~\eqref{b035-3eq} to any front of $T$ containing $T_{\geq t} \cap \Lv_{n}(T)$, e.g.\ $\Fr_n(T)$.
\end{PROOF}

We will prove later that $\Pi$ is surjective (see~\autoref{pp14}), however, it is evident that $\Pi$ is not one-to-one, as it does not matter how a probability tree is defined above a measure-zero node. While it might seem reasonable to eliminate measure-zero nodes and restrict to probability trees with strictly positive measures, this approach is not ideal, as the applications often involve trees with measure-zero nodes. Instead, we will isolate the \emph{positive part} of a tree with respect to measures in $\TP$ and $\IP$. 

\begin{definition}\label{pp40}
    Let $\bar{\mu} \in \TP$ and $\Xi \in \IP$. 

    \begin{enumerate} 
        \item We say that $\bar{\mu}$ is \emph{positive} if, for any $t \in T_{\bar\mu} \menos \max(T_{\bar\mu})$, and $s \in \suc_{T_{\bar\mu}}(t)$, $\mu_{t}(\{ s \})$ is positive. Similarly, $\Xi$ is \emph{positive} if, for any $t \in T_\Xi$, $\Xi(\{ t \}) > 0$. 
        
        \item If $\langle T, \bar{\mu} \rangle$ is a probability tree, we say that it is \emph{positive} if $\bar{\mu}$ is positive.

        \item $N_{\Xi} \coloneqq \{ t \in T \colon \Xi(\{ t \}) = 0 \}$ and $N_{\bar\mu}\coloneqq N_{\Xi^{\bar\mu}}$, which are called the \emph{null part} of $\Xi$ and $\bar\mu$, respectively.

        \item $T_{\bar{\mu}}^{+} \coloneqq T_{\bar{\mu}} \menos N_{\bar{\mu}}$ and $T_{\Xi}^{+} \coloneqq T_{\Xi} \menos N_{\Xi}$, which are called the \emph{positive part} of $\bar\mu$ and $\Xi$, respectively. 

        \item For $S \subseteq T_{\bar{\mu}}$, %and $R \subseteq T_{\Xi}$.   
        set $\bar{\mu} \rest S \coloneqq \langle \mu_{t} \rest [S \cap \suc_{T}(t)] \colon t \in S \menos \max(S) \rangle$. %Similarly,  $\Xi \rest R \coloneqq \langle \Xi_{n} \rest [R \cap \Lv_{n}(T)] \colon n  < \alt(T) \text{ and } R\cap \Lv_n(T)\neq\emptyset \rangle.$

        \item $\bar{\mu}_{+} \coloneqq \bar{\mu} \rest T_{\bar{\mu}}^{+}$ and $\Xi_{+} \coloneqq \Xi \rest T_{\Xi}^{+}$. 

        % \item $\bar{\mu}_{+} \coloneqq \langle \mu_{t, +} \colon t \in T \menos \max(T) \rangle$, where for all $t \in T \menos \max(T)$, $\mu_{t, +} \coloneqq \mu_{t} \rest T_{\bar{\mu}}^{+}$. Similarly, $\Xi_{+} \coloneqq \langle \Xi_{n, +} \colon n < \alt(T_{\Xi}) \rangle$, where for any $n < \alt(T)$, $\Xi_{n, +} \coloneqq \Xi_{n} \rest T_{\Xi}^{+}$.  \Andres{Arreglar. Debe indizar con los positivos.} \Andres{Además, cambiar el orden para definir $\mu_{+}$ con las restricciones.}

        \item $\TP_{+}$ is the collection of all positive $\bar{\mu} \in \TP$. Similarly, $\IP_{+}$ is the collection of all positive $\Xi \in \IP$. 
        
        \item We define the functions:
        \begin{align*}
            \varphi_{\TP} \colon & \TP \to \TP_{+} &&  \text{by $\varphi_{\TP}(\bar{\mu}) \coloneqq \bar{\mu}_{+}$,}\\
            \varphi_{\IP} \colon & \GP \to \IP_{+} && \text{by $\varphi_{\IP}(\Xi) \coloneqq \Xi_{+}$, and}\\
            \Pi_{+} \colon & \TP_{+} \to \IP_{+}  && \text{is  $\Pi \rest \TP_{+}$.}
        \end{align*}
        \autoref{pp45-1} below justifies that the co-domains of these functions, and the domain of $\Pi_+$, are as indicated.   
    \end{enumerate}
\end{definition}

Next, we list some basic properties of the notions introduced in~\autoref{pp40}. 

\begin{fact}\label{pp45-1}
    Let $\bar{\mu} \in \TP$ and $\Xi \in \IP$. Then:  

    \begin{enumerate}[label=\normalfont(\alph*)] 
        \item\label{pp45.0} $T^+_{\bar{\mu}} = T^+_{\Xi^{\bar{\mu}}}$.
    
        \item\label{pp45.1} 
        $T^{+}_{\Xi}$ is a subtree of $T_{\Xi}$ (so $T^+_\Xi\in\calT$) and $\max T^+_\Xi= T^+_\Xi\cap \max T_\Xi$. Moreover, if $T_{\Xi}$ is well-pruned then so is $T^+_{\Xi}$ and $\hgt(T^{+}_{\Xi}) = \hgt(T_{\Xi})$. A similar result holds for $\bar\mu$.

        \item\label{pp45.5} $ \bar{\mu}_{+} \in \TP$ and $\Xi_{+} \in \IP$, also  $T_{\bar{\mu}_{+}} = T_{\bar{\mu}}^{+}$ and $T_{\Xi_+} = T_{\Xi}^{+}$. 

        \item\label{pp45.6} $\Xi^{\bar{\mu}_{+}}=\Xi_{+}^{\bar{\mu}}$.  

        \item\label{pp45.7} $\Xi$ is positive iff $\Xi = \Xi_{+}$ iff $T_{\Xi} = T_{\Xi}^{+}.$
    
        \item\label{pp45.2} $\bar{\mu}$ is positive iff $\bar{\mu} = \bar{\mu}_{+}$ iff $T_{\bar{\mu}} = T_{\bar{\mu}}^{+}$.

        \item\label{pp45.3} $\bar{\mu}$ is positive iff\/ $\Xi^{\bar{\mu}}$ is positive.

        %\item\label{pp45.4} Let $\Xi' \in \IP$. If $T_{\Xi} = T_{\Xi'}$ and $\Xi =\Xi'_{+}$, then $\Xi = \Xi'$. \Diego{Redundante con~\ref{pp45.7}}
    \end{enumerate}
\end{fact}

\begin{PROOF}{\ref{pp45-1}}
    \ref{pp45.0}: Immediate because $T_{\bar\mu}=T_{\Xi^{\bar\mu}}$ and $N_{\bar\mu}=N_{\Xi^{\bar\mu}}$. 

    \ref{pp45.1}: We just deal with $\Xi$ since it implies the result for $\bar\mu$ by~\ref{pp45.0}. First notice that $\Xi(\la\,\ra)=1$, so $\la\,\ra\in T^+_{\Xi}$, and it is clear that $s\subseteq t$ in $T_\Xi$ implies $\Xi(t)\leq\Xi(s)$. Therefore, $T^+_\Xi$ is a subtree of $T_\Xi$.

    Now let $t \in T_{\Xi}^{+} \menos \max(T_{\Xi})$, so $\Xi(\{t\})>0$. Since we have that  $\Xi(\{t\})=\sum_{s\in\suc_{T_{\Xi}}(t)}\Xi(\{s\})$, there must be some $s\in\suc_{T_{\Xi}}(t)$ such that $\Xi(\{s\})>0$, so $s\in T^+_{\Xi}$ and, hence $t\notin\max T^+_\Xi$. This indicates that $\max T^+_\Xi \subseteq \max T_\Xi$ (the contention $T^+_\Xi\cap \max T_\Xi\subseteq \max T^+_\Xi$ is trivial). This implies that, whenever $T_\Xi$ is well-pruned, $T^+_{\Xi}$ is well-pruned with the same height. 
    
    %such that $T_{\bar{\mu}}^{+} \cap \suc_{T}(t) = \suc_{T_{\bar{\mu}}^{+}}(t) = \emptyset$. Therefore, for any $s \in \suc_{T}(t)$, we have that $\mu_{t}(\{ s \}) = 0$, which is not possible because $\suc_{T}(t) \neq \emptyset$ and $\mu_{t}$ is a probability measure. Thus $T_{\Xi}$ is well-pruned. On the other hand, $ \alt(T_{\bar{T}}) \leq \alt(T_{\Xi}^{+})$ holds because, otherwise we can find $t \in T_{\bar{\mu}}^{+}$ such that $\suc_{T}(t) \neq \emptyset$ and $\mu_{t}(\suc_{T}(t)) = 0$. 

    \ref{pp45.5}: It is clear that $\Xi_+$ is a measure on $\pts(T^+_\Xi)$ and $\Xi_+(\{\la\,\ra\})=\Xi(\{\la\,\ra\})=1$. Now, for $t\in T^+_\Xi\menos\max T_\Xi$, since the nodes in $\suc_{T_\Xi}(t)\menos \suc_{T^+_\Xi}(t)$ have measure zero,
    $$
    \Xi_{+}(\{ t \}) = \Xi(\suc_{T_\Xi}(t)) = \Xi(\suc_{T_{\Xi}^{+}}(t)) = \Xi(\suc_{T_{\Xi}^{+}}(t)) = \Xi_+(\suc_{T_{\Xi}^{+}}(t)).
    $$
    Thus, $\Xi_{+} \in \IP$ and $T_{\Xi_+}=T^+_\Xi$. The proof for $\bar\mu_+$ is similar.

    \ref{pp45.6}: By~\ref{pp45.5} and~\ref{pp45.0}, $T_{\Xi^{\bar\mu}_+}= T^+_{\Xi^{\bar\mu}} = T^+_{\bar\mu} = T_{\bar\mu_+}= T_{\Xi^{\bar\mu_+}}$. The equality $\Xi^{\bar\mu_+}(\{t\}) = \Xi^{\bar\mu}(\{t\})$ for $t\in T_{\bar\mu_+}$ is a straightforward calculation using \autoref{b036.1}.

    Item~\ref{pp45.7} and~\ref{pp45.2} clear by \autoref{pp40}, and~\ref{pp45.3} follows by~\ref{pp45.7},~\ref{pp45.2}, and~\ref{pp45.0}.
\end{PROOF}

We can use $\mu_{+}$ and $\nu_{+}$ to characterize when $\Pi(\bar{\mu})$ and $\Pi(\bar{\nu})$ are equal. Furthermore, this characterization enables us to show that, by restricting to positive trees, $\Pi_{+}$ establishes a bijection between $\TP_{+}$ and $\IP_{+}$ (see also~\autoref{ppp17}).

\begin{lemma}\label{pp45}
    Let $\bar{\mu}, \bar{\nu} \in \TP$. Then $\Pi(\bar{\mu}) = \Pi(\bar{\nu})$ iff $\bar{\mu}_{+} = \bar{\nu}_{+}$ and $T_{\bar{\mu}} = T_{\bar{\nu}}$. As a consequence,  $\Pi_{+}$ is a one-to-one function. 
\end{lemma}

\begin{PROOF}{\ref{pp45}}
    Assume that $\Pi(\bar{\mu}) = \Xi^{\bar{\mu}} = \Xi^{\bar{\nu}} =  \Pi(\bar{\nu})$. It is clear that $T_{\bar{\mu}} = T_{\bar{\nu}}$ and $N_{\Xi^{\bar{\mu}}} = N_{\Xi^{\bar{\nu}}}$, therefore  $T^+ \coloneqq T_{\Xi^{\bar{\mu}}}^{+} = T_{\Xi^{\bar{\nu}}}^{+}$. On the other hand, let $t \in T^+ \menos \max(T^+)$ and $s \in \suc_{T^+}(t)$. By~\autoref{pp13} we have 
    $$\mu_{t}(\{ s \}) = \frac{\Xi^{\bar{\mu}}(\{ s \})}{\Xi^{\bar{\mu}}(\{ t \})} = \frac{\Xi^{\bar{\nu}}(\{ s \})}{\Xi^{\bar{\nu}}(\{ t \})} = \nu_{t}(\{ s \}).$$ Thus, $\bar{\mu}_{+} = \bar{\nu}_{+}$. 

    Conversely, if $T\coloneqq T_{\bar\mu} = T_{\bar\nu}$ and $\bar\mu_+ = \bar\nu_+$, then $\Xi_+^{\bar\mu} = \Xi^{\bar\mu_+} = \Xi^{\bar\nu_+} = \Xi_+^{\bar\nu}$ by \autoref{pp45-1}. Hence $T^+\coloneqq T^+_{\bar\mu} = T^+_{\bar\nu}$ and $N\coloneqq N_{\bar\mu} = N_{\bar\nu}$. Since the nodes in $N$ have measure zero with respect both $\Xi^{\bar\mu}$ and $\Xi^{\bar\nu}$, we conclude that $\Xi^{\bar\mu}=\Xi^{\bar\nu}$.
    %
    % To prove the converse, we prove by induction on $n<\alt(T)$ that $\Xi^{\bar\mu}(\{t\}) = \Xi^{\bar\nu}(\{t\})$ for all $t\in\Lv_n(T)$, where $T \coloneqq T_{\bar{\mu}} = T_{\bar{\nu}}$. Now, assume that for any $n < \alt(T)$ and $s \in \Lv_{n}(T)$, $\Xi_{n}^{\bar{\mu}}(\{ s \}) = \Xi_{n}^{\bar{\nu}}(\{ s \})$ and let $t \in \Lv_{n + 1}(T)$. If $t \rest n \in N_{\bar{\mu}}$ then $\Xi_{n + 1}^{\bar{\mu}}(\{ t \}) = 0 = \Xi_{n+1}^{\bar{\nu}}(\{ t \})$, so assume that $t \rest n \notin N_{\bar{\mu}}$.  By~\autoref{pp13} and the induction hypothesis, 
    % $$
    % \Xi_{n+1}^{\bar{\mu}}(\{ t \}) = \Xi_{n}^{\bar{\mu}}(\{ t \}) \mu_{t \rest n} (\{ t \rest n \}) = \Xi_{n}^{\bar{\nu}}(\{ t \}) \nu_{t \rest n} (\{ t \}) = \Xi_{n+1}^{\bar{\nu}}(\{ t \}). $$ 
    % This shows that $\Pi(\bar{\mu}) = \Pi(\bar{\nu})$. 
\end{PROOF}

The condition in~\autoref{pp45} that characterizes when $\Pi(\bar{\mu})$ and $\Pi(\bar{\nu})$ are equal, motivates the definition of the following equivalence relations on $\TP$ and $\IP$.

\begin{definition}\label{pp45-i}
    \ 
    \begin{enumerate}[label=\normalfont(\arabic*)]
        \item We say that $\bar{\mu}, \bar{\nu} \in \TP$ are \emph{positive equivalent}, denoted by $\bar{\mu} \equiv_{\TP} \bar{\nu}$, iff $\bar{\mu}_{+} = \bar{\nu}_{+}$. It is clear that this is an equivalence relation. Denote by $[\bar{\mu}]_{\equiv_{\TP}} \coloneqq [\bar{\mu}]_{\TP}$ the $\equiv_{\TP}$-equivalence class of $\bar{\mu}$. 

        \item Similarly,  $\Xi, \bar{\rho} \in \IP$ are \emph{positive equivalent}, denoted by $\Xi \equiv_{\IP} \bar{\rho}$, iff $\Xi_{+} = \bar{\rho}_{+}$. It is clear that this is an equivalence relation. Denote by $[\Xi]_{\equiv_{\IP}} \coloneqq [\Xi]_{\IP}$ the $\equiv_{\IP}$-equivalence class of $\Xi$.

        \item  $\pi_{\TP} \colon \TP \to \TP / \equiv_{\TP}$  and $\pi_{\IP} \colon \IP \to \IP / \equiv_{\IP}$  are defined by $\pi_{\TP}(\bar{\mu}) \coloneqq [\bar{\mu}]_{\TP}$ and $\pi_{\IP}(\Xi) \coloneqq [\Xi]_{\IP}$, respectively. Furthermore, $\pi_{\TP}^{+} \coloneqq \pi_{\TP} \rest \TP_{+}$ and similarly, $\pi_{\IP}^{+} \coloneqq \pi_{\IP} \rest \IP_+$.
        
        \item Define $\hat{\Pi} \colon \TP / \equiv_{\TP} \to \IP / \equiv_{\IP}$  by $\hat{\Pi}([\bar{\mu}]_{\TP}) \coloneqq [\Pi(\bar{\mu})]_{\IP}$, which is well-defined by virtue of~\autoref{pp45-j} below. 
    \end{enumerate}
\end{definition}

Notice that, trivially $\bar{\mu}_{+} \equiv_{\TP} \bar{\mu}$ and $\Xi_+\equiv_\IP \Xi$ for $\bar\mu\in\TP$ and $\Xi\in\IP$. In fact, $\bar{\mu}_{+}$ is the canonical class-representative of $\bar{\mu}$. Similarly for inductive probability measures.  

We can characterize when $\bar{\mu}$ and $\bar{\nu}$ are positively equivalent in terms of $\Xi^{\bar{\mu}}$ and $\Xi^{\bar{\nu}}$.

\begin{lemma}\label{pp45-j}
   If $\bar{\mu}, \bar{\nu} \in \TP$, then the following statements are equivalent.
       \begin{multicols}{2}
       \begin{enumerate}[label=\normalfont(\roman*)]
           \item\label{pp45-j.1} $\bar{\mu} \equiv_{\TP} \bar{\nu}$.
           
           \item\label{pp45-j.2} $\Xi^{\bar{\mu}} \equiv_{\IP} \Xi^{\bar{\nu}}$. 
           
           \item\label{pp45-j.3} $\hat{\Pi}([\bar{\mu}]_{\TP}) = \hat{\Pi}([\bar{\nu}]_{\TP})$.
       \end{enumerate}
       \end{multicols}
    As a consequence, $\hat{\Pi}$ is well-defined and one-to-one. 
\end{lemma}

\begin{PROOF}{\ref{pp45-j}}
    Let $\bar{\mu}, \bar{\nu} \in \TP$. Since $\Pi_{+}$ is one-to-one and by~\autoref{pp45-1}~\ref{pp45.6}, we have: 
    $$ \bar{\mu} \equiv_{\TP} \bar{\nu} \Leftrightarrow
    \bar{\mu}^{+} = \bar{\nu}^{+} \Leftrightarrow \Xi^{\bar{\mu}_{+}} = \Xi^{\bar{\nu}_{+}} \Leftrightarrow \Xi^{\bar{\mu}}_{+} = \Xi_{+}^{\bar{\nu}} \Leftrightarrow \Xi^{\bar{\mu}} \equiv_{\IP} \Xi^{\bar{\nu}} \Leftrightarrow \Pi(\bar{\mu}) \equiv_{\IP} \Pi(\bar{\nu}). 
    $$ 

    On the other hand,~\ref{pp45-j.2} $\Leftrightarrow$ \ref{pp45-j.3} follows immediately from the definition of $\hat{\Pi}$. 
\end{PROOF}

We explore the close relationship between $\IP$ and $\TP$, and even show that probability trees can be identified with inductive probability measures. To this end, we show how to construct probability trees from inductive probability measures. 

\begin{definition}\label{pp47}
    \ 
    \begin{enumerate}[label = \normalfont (\arabic*)]
        \item  Denote by $\GP$ the collection of pairs $(\Xi_{+},\langle \nu_{t} \colon t \in N_{\Xi} \rangle)$ such that $\Xi \in \IP$ and, for any $t \in N_{\Xi }$, $\nu_{t}$ is a probability measure on $\calP(\suc_{T}(t))$. %Notice that each $\bar{\eta} \in \GP$ is  uniquely determined by  $\Xi \in \IP$  and $\langle \nu_{t} \colon t \in N_{\Xi} \rangle$. 
        %Furthermore,  $\IP_{+} \subseteq \GP$. 
        Since $(\Xi,\la\,\ra)\in\GP$ for all $\Xi\in\IP_{+}$, we often identify a positive $\Xi$ with $(\Xi,\la\,\ra)$ and claim that $\IP_+\subseteq \GP$.

        \item\label{pp47.2}  Given $\bar{\eta} \in \GP$ determined by $\Xi \in \IP$ and $\langle \nu_{t} \colon t \in T_{\Xi} \rangle$, we define $\bar{\mu}^{\bar{\eta}} \in \TP$ as follows. For any $t \in T_{\Xi} \menos \max(T_{\Xi})$ and $s \in \suc_{T}(t)$,  
        $$ \mu_{t}^{\bar{\eta}}( \{ s \}) \coloneqq  
        \begin{cases} 
        \frac{\Xi( \{ s \})}{ \Xi(\{ t \})}, & \text{if $t \notin N_{\Xi}$,} \\[1ex] 
        \nu_{t}(\{ s \}), & \text{if  $t \in N_{\Xi}$.} 
        \end{cases}$$
        When $\Xi\in\IP_+$, denote $\bar\mu^\Xi\coloneqq \bar\mu^{(\Xi,\la\,\ra)}$.

        \item Define the function $\Delta \colon \GP \to \TP$ by $\Delta(\bar{\eta}) \coloneqq \bar{\mu}^{\eta}$, which is well-defined by virtue of~\autoref{pp14}~\ref{pp14.1} below.  
        
        \item Define $\gamma_{\GP} \colon \GP \to \IP$,  $\Delta_{+} \colon \GP \to \TP_{+}$ and $\gamma_{\GP}^{+} \colon \GP \to \IP^{+}$ by $\gamma_{\GP}(\bar{\eta}) \coloneqq \Xi$,  $\Delta_{+} (\bar{\eta}) \coloneqq \bar{\mu}^{\bar{\eta}}_{+}$, and $\gamma_{\GP}^{+}(\bar{\eta}) \coloneqq \Xi_{+}$, respectively, where $\bar{\eta}$ is determined by $\Xi$. 
    \end{enumerate}
\end{definition}

%The sequence $\bar{\mu}^{\bar{\eta}}$ above can be constructed by recursion. 

Given $\bar{\eta} \in \GP$ determined by $\Xi$, we can use $\bar{\mu}^{\bar{\eta}}$ to induce a probability tree structure on $T_{\Xi}$. This  will allow us to define a bijection between $\GP$ and $\TP$ 

\begin{lemma}\label{pp14}
    Let $\bar\eta \in \GP$ be determined by $\Xi \in \IP$ and $\seq{\nu_t}{t\in N_\Xi}$, and let $\bar\mu\in\TP$.
    
   \begin{multicols}{2}
        \begin{enumerate}[label=\normalfont(\alph*)]
        \item\label{pp14.1}  $\langle T_{\Xi}, \bar{\mu}^{\bar{\eta}} \rangle$ is a probability tree.  
        \item\label{pp14.2}$\Xi^{\bar{\mu}^{\bar{\eta}}} = \Xi$.

        \item\label{pp14.2.1} $\Pi$ is surjective. 

        \item\label{pp14.3} $(\Xi^{\bar\mu}_+,\bar\mu\frestr N_{\bar\mu})\in\GP$, $\bar\mu^{(\Xi^{\bar\mu}_+,\bar\mu\frestr N_{\bar\mu})} = \mu$.

        \item\label{pp14.2.f} If $\bar{\mu} \in \TP_{+}$ then $\bar{\mu}^{\Xi^{\bar{\mu}}} = \bar{\mu}$. 

        \item\label{pp14.6} If $\Xi'\in\IP_+$ then $\bar\mu^{\Xi'}\in\TP_+$.

        \item\label{pp14.5}$\Delta$ is bijective. 
        
        \item\label{pp14.4} $\Delta^{-1}(\bar\mu) = (\Xi^{\bar\mu}_+,\bar\mu\frestr N_{\bar\mu})$. 
    \end{enumerate}
   \end{multicols}
\end{lemma} 

\begin{PROOF}{\ref{pp14}}
    \ref{pp14.1}: Let $t \in T_{\Xi} \menos \max(T_{\Xi})$. On the one hand, if $t \notin N_{\Xi}$ then, by~\autoref{pp47}~\ref{pp47.2}, we have:
    \begin{equation*}
            \mu_{t}^{\bar{\eta}}(\suc_{T}(t)) =  \sum_{s \in \suc_{T_{\Xi}}(t)} \frac{\Xi( \{ s \})}{ \Xi(\{ t \})} = \frac{1}{ \Xi(\{ t \})} \sum_{s \in \suc_{T_{\Xi}}(t)} \Xi(\{ s \}) = 1.
    \end{equation*}

    On the other hand, if $\Xi(\{ t \}) = 0$ then 
    $ \mu^{\bar{\eta}}_{t}(\suc_{T}(t)) = \nu_{t}(\suc_{T}(t)) = 1$ because $\nu_{t}$ is a probability measure on  $\suc_{T}(t)$. 

    \ref{pp14.2}: It is clear that $T_{\Xi^{\bar\mu^{\bar\eta}}} = T_{\bar\mu^{\bar\eta}} = T_\Xi$. By induction on $n < \alt(T_{\Xi})$, we show that $\Xi^{\bar\mu^{\bar\eta}}(\{t\}) = \Xi(\{t\})$ for all $t\in\Lv_n(T_\Xi)$.  If $n = 0$ and $t \in \Lv_{0}(T_{\Xi})$ then $t = \langle \, \rangle$, so $\Xi^{\bar{\mu}^{\bar{\eta}}}(\{ t \}) = 1 = \Xi( \{ t \})$. Now assume that, for any $s \in \Lv_{n}(T_{\Xi})$, $\Xi^{\bar{\mu}^{\bar{\eta}}}(\{ s \}) = \Xi(\{ s \})$. For $t \in \Lv_{n+1}(T_{\Xi})$, by~\autoref{pp13} and induction hypothesis, we have that:    $$\Xi^{\bar{\mu}^{\bar{\eta}}}(\{ t \}) = \Xi^{\bar{\mu}^{\bar{\eta}}}(\{ t \rest n \}) \mu_{t \rest n}^{\bar{\eta}}(\{ t \}) =  \Xi( \{ t \rest n \}) \mu^{\bar{\eta}}_{t \rest n}(\{ t \}).$$
    Consider two possible cases. On the one hand, if $t \rest n \notin N_{\Xi}$ then, by the definition of $\bar{\mu}^{\bar{\eta}}$,
    $$ \Xi^{\bar{\mu}^{\bar{\eta}}}(\{ t \}) =  \Xi( \{ t \rest n \}) \mu_{t \rest n}^{\bar{\eta}}(\{ t \}) = \Xi(\{ t \rest n \}) \frac{\Xi(\{ t \})}{\Xi(\{ t \rest n \})} = \Xi(\{ t \}). 
    $$
    On the other hand, if $t \rest n \in N_{\Xi}$, then  $ \Xi^{\bar{\mu}^{\bar{\eta}}}(\{ t \}) = 0 = \Xi(\{ t \}), $
    where the last equality holds because $\Xi(\suc_{T}(t\frestr n)) = \Xi(\{ t \rest n \})$ and, thus, $\Xi(\{ s \}) = 0$ for any $s \in \suc_{T}(t \rest n)$. 

    \ref{pp14.2.1}: For $\Xi' \in \IP$ consider any $\seq{ \nu_{t} }{t\in N_{\Xi'}}$ such that $ \bar{\sigma} \coloneqq  (\bar{\Xi}'_{+} , \seq{ \nu_{t} }{ t \in N_{\Xi'} } \in \GP$. Then, by~\ref{pp14.1} and~\ref{pp14.2}, $\bar{\mu}^{\bar{\sigma}} \in \TP$ and $\Pi(\bar{\mu}^{\bar{\sigma}}) = {\Xi'}^{\bar{\mu}^{\bar{\sigma}}} = \Xi'$. This shows that $\Pi$ is surjective. 

    \ref{pp14.3}: It is clear that $(\Xi^{\bar{\mu}}_+,\bar\mu\frestr N_\Xi) \in \IP$ is represented by $\Xi^{\bar\mu}$, so $T_{(\Xi^{\bar{\mu}}_+,\bar\mu\frestr N_\xi)} = T_{\Xi^{\bar\mu}} = T_{\bar\mu}$. For $t \in T_{\bar{\mu}} \menos \max(T_{\bar{\mu}})$ and $s \in \suc_{T}(t)$, by~\autoref{pp47}~\ref{pp47.2}, if $\Xi^{\bar\mu}(\{t\})>0$ then $\mu_{t}^{{(\Xi^{\bar{\mu}}_+,\bar\mu\frestr N_\xi)}}(\{ s \}) = \frac{\Xi^{\bar{\mu}}(\{ s \})}{\Xi^{\bar{\mu}}(\{ t \})} = \mu_{t}(\{ s \})$, otherwise $\mu^{{(\Xi^{\bar{\mu}}_+,\bar\mu\frestr N_\xi)}}(\{s\})=\mu_t(\{s\})$.

    \ref{pp14.2.f}: By~\ref{pp14.3} because $\Xi^{\bar\mu}$ is positive when $\bar\mu$ is positive (so it is identified with $(\Xi^{\bar\mu},\la\,\ra)$).
 
    \ref{pp14.6}: By~\ref{pp14.2} and \autoref{pp45-1}~\ref{pp45.3}.

    \ref{pp14.5},~\ref{pp14.4}: By~\ref{pp14.3} and because $(\Xi^{\bar\mu^{\bar\eta}}_+,\bar\mu^{\bar\eta}\frestr N_{\Xi^{\bar\mu^{\bar\eta}}}) = (\Xi_+,\seq{\nu_t}{t\in N_\Xi}) = \bar\eta$ by~\ref{pp14.2}.
    %
    %If $\bar{\mu} \in \TP$ then $\Xi^{\bar{\mu}}_{+} \cup \langle \mu_{t} \colon t \in N_{\Xi^{\bar{\mu}}} \rangle$ is a preimage of $\mu$ under $\Delta$. This shows that $\Delta$ is surjective. On the other hand, let $\bar{\eta}, \bar{\sigma} \in \GP$ determined by $\Xi$, $\langle \nu_{t} \colon t \in N_{\Xi} \rangle$, and $\bar{\rho}$, $\langle \mu_{t} \colon t \in T_{\bar{\rho}} \rangle$ respectively, such that $\bar{\mu}^{\bar{\eta}} = \bar{\mu}^{\bar{\sigma}}$. By~\ref{pp14.2}, it follows that $\Xi = \bar{\rho}$ and therefore, $N_{\Xi} = N_{\bar{\rho}}$ and $T_{\Xi} = T_{\bar{\rho}}$, hence $T_{\Xi}^{+} = T_{\bar{\rho}}^{+}$.  As a consequence, $\langle \nu_{t} \colon t \in N_{\Xi} \rangle = \langle \mu_{t} \colon t \in N_{\bar{\rho}} \rangle$. Thus, $\bar{\eta} = \bar{\sigma}$. Thus, $\Delta$ is bijective. 
\end{PROOF}

As a consequence of~\autoref{pp14}~\ref{pp14.2.1} and \autoref{pp45-1}~\ref{pp45.3}, $\hat{\Pi}$ and $\Pi_{+}$ are surjective functions. Therefore, by~\autoref{pp45} and~\ref{pp45-j}, they are bijections. Furthermore,~\autoref{pp14}~\ref{pp14.2},~\ref{pp14.2.f} and~\ref{pp14.6} allow us to define explicitly the inverse of $\Pi_{+}$. 

\begin{corollary}\label{ppp17}
    The function $\Pi_{+} \colon \TP_{+} \to \IP_{+}$ is bijective and its inverse is $\Delta \rest \IP_{+}$. Moreover, $\hat{\Pi}$ is also bijective. 
\end{corollary}

\begin{corollary} 
    Let $\langle T, \bar{\mu} \rangle$ be a probability tree.  If $\Xi \in \IP$ is positive, then $\bar{\mu}^{\Xi} = \bar{\mu} $ iff $ \Xi = \Xi^{\bar{\mu}}$. 
\end{corollary}

\begin{figure}
    \centering
      \begin{tikzpicture}[>=stealth, node distance=3cm, scale=1, every node/.style={transform shape}]
    % Nodos con nombres correctos
    % Nodo nivel 1 (aumentado en Y)
    \node (B) at (-4, 4) {\footnotesize$\TP$};      % Nodo izquierdo nivel 2
    \node (IP) at (4, 4) {\footnotesize$\IP$};       % Nodo derecho nivel 2
    \node (TPc) at (-4, -3) {\footnotesize$\TP/_{\equiv_{\TP}}$};  % Nodo izquierdo nivel 4
    \node (IPc) at (4, -3) {\footnotesize$\IP/_{\equiv_{\IP}}$};    % Nodo derecho nivel 4
 
    \node (TP+) at (-2, -1) {\footnotesize$\TP_{+}$};  % Nodo izquierdo nivel 3 (punto medio en Y)
    
    \node (IP+) at (2, -1) {\footnotesize$\IP_{+}$};   % Nodo derecho nivel 3 (punto medio en Y)
    
    % Flechas con funciones
 % De nivel 1 a nodo derecho nivel 2
    \draw[->] (B) -- (TPc) node[midway, left] {$\pi_{\TP}$};            % De nodo izquierdo nivel 2 a nodo izquierdo nivel 4
    \draw[->] (IP) -- (IPc) node[midway, right] {$\pi_{\IP}$};           % De nodo derecho nivel 2 a nodo derecho nivel 4
            % De nodo izquierdo nivel 4 a nodo nivel 5
       % De nodo derecho nivel 4 a nodo nivel 5
    \draw[->] (B) -- (TP+) node[midway, right] {$\varphi_{\TP}$};            % De nodo izquierdo nivel 2 a nodo izquierdo nivel 3
    \draw[->] (IP) -- (IP+) node[midway, left] {$\varphi_{\IP}$};           % De nodo derecho nivel 2 a nodo derecho nivel 3
    \draw[->] (TP+) -- (TPc) node[midway, right] {$\pi_{\TP}^{+}$};            % De nodo izquierdo nivel 3 a nodo izquierdo nivel 4
    \draw[->] (IP+) -- (IPc) node[midway, left] {$\pi_{\IP}^{+}$};        % De nodo derecho nivel 3 a nodo derecho nivel 4
    
    % Flecha bidireccional entre G y H (nivel 3)
    \draw[<-] (IP+) -- (TP+) node[midway, above] {$\Pi_{+}$};         % Conexión bidireccional entre los nodos del nivel 3
    
    % Flecha bidireccional entre los nodos del nivel 5 (F)
    \draw[->] (TPc) -- (IPc) node[midway, above] {$\hat{\Pi}$};         % Conexión bidireccional entre los nodos del nivel 5
    
    % Flecha de izquierda a derecha entre los nodos del nivel 2 (B y C)
    \draw[->] (B) -- (IP) node[midway, above] {$\Pi$};          % Conexión de izquierda a derecha entre los nodos del nivel 2

    % Lo relacionado a GP 

      % Nodo del nivel 3.5 (entre nivel 2 y 3)
  \node (GP) at (0, 1.5) {\footnotesize$\mathcal{GP}$};  % Nodo nivel 3.5 (debajo de TP e IP)
  
  % Flechas con funciones
  \draw[<-] (B) -- (GP) node[midway, right] {$\ \ \Delta$};  % De nodo izquierdo nivel 2 a nodo nivel 3.5
  \draw[<-] (IP) -- (GP) node[midway, left] {$\gamma_{\GP}$};  % De nodo derecho nivel 2 a nodo nivel 3.5
  \draw[->] (GP) -- (TP+) node[midway, left] {$\Delta_{+}$};  % De nodo nivel 3.5 a nodo nivel 3
  \draw[<-] (IP+) -- (GP) node[midway, right] {$\gamma_{\GP}^{+}$};  % De nodo nivel 3.5 a nodo nivel 3

    \end{tikzpicture}
        \caption{Connections between the classes associated with $\TP$, $\IP$ and $\GP$.}
        \label{f45}
    \end{figure}
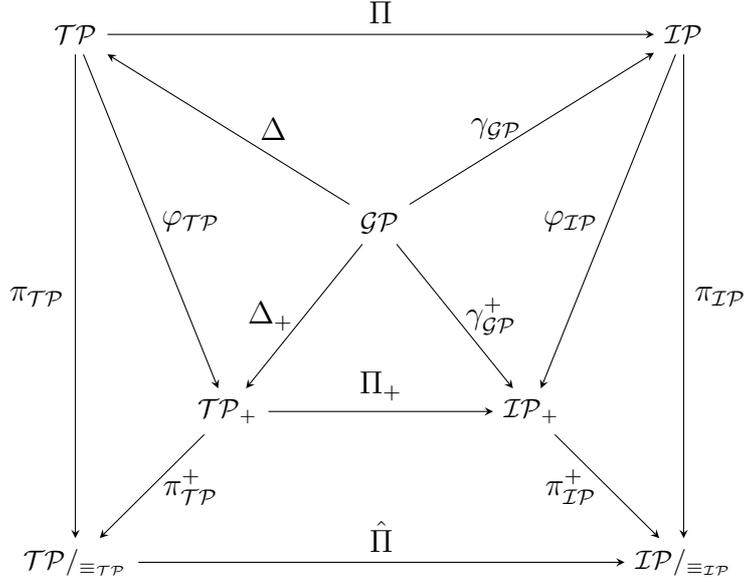

Finally, we can summarize the connections between $\TP$, $\IP$, and $\GP$ in~\autoref{f45}.

\begin{theorem}\label{pp51}
    \ 
    \begin{enumerate}[label=\normalfont(\alph*)]
        \item\label{pp51.2} $\Delta$, $\Pi_{+}$, $\hat{\Pi}$, $\pi_{\IP}^{+}$ and $\pi_{\TP}^{+}$ are bijective. 
        
        \item\label{pp51.1} All the diagrams in~\autoref{f45} are commutative.  As a consequence, any pair of paths with the same starting and ending points produce the same function.
        
        \item\label{pp51.0}  $\Pi$, $\pi_{\TP}$, $\pi_{\IP}$, $\Delta_{+}$, $\varphi_{\IP}, \varphi_{\TP}$, $\gamma_{\GP}$  and $\gamma_{\GP}^{+}$ are surjective.
    \end{enumerate}
\end{theorem}

\begin{PROOF}{\ref{pp51}}
    \ref{pp51.2}: It is clear that $\pi_{\IP}^{+}$ and $\pi_{\TP}^{+}$ are bijective functions. The result for $\Pi_{+}$, $\hat{\Pi}$ and $\Delta$ are from~\autoref{pp14}~\ref{pp14.5} and~\autoref{ppp17}.
    
    \ref{pp51.1}: It is enough to prove that the seven smallest sub-diagrams are commutative. Fix $\bar{\mu} \in \TP$, $\Xi \in \IP$ and $\bar{\eta} \in \GP$. The far left and right sub-diagrams commute because $\mu_+\equiv_\TP \mu$ and $\Xi_+\equiv_\IP \Xi$. It is clear that the sub-diagram at the bottom commutes.
    
    %To deal with the sub-diagram on the far left, notice that $ \pi_{\TP}^{+}(\bar{\mu}_{+}) = [\bar{\mu}_{+}]_{\TP} = [\bar{\mu}]_{\TP} = \pi_{\TP}(\bar{\mu}).$ For the sub-diagram on the far right,  we have $\pi_{\IP}^{+}(\Xi_{+}) = [\Xi_{+}]_{\IP} = [\Xi]_{\IP} = \pi_{\IP}(\Xi)$. For the bottom sub-diagram, by ~\autoref{pp14}~\ref{pp14.2}, if $\Xi \in \IP^{+}$ then  $\hat{\Pi}(\pi_{\TP}^{+}(\Pi_{+}^{-1}(\Xi))) = \hat{\Pi}([\bar{\mu}^{\Xi}]_{\TP}) = [\Xi^{\bar{\mu}^{\Xi}}]_{\IP} = [\Xi]_{\IP} = \pi_{\IP}^{+}(\Xi)$.
    
    Now, we deal with the four sub-diagrams bounded by the trapeze with vertices $\TP$, $\IP$, $\TP_{+}$ and $\IP_{+}$. For the left sub-diagram, $\varphi_{\TP}(\Delta(\bar{\eta})) = \varphi_{\TP}(\bar{\mu}^{\eta}) = \bar{\mu}^{\bar{\eta}}_{+} = \Delta_{+}(\bar{\eta})$. For the right, $\varphi_{\IP}(\gamma_{\GP}(\bar{\eta})) = \gamma_\GP^+(\bar\eta)$ by the definition of the maps. For the bottom sub-diagram, $\Pi_+(\Delta_+(\bar\eta))=\Xi^{\bar\mu^{\bar\eta}}_+ = \gamma^+_\GP(\bar\eta)$ by \autoref{pp14}~\ref{pp14.2} and \autoref{pp45-1}~\ref{pp45.6}.  Finally, for the top sub-diagram, $\Pi(\Delta(\bar{\eta})) = \Pi(\bar{\mu}^{\bar{\eta}}) = \Xi^{\bar{\mu}^{\bar{\eta}}} =  \gamma_{\GP}(\bar{\eta})$ by \autoref{pp14}~\ref{pp14.2}. 

    \ref{pp51.0}: It is clear that  $\pi_{\TP}$, $\pi_{\IP}$, $\varphi_{\IP}, \varphi_{\TP}$ are surjective, and $\Pi$ is surjective by~\autoref{pp14}. This implies, also using~\ref{pp51.2} and~\ref{pp51.1}, that $\Delta^+$, $\gamma_\GP$ and $\gamma^+_\GP$ are surjective.
\end{PROOF}

\subsection{Borel probability measures}\

So far, we have defined three distinct classes associated with probability trees, namely $\TP$, $\IP$, and $\GP$. However, there is a fourth which arises naturally by noting that every probability measure on $\calB([T])$ induces an inductive probability measure in $T$ (see~\autoref{pp50} and~\autoref{b120}). The new class is then the class of probability measures  on $\calB([T])$, which we introduce below.

\begin{definition}\label{pp160}
    We define $\BP$ as the collection of all probability measures on $\calB([T])$ for some $T \in \calT$. Notice that $T$ is uniquely determined by $\lambda$, so it will be denoted by $T_{\lambda}$. 
\end{definition}

Note that, when $[T] = \max(T)$ (e.g.\ in the case $\alt(T) < \omega$),  $\calB([T]) = \calP(\max(T))$. 

Now, we construct inductive probability measures derived from measures in $\BP$.  

\begin{definition}\label{pp50}
    For $\lambda \in \BP$, define the measure $\Xi^{\lambda}$ on $\pts(T_\lambda)$ determined by $\Xi^{\lambda}(\{ t \}) \coloneqq \lambda([t]_{T})$ for all $t\in T_\lambda$.  
     Furthermore, we define the function $\Psi \colon \mathcal{BP} \to \IP$ such that, for any $\lambda \in \mathcal{BP}$, $\Psi(\lambda) \coloneqq \Xi^{\lambda}$.
\end{definition}

We will see that $\Xi^{\lambda} \in \IP$ and therefore, $\Psi$ is well-defined in the sense that $\ran \Psi \subseteq \IP$. 

\begin{lemma}\label{b120}
    If $\lambda \in \BP$ then $\Xi^{\lambda}$ is an inductive probability measure in $T_{\lambda}$.  
\end{lemma}

\begin{PROOF}{\ref{b120}}
    Assume that $\lambda \in \BP$. For $t\in T_\lambda\menos \max(T_\lambda)$, since $[t]=\bigcup_{s\in\suc_{T_\lambda}(t)}[s]$ is a disjoint union, 
    $$ \Xi^{\lambda}(\{ t \}) = \lambda([t]) = \sum_{s \in \suc_{T}(t)} \lambda([s]) = \sum_{s \in \suc_{T}(t)} \Xi^{\lambda}(\{ s \}) = \Xi^{\lambda}(\suc_{T}(t)), 
    $$
    which proves that $\Xi^{\lambda} \in \IP$. 
\end{PROOF} 

We aim to expand~\autoref{f45} by incorporating the new class $\BP$, for which we need to establish connections between $\BP$ and the classes $\IP$ and $\TP$. The relationship between $\BP$ and $\IP$ poses no issues, as it is defined by $\Psi$. We will show that $\Psi$ is a bijection: one-to-one will follow by \autoref{b070}, while surjectivity is consequence of a connection between $\TP$ and $\BP$, namely, for $\bar\mu\in\TP$, we will construct a $\lambda^{\bar\mu}$ satisfying $\Xi^{\lambda^{\bar\mu}} = \Xi^{\bar\mu}$. The definition of $\lambda^{\bar\mu}$ is easy when $\lim T_{\bar{\mu}} =\emptyset$ (i.e.\ $[T_{\bar\mu}]=\max(T_{\bar\mu})$), as the only possible measure is $\lambda^{\bar{\mu}} \coloneqq \Xi^{\bar{\mu}}\frestr\max(T_{\bar\mu}) $, which is in $\BP$ with $T_{\lambda^{\bar\mu}} = T_{\bar\mu}$. However, we need more tools for the construction of $\lambda^{\bar\mu}$ when $\lim T_{\bar\mu} \neq\emptyset$ and to prove $\Xi^{\lambda^{\bar\mu}} = \Xi^{\bar\mu}$ (even in the case $\lim T_{\bar\mu} = \emptyset$). We present two ways to do this: the first uses \autoref{b035-3}, while the second is a concrete construction using a connection with the Lebesgue measure of the unit interval (see~\autoref{5.2},~\autoref{b060}). The second construction is one of the main applications of this paper, as it not only addresses the problem at hand but also provides an interesting connection between probability trees and the real line, which will also have consequences for representing  cardinals invariants of the continuum (see~\autoref{7}).

\begin{theorem}\label{ppp10}
    For every $\bar{\mu}\in \TP$ there is a unique $\lambda^{\bar{\mu}} \in \BP$ with $T_{\lambda^{\bar\mu}} = T_{\bar\mu}$ such that, for any $t \in T_{\bar\mu}$, $\lambda^{\bar{\mu}}([t]_{T}) = \Xi^{\bar{\mu}}(\{ t \})$, i.e.\ $\Xi^{\lambda^{\bar\mu}} = \Xi^{\bar\mu}$. 
\end{theorem}

\begin{PROOF}{\ref{ppp10}}
    Set $T\coloneqq T_{\bar\mu}$. Define
    \begin{align}\label{ppp10FT}
         \Fwf_{T} & \coloneqq \Fwf = \bigset{\bigcup_{t\in F}[t]}{F\subseteq\Fr_n(T),\ n<\alt(T)}.
    \end{align}
    This is an algebra of sets over $[T]$ and every $C\in\Fwf$ is clopen in $[T]$, so $\Fwf\subseteq \Bwf_T$. It is clear that $[t]\in \Fwf$ for all $t\in T$, so $\sigma(\Fwf) = \Bwf_T$. We first define $\lambda^{\bar\mu}$ on $\Fwf$ by $\lambda^{\bar\mu}(A)\coloneqq \Xi^{\bar\mu}(F)$ whenever $A=\bigcup_{t\in F}[t]$ for some $F\subseteq\Fr_n(T)$ and $n<\alt(T)$. We show that this function is well-defined and that it is $\sigma$-additive, as it is clear that $\lambda^{\bar\mu}(\emptyset) = \Xi^{\bar\mu}(\emptyset) = 0$ (the only $F$ representing $\emptyset$ is $\emptyset$). 

    To see that the map is well-defined, assume that $A=\bigcup_{t\in F}[t] = \bigcup_{t'\in F'}[t']$ where $F\subseteq\Fr_n(T)$, $F'\subseteq\Fr_{n'}(T)$ and $n,n'<\alt(T)$. Without loss of generality, consider the case $n\leq n'$. Then, for any $t\in F$, since $[t]\subseteq \bigcup_{t'\in F'}[t']$ we must have that $T_{\geq t}\cap F' = T_{\geq t}\cap \Fr_{n'}(T)$. Therefore, by~\eqref{b035-3eq}, $\Xi^{\bar\mu}(\{t\}) = \Xi(\set{t'\in F'}{t\subseteq t'})$. This implies that 
    $$\Xi^{\bar\mu}(F) = \sum_{t\in F}\Xi^{\bar\mu}(\{t\}) = \sum_{t\in F}\Xi^{\bar\mu}(\set{t'\in F'}{t\subseteq t'}) = \Xi^{\bar\mu}(F'),$$
    where the last equality hold because $F'=\set{t'\in F'}{\exists t\in F\ (t\subseteq t')}$, which follows by $\bigcup_{t\in F}[t] = \bigcup_{t'\in F'}[t']$.

    To show that $\lambda^{\bar\mu}$ is $\sigma$-additive in $\Fwf$, for $k<\omega$ let $n_k<\alt(T)$, $F_k\subseteq \Fr_{n_k}(T)$, and $A_k\coloneqq \bigcup_{t\in F_k}[t]$, and assume that $\set{A_k}{k<\omega}$ is pairwise disjoint and $\bigcup_{k<\omega}A_k = \bigcup_{t\in F}[t]$ for some $F\subseteq\Fr_n(T)$ and $n<\alt(T)$. We prove that $\Xi^{\bar\mu}(F)=\sum_{k<\omega}\Xi^{\bar\mu}(F_k)$. We may assume that $n_k\geq n$ for all $k<\omega$ because we can modify $n_k,F_k$ without affecting $A_k$ by using that, whenever $n_k\leq n'<\alt(T)$, $A_k=\bigcup_{t\in F_k}[t] = \bigcup_{t'\in F'}[t']$ where $F'$ is the set of nodes $t'\in\Fr_{n'}(T)$ above some $t\in F_k$. 

    Notice that any pair of nodes in each $F_k$ are pairwise incompatible, as well as nodes coming from different $F_k$ and $F_{k'}$. 
    Set $F_\omega\coloneqq \bigcup_{k<\omega}F_k$. Since $n_k\geq n$ for all $k<\omega$, we get that $F'_\omega\coloneqq F_\omega\cup(\Fr_n(T)\menos F)$ is a front of $T$. Hence,
    for $t\in F$, we must have $T_{\geq t}\cap F_\omega = T_{\geq t}\cap F'_\omega$. Therefore, by~\eqref{b035-3eq}, 
    \begin{align*}
        \Xi^{\bar\mu}(F) & = \sum_{t\in F}\Xi^{\bar\mu}(\set{t'\in F_\omega}{t\subseteq t'}) = \Xi^{\bar\mu}(F_\omega) = \sum_{k<\omega}\Xi^{\bar\mu}(F_k) = \sum_{k<\omega}\lambda^{\bar\mu}(A_k).
    \end{align*}
    This shows that $\lambda^{\bar\mu}$ is a measure on $\Fwf$. Therefore, by \autoref{b070} this is extended by a unique measure on $\Bwf_T = \sigma_{[T]}(\Fwf)$, which we still denote by $\lambda^{\bar\mu}$. This is a probability measure because $\lambda^{\bar\mu}([T]) = \Xi^{\bar\mu}(\{\la\,\ra\}) = 1$. To show the uniqueness, if $\lambda\in\BP$ satisfies that $\Xi^\lambda = \Xi^{\bar\mu}$, then $T_\lambda = T_{\Xi^\lambda} = T_{\bar\mu} = T$ and, for $n<\alt(T)$ and $F\subseteq \Fr_n(T)$, $\lambda\left(\bigcup_{t\in F}[t]\right) = \sum_{t\in F}\lambda([t]) = \sum_{t\in F}\Xi^{\bar\mu}(\{t\}) = \Xi^{\bar\mu}(F)$. Hence, $\lambda$ extends the measure we already defined on $\Fwf$, so $\lambda=\lambda^{\bar\mu}$ by the uniqueness of the extension.
    %
    % If $\alt(T) < \omega$, define $\lambda^{\bar{\mu}} \coloneqq \Xi_{h}^{\bar{\mu}}$, where $h \coloneqq \alt(T) - 1$. This satisfies the required condition by virtue of~\autoref{ppp12}. We deal with the infinite height case in~\autoref{5.2} (see in particular~\autoref{b060}). 
    %
    % We now prove the uniqueness. Assume that $\gamma$ is a measure on $\Bwf_T$ such that for all $t\in T$, $\gamma([t]) = \Xi_{|t|}^{\bar{\mu}}(\{t\})$. Define 
    % \[\Fwf \coloneqq \bigset{\bigcup_{t\in F}[t]}{F\subseteq \Lv_n(T),\ n<\omega}.\]
    % This is an algebra of sets over $[T]$ and every $C\in\Fwf$ is clopen in $[T]$, so $\Fwf\subseteq \Bwf_T$. It is clear that $[t]\in \Fwf$ for all $t\in T$, so $\sigma(\Fwf) = \Bwf_T$. On the other hand, if $n<\omega$ and $F\subseteq \Lv_n(T)$, then
    % \[\gamma\left(\bigcup_{t\in F}[t]\right) = \sum_{t\in F}\gamma([t]) = \sum_{t\in F} \Xi_{\vert t \vert}^{\bar{\mu}}(\{ t \}) =  \sum_{t\in F} \lambda^{\bar{\mu}}([t]) = \lambda^{\bar{\mu}}\left(\bigcup_{t\in F}[t]\right).\]
    % Therefore, $\nu \coloneqq \gamma\frestr\Fwf = \lambda^{\bar{\mu}}\frestr\Fwf$, which is a measure on $\Fwf$. Both $\gamma$ and $\lambda^{\bar{\mu}}$ are probability measures on $\sigma(\Fwf)$ extending $\nu$, so $\gamma= \lambda^{\bar{\mu}}$ by~\autoref{b070}.
\end{PROOF}

As a consequence, we have a connection between $\TP$ and $\BP$. 

\begin{definition}\label{ppp3}
    We define the function $\Lambda \colon \TP \to \BP$ by $\Lambda(\bar{\mu}) \coloneqq \lambda^{\bar{\mu}}$. 
\end{definition}

We list some properties of $\Lambda$ below. 

\begin{lemma}\label{pp129}
    \
    \begin{enumerate}[label=\normalfont(\alph*)]
        %\item\label{pp129.1.1} If $\langle T, \bar{\mu} \rangle$ is a probability tree, then $\Xi^{\bar{\mu}} = \Xi^{\lambda^{\bar{\mu}}}$. 
    
        \item\label{pp129.2}  $\Psi$ is a bijective function and, for $\Xi\in\IP$, $\Psi^{-1}(\Xi)=\lambda^{\bar\mu}$ where $\bar\mu\in\TP$ satisfies $\Xi^{\bar\mu}=\Xi$ (i.e.\ $\Psi^{-1}(\Xi)$ does not depend on this $\mu$). 

        \item\label{pp129.3} $\Lambda$ is surjective, and for $\bar{\mu}, \bar{\nu} \in \TP$, $\Lambda(\bar{\mu}) = \Lambda(\bar{\nu})$ iff\/ $\Pi(\bar{\mu}) = \Pi(\bar{\nu})$. 
    \end{enumerate}
\end{lemma}

\begin{PROOF}{\ref{pp129}}
    %\ref{pp129.1.1}: 
    %Assume that $\langle T, \bar{\mu} \rangle$ is a probability tree and let $t \in T$. By~\autoref{pp50} and~\autoref{ppp10},  $\Xi_{\vert t \vert}^{\lambda^{\bar{\mu}}}(\{ t \}) = \lambda^{\bar{\mu}}([t]) = \Xi_{\vert t  \vert}^{\bar{\mu}}(\{ t \})$. Thus, $\Xi^{\bar{\mu}} = \Xi^{\lambda^{\bar{\mu}}}$.
    % On the other hand, assume that $\alt(T_{\lambda^{\bar{\mu}}}) < \omega$. In this case, we proceed by induction on $\alt(T_{\bar{\mu}})$. The base step is clear because by the definition of $\Xi^{\lambda}$, for any $t \in \Lv_{h}(T)$,  we have that  $\Xi_{h}^{\Xi_{h}^{\bar{\mu}}} = \Xi_{h}^{\bar{\mu}}(\{ t \})$. Now, assume that for any $t \in \Lv_{n}(T)$, $\Xi_{n}^{\Xi_{h}^{\bar{\mu}}} = \Xi_{n}$ and let $t \in \Lv_{n-1}(T)$. Then, $ \Xi_{n-1}^{\Xi_{h}^{\bar{\mu}}}(\{ t \}) = \Xi_{n}^{\Xi_{h}^{\bar{\mu}}}(\suc_{T}(t)) = \Xi_{n}(\suc_{T}(t)) = \Xi_{n-1}(\{ t \}),$ whence follows the result.  

     \ref{pp129.2}: Let $\lambda, \lambda' \in \mathcal{BP}$ be such that $\Xi^{\lambda} = \Xi^{\lambda'}$. Then $T_\lambda = T_{\Xi^\lambda} = T_{\lambda'}$ and, by~\autoref{pp50}, for any $t \in T_{\lambda}$, $ \lambda([t]) = \Xi^{\lambda}(\{t\}) = \Xi^{\lambda'}(\{t\}) = \lambda'([t])$. This implies that $\lambda\frestr \Fwf_{T_\lambda}= \lambda'\frestr \Fwf_{T_\lambda}$ (see~\eqref{ppp10FT}), which implies $\lambda = \lambda'$ by \autoref{b070}. This proves that $\Psi$ is one-to-one. On the other hand, the surjectivity follows easily from~\autoref{ppp10}: for $\Xi \in \IP$, since $\Pi$ is surjective, we can find a $\bar{\mu} \in \TP$ such that $\Xi^{\bar{\mu}} = \Xi$, so $\Psi(\lambda^{\bar{\mu}}) = \Xi^{\lambda^{\bar{\mu}}} = \Xi^{\bar{\mu}} =\Xi$. %Thus, $\Psi$ is bijective. 
     This also shows that $\Psi^{-1}(\Xi) = \lambda^{\bar\mu}$.

    \ref{pp129.3}: By \ref{pp129.2}, $\Lambda = \Psi^{-1}\circ \Pi$, which is a composition of two surjective functions. Moreover, for  $\bar{\mu}, \bar{\nu} \in \TP$, $\Pi(\bar\mu) = \Pi(\bar\nu)$ iff $\Psi^{-1}(\Pi(\bar\mu))) = \Psi^{-1}(\Pi(\bar\nu))$, which means that $\Lambda(\bar\mu) = \Lambda(\bar\nu)$.
\end{PROOF}

Similar to the case of $\TP$ and $\IP$, we can define an equivalence relation on $\BP$. %this time leveraging the one already established on $\IP$.

\begin{definition}\label{pp134}
    \ 
    \begin{enumerate}
        \item Let $\BP_+$ be the class of $\lambda\in\BP$ such that $\lambda([t])>0$ for all $t\in T_\lambda$, i.e.\ the only measure zero open subset of $[T]$ is the empty set.
    
        \item For $\lambda\in\BP$, set $N_\lambda\coloneqq \set{t\in T_\lambda}{\lambda([t])=0}$, $T_{\lambda}^+\coloneqq T_\lambda\menos N_\lambda$, and $\lambda_+\coloneqq \lambda\frestr\calB([T^+_\lambda])$. Also set $N^*_\lambda \coloneqq \bigcup_{t\in N_\lambda}[t]$.
        
        \item We say that $\lambda, \lambda' \in \BP$ are \emph{positive equivalent}, denoted by $\lambda \equiv_{\BP} \lambda'$, iff $\lambda^+ = \lambda'_+$. Denote the equivalence class of $\lambda$ by $[\lambda]_\BP$.

        \item  Define the maps
            \begin{align*}
                \Psi_+ & \colon \BP_+ \to \IP_+ && \text{by $\Psi_+\coloneqq \Psi\frestr \BP_+$,}\\
                \Lambda_+ & \colon \TP_+ \to \BP_+ && \text{by $\Lambda_+\coloneqq \Lambda\frestr \TP_+$,}\\
                \pi^+_\BP & \colon \BP_+ \to \BP/_{\equiv_\BP} && \text{by $\pi^+_\BP(\lambda) \coloneqq [\lambda]_\BP$,}\\
                \hat{\Psi} &\colon \BP /_{\equiv_{\BP}} \to \IP /_{\equiv_{\IP}} && \text{by $\hat{\Psi}([\lambda]_{\BP}) \coloneqq [\Psi(\lambda)]_{\IP}$, and}\\ 
                \hat{\Lambda} & \colon \TP / _{\equiv_{\TP}} \to \BP /_{\equiv_{\BP}} && \text{by $\hat{\Lambda}([\bar{\mu}]_{\TP}) \coloneqq [\Lambda(\bar{\mu})]_{\BP}$.}
            \end{align*}
        These maps are well-defined thanks to the following result.
    \end{enumerate}
\end{definition}

\begin{fact}\label{pp140}
    Let $\bar\mu,\bar\nu\in\TP$ and $\lambda,\lambda'\in\BP$. Then: 
    \begin{enumerate}[label = \normalfont (\alph*)]
        \item\label{pp140a} $N_\lambda = N_{\Xi^\lambda}$ and $T^+_\lambda = T^+_{\Xi^\lambda}$.
        \item\label{pp140a1} $N^*_\lambda$ is the largest open measure zero subset of $[T_\lambda]$ and $[T^+_\lambda]= [T_\lambda]\menos N^*_\lambda$.
        \item\label{pp140b} $\lambda_+\in\BP_+$, $T_{\lambda_+} = T^+_{\lambda}$ and $\Xi^{\lambda_+} = \Xi^\lambda_+$.
        \item\label{pp140c} $\lambda\in\BP_+$ iff $T^+_\lambda = T_\lambda$ iff $\lambda_+=\lambda$. 
        \item\label{pp140d} $\lambda_+ = \lambda'_+$ iff $\Xi^\lambda \equiv_\IP \Xi^{\lambda'}$ (so $\hat{\Psi}$ is well-defined). 
        \item\label{pp140e} $\lambda^{\bar\mu_+} = \lambda^{\bar\mu}_+$ and $T^{\bar\mu_+} = T^+_{\bar\mu}$.
        \item\label{pp140f} $\bar\mu_+=\bar\nu_+$ iff $\lambda^{\bar\mu}_+ = \lambda^{\bar\nu}_+$ (so $\hat{\Lambda}$ is well-defined).
    \end{enumerate}
\end{fact}
\begin{PROOF}{\ref{pp140}}
    \ref{pp140a}: Immediate by the definition of $\Xi^\lambda$. 

    \ref{pp140a1}: Clear because $N^*_\lambda$ is composed by all the measure zero basic clopen sets. 

    \ref{pp140b}: By~\ref{pp140a}, $T^+_\lambda\in \calT$. It is clear that $\lambda_+$ is a measure on $\Bwf([T^+_\lambda])$ and $\lambda([T^+_\lambda]) = \lambda([T_\lambda])=1$ by~\ref{pp140a1}, so $\lambda_+\in\BP$ and $T_{\lambda_+}=T^+_\lambda$. Moreover, for $t\in T^+_\lambda$, $[t]_{T^+_\lambda} = [t]_{T_\lambda}\menos N^*_\lambda$, so $\lambda([t]_{T^+_\lambda}) = \lambda([t]_{T_\lambda})>0$. Thus $\lambda_+\in\BP_+$. 
    
    Finally, $T_{\Xi^{\lambda^+}} = T_{\lambda^+}= T^+_\lambda = T^+_{\Xi^\lambda}$ and, for $t\in T^+_\lambda$, $\Xi^{\lambda_+}(\{t\}) = \lambda([t]_{T^+_\lambda}) = \lambda([t]_{T_\lambda}) = \Xi(\{t\})$. Hence $\Xi^{\lambda_+} = \Xi^\lambda_+$.

    \ref{pp140c}: Clear by~\ref{pp140b} and the definitions.

    \ref{pp140d}: Clear by~\ref{pp140b} and because $\Psi$ is a bijection.

    \ref{pp140e}: By \autoref{ppp10}, $\Xi^{\lambda^{\bar\mu_+}}=\Xi^{\bar\mu_+}=\Xi^{\bar\mu}_+ = \Xi^{\lambda^{\bar\mu}}_+ =\Xi^{\lambda^{\bar\mu}_+}$, so the result follows because $\Psi$ is a bijection.

    \ref{pp140f}: The implication from left to right is clear by~\ref{pp140e}. For the converse, if $\lambda^{\bar\mu}_+ = \lambda^{\bar\nu}_+$ then, by~\ref{pp140b} and \autoref{ppp10}, $\Xi^{\bar\mu_+} = \Xi^{\lambda^{\bar\mu_+}} = \Xi^{\lambda^{\bar\nu_+}} = \Xi^{\bar\nu_+}$, so $\bar\mu_+=\bar\nu_+$ because $\Pi_+$ is a bijection.
\end{PROOF}

\begin{corollary}\label{pp129.1}
    The functions $\hat{\Psi} $ and $\hat{\Lambda}$ are bijections.
\end{corollary}

Finally, we can expand~\autoref{f45} by including $\BP$. 

\begin{theorem}\label{ppp14.T}
    Consider the diagram in~\autoref{f50}. Then:
   \begin{enumerate}[label=\normalfont(\alph*)]
        \item\label{ppp14.T.1} All the sub-diagrams in~\autoref{f50} are commutative.  As a consequence, any pair of paths that start at the same point and end at the same point produce the same function. Moreover: 
        \begin{enumerate}[label = \normalfont (\alph{enumi}.\arabic*)]
            \item any path from $\BP$ to $\BP_+$ is the map $\lambda\mapsto \lambda_+$, and
            \item any path from $\BP$ to $\BP/_{\equiv_\BP}$ is the map $\lambda\mapsto[\lambda]_\BP$.
        \end{enumerate}

        \item\label{ppp14.T.0}  $\Pi$, $\Lambda$, $\pi_{\TP}$, $\pi_{\IP}$, $\Delta_{+}$, $\varphi_{\IP}, \varphi_{\TP}$, $\gamma_{\GP}$ and $\gamma_{\GP}^+$ are surjective.

        \item\label{ppp14.T.2} $\Psi$, $\Delta$, $\Pi_{+}$, $\Psi_+$, $\Lambda_+$, $\hat{\Pi}$, $\hat{\Psi}$, $\hat{\Lambda}$,  $\pi_{\IP}^{+}$, $\pi_{\TP}^{+}$ and $\pi^+_\BP$ are bijective.
    \end{enumerate}
    \begin{figure}
    \centering
    \begin{tikzpicture}[>=stealth, node distance=3cm, scale=1, every node/.style={transform shape}]
    % Nodos con nombres correctos
    % Nodo nivel 1 (aumentado en Y)
    
    \node (BP) at (0, 6) {\footnotesize$\mathcal{BP}$};
    
    \node (TP) at (-4, 4) {\footnotesize$\TP$};   
    \node (IP) at (4, 4) {\footnotesize$\IP$};    

    \node (GP) at (0, 2) {\footnotesize $\mathcal{GP}$}; 

    \node (TP+) at (-2, 0) {\footnotesize$\TP_{+}$};  
    \node (IP+) at (2, 0) {\footnotesize$\IP_{+}$};  

    \node (BP+) at (0, -2) {\footnotesize$\BP_{+}$};

    \node (BPc) at (0, -4) {\footnotesize$\mathcal{BP}/_{\equiv_{\mathcal{BP}}}$};
    
    \node (TPc) at (-4, -6.5){\footnotesize$\TP/_{\equiv_{\TP}}$};
    \node (IPc) at (4, -6.5) {\footnotesize$\IP/_{\equiv_{\IP}}$}; 
 
    % Flechas con funciones
 % De nivel 1 a nodo derecho nivel 2
    \draw[->] (TP) -- (TPc) node[midway, left] {$\pi_{\TP}$};            % De nodo izquierdo nivel 2 a nodo izquierdo nivel 4
    \draw[->] (IP) -- (IPc) node[midway, right] {$\pi_{\IP}$};           % De nodo derecho nivel 2 a nodo derecho nivel 4
            % De nodo izquierdo nivel 4 a nodo nivel 5
       % De nodo derecho nivel 4 a nodo nivel 5
    \draw[->] (TP) -- (TP+) node[midway, right] {$\varphi_{\TP}$};            % De nodo izquierdo nivel 2 a nodo izquierdo nivel 3
    \draw[->] (IP) -- (IP+) node[midway, left] {$\varphi_{\IP}$};           % De nodo derecho nivel 2 a nodo derecho nivel 3
    \draw[->] (TP+) -- (TPc) node[midway, right] {$\pi_{\TP}^{+}$};            % De nodo izquierdo nivel 3 a nodo izquierdo nivel 4
    \draw[->] (IP+) -- (IPc) node[midway, left] {$\pi_{\IP}^{+}$};        % De nodo derecho nivel 3 a nodo derecho nivel 4
    
    % Flecha bidireccional entre G y H (nivel 3)
    \draw[->] (TP+) -- (IP+) node[midway, above] {$\Pi_{+}$};         % Conexión bidireccional entre los nodos del nivel 3
    
    % Flecha bidireccional entre los nodos del nivel 5 (F)
    \draw[->] (TPc) -- (IPc) node[midway, above] {$\hat{\Pi}$};         % Conexión bidireccional entre los nodos del nivel 5
    
    % Flecha de izquierda a derecha entre los nodos del nivel 2 (B y C)
    \draw[->] (TP) -- (IP) node[midway, above] {$\Pi$};          % Conexión de izquierda a derecha entre los nodos del nivel 2

    % Lo relacionado a GP 

  % Flechas con funciones
  \draw[<-] (TP) -- (GP) node[midway, right] {$\ \ \Delta$};  % De nodo izquierdo nivel 2 a nodo nivel 3.5
  \draw[<-] (IP) -- (GP) node[midway, left] {$\gamma_{\GP}$};  % De nodo derecho nivel 2 a nodo nivel 3.5
  \draw[->] (GP) -- (TP+) node[midway, left] {$\Delta_{+}$};  % De nodo nivel 3.5 a nodo nivel 3
  \draw[->] (GP) -- (IP+) node[midway, right] {$\gamma_{\GP}^+$};  % De nodo nivel 3.5 a nodo nivel 3

    \draw[<-] (BP) -- (TP) node[midway, above, sloped] {$\Lambda$}; 

    \draw[->] (BP) -- (IP) node[midway, above, sloped] {$\Psi$}; 

    \draw[->] (TPc) -- (BPc) node[midway, left] {$\hat{\Lambda} \ \ $};

    \draw[<-] (IPc) -- (BPc) node[midway, right] {$\ \ \hat{\Psi}$};  

    % Lo relacionado a BP+

    \draw[->] (TP+) -- (BP+) node[midway, left] {$\ \ \Lambda_{+}$};  

    \draw[->] (IP+) -- (BP+) node[midway, right] {$\ \Psi^{-1}_{+}$}; 

    \draw[->] (BP+) -- (BPc) node[midway, right] {$\pi_{\BP}^{+}$}; 

    \end{tikzpicture}
    
        \caption{Connections between the classes associated to $\TP$, $\IP$, $\GP$ and $\BP$.}
        \label{f50}
    \end{figure}
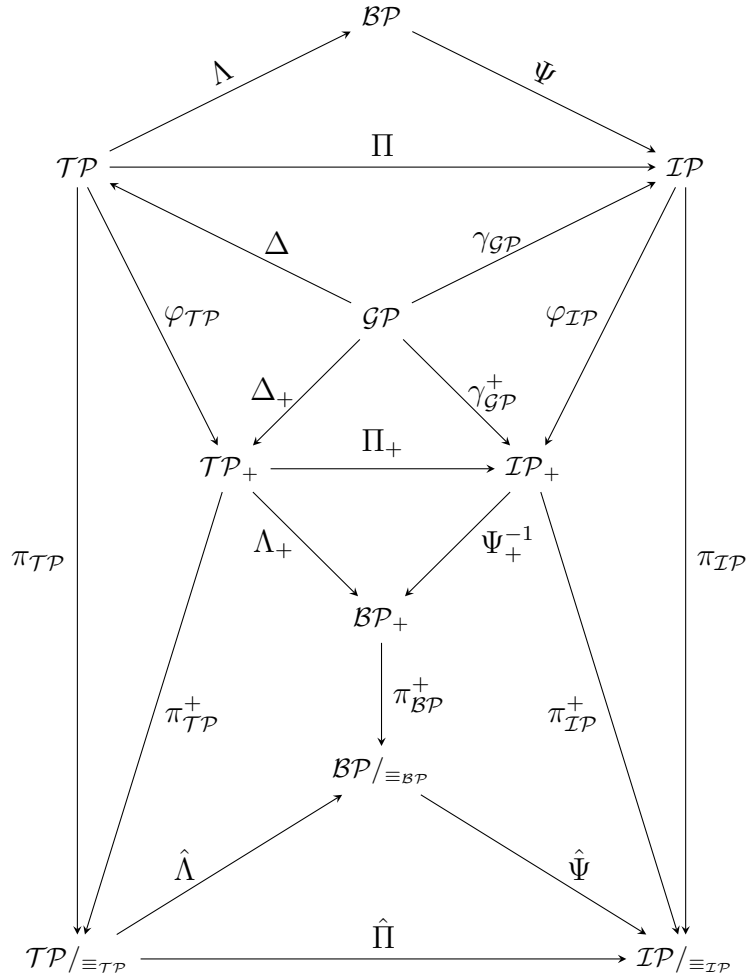

\end{theorem}

\begin{PROOF}{\ref{ppp14.T}} 
    The commutativity of the sub-diagrams of \autoref{f50} follows easily by \autoref{pp51} and~\ref{ppp10}. 
    
    We have~\ref{ppp14.T.0} and~\ref{ppp14.T.2} as a consequence of our previous efforts in~\autoref{pp51}, \autoref{pp129}, \autoref{pp140} and \autoref{pp129.1}.  
\end{PROOF}

\begin{corollary}\label{b135}
     Let $\langle T, \bar{\mu} \rangle$ be a probability tree. 
     \begin{enumerate}[label=\normalfont(\alph*)]
        \item\label{b135.1} 
       If $\lambda \in \BP$ and  $\Xi^{\bar\mu}=\Xi^\lambda$ then $\lambda = \lambda^{\bar{\mu}}$. As a consequence, $\lambda = \lambda^{\bar{\mu}^{\Xi^{\lambda}}}$ whenever $\bar{\mu}$ is positive.  
     
        \item\label{b135.2} If $ \bar{\mu}$ is  positive  then $\bar{\mu}^{\Xi^{\lambda^{\bar{\mu}}}} = \bar{\mu}$.
     \end{enumerate}
\end{corollary}

% \begin{PROOF}{\ref{pp1}}
%     On the one hand, we have that for any $s \in S \menos \max(S)$, $\tau_{s}$ is well-defined because since $\varphi \colon T \to S$ is an order-isomorphism,  $r \in \suc_{S}(s)$ iff $\varphi^{-1}(r) \in \suc_{T}(\varphi^{-1}(s))$.  On the other hand, for $s \in S \menos \max(S)$,
%     \begin{equation*}
%         \begin{split}
%             \tau_{s}(\suc_{S}(s)) & = \sum_{r \in \suc_{S}(s)} \tau_{s}(\{ r \}) = \sum_{r \in \suc_{S}(s)} \mu_{\varphi^{-1}(s)}(\{ \varphi^{-1}(r) \})\\ 
%             & = \sum_{t \in \suc_{T}(\varphi^{-1}(s))} \mu_{\varphi^{-1}(s)}(\{ t \}) = \mu_{\varphi^{-1}(s)}(\suc_{T}(\varphi^{-1}(s))) = 1,
%         \end{split}
%     \end{equation*} 
%     where the last equality holds because $\langle T, \bar{\mu} \rangle$ is a probability tree. 
% \end{PROOF}

\subsection{Relative expected value in probability trees}\label{4.2}\

In this subsection, we are going to introduce a notion of expected value on probability trees that we call \emph{relative expected value}. 

For all this section fix a probability tree $\langle T, \bar{\mu} \rangle$. Recall that, for any $t \in T$, $T_{\geq t}$ is the set of all nodes in $T$ above $t$. Notice that, if $0 \leq n < \hgt(T_{\geq t})$, then $\Lv_{n}(T_{\geq t}) \subseteq \Lv_{|t| + n}(T)$. Furthermore, $T_{\geq t}$ inherits a probability tree structure from $\langle T, \bar{\mu} \rangle$: 

\begin{lemma}\label{b038.9}
    Let $\langle T, \bar{\mu}^{T} \rangle$ be a probability tree.  Then, for any $t \in T$, $T_{\geq t}$ inherits probability tree structure from $T$ in a natural way, that is, $\langle T_{\geq t}, \bar{\mu}^{T_{\geq t}} \rangle$ is a probability tree, where $\bar{\mu}^{T_{\geq t}} \coloneqq \langle \mu_{s}^{T} \colon s \in T_{\geq t} \menos \max(T_{\geq t}) \rangle$. 
\end{lemma}

Since for any $s \in T_{\geq t}$, $\mu^{T_{\geq t}} =  \mu_{s}^{T}$, we can abuse of the notation and write ``$\langle T_{\geq t}, \bar{\mu} \rangle$'' instead of ``$\langle T, \bar{\mu}^{T_{\geq t}} \rangle$''. This will be widely applied in this subsection. 

We can now define the \emph{relative expected value} in probability trees. 

\begin{definition}\label{p60}
     Let $t\in T$, $m\coloneqq |t|$, $m \leq n < \alt(T)$ such that $T_{\geq t}\cap \Lv_n(T) = T_{\geq t}\cap \Fr_n(T)$, and let $X$ be a random variable on the probability space $\la\Fr_{n}(T),\Xi^{\bar\mu}\frestr \pts(\Fr_n(T))\ra.$ Then, we define: $$\E_{\Lv_{n}(T)}[X \colon s \rest m = t] \coloneqq \E_{\Lv_{n-m}(T_{\geq t})}[X \rest \Lv_{n-m}(T_{\geq t})],$$ and call it the \emph{relative expected value of $X$ with respect to $t$}. Here, $X \rest \Lv_{n-m}(T_{\geq t})$ is interpreted as a random variable on $\la\Lv_{n-m}(T_{\geq t}),\Xi^{\bar\mu^{T_{\geq t}}}\frestr\pts(\Lv_{n-m}(T_{\geq t}))\ra$. 
     When the context is clear, we simply write $E_{n}[X \colon s \rest m = t]$ or even $\E[X \colon s \rest m = t]$, instead of $\E_{\Lv_{n}(T)}[X \colon s \rest m = t].$
\end{definition}

Notice that the ``$s$'' above is a dummy variable, that is, the expected value of $X$ is calculated by varying $s$ over the nodes in $T$ at level $n$ that extend $t$. %Furthermore,  the relative expected value is well-defined because, by~\autoref{b038.9}, $T_{\geq t}$ inherits probability space structure from $T$ and, on the other hand, if $X$ is a random variable on $\Lv_{m + n}(T),$ then we have that $\Lv_{n}(T_{\geq t}) \subseteq \Lv_{m+n}(T),$ hence $X \rest \Lv_{n}(T_{\geq t})$ is a random variable on $\Lv_{n}(T_{\geq t}),$ so calculating its expected value makes sense. \Diego{Creo que esto es claro}
Since the relative expected value is defined in terms of the typical expected value, it is clear that it is linear, i.e.\ for $a,b\in\R$ and any random variables $X,Y$ on $\Fr_n(T)$,
$$\E_{n}[a X + b Y \colon s \rest m = t] = a\, \E_{n}[X \colon s \rest m = t] + b\, \E_{n}[Y \colon s \rest m = t].$$

% \begin{lemma}\label{p61}
%       Let $\langle T, \bar{\mu} \rangle$ be a probability tree $m < \omega$, and $t \in \Lv_{m}(T).$ Consider $0 < n < \omega,$ two random variables $X, Y$ on $\Lv_{m + n}(T)$ and $r, s \in \R.$ Then,  $$\E_{m + n}[r X + s Y \colon s \rest m = t] = r \E_{m + n}[X \colon s \rest m = t] + s \E_{m+n}[Y \colon s \rest m = t].$$
% \end{lemma}

The following result allows us to decompose the probability of the successors of $t$ at the level $m+n$ of $T,$ in terms of the probability at the level $n$ of $T_{\geq t}$: 

\begin{lemma}\label{p62}
    For $t\leq s$ in $T$, $\Xi^{\bar{\mu}}(\{ s \}) = \Xi^{\bar{\mu}^{T_{\geq t}}}_{n}(\{ s \})  \Xi^{\bar{\mu}}(\{ t \}).$
\end{lemma}

% \begin{PROOF}[\textbf{Proof}]{\ref{p62}}
%     Since by~\autoref{b038.9}, $\langle T_{\geq t}, \bar{\mu} \rangle$ is a probability tree, by~\autoref{b038.9}, $\Lv_{n}(T_{\geq t})$ is a probability space, that is, considering $\Xi^{\bar{\mu} \rest T_{\geq t}}_{n}$ makes sense. Now, $s \in \Lv_{k}(T_{\geq t})$  implies $s \in \Lv_{m +k}(T),$ and therefore, 
%         \begin{equation*}
%             \begin{split}
%                 \Xi_{m+n}^{\bar{\mu}}(\{ s \}) & = \prod_{0 \leq i \leq m + n -1} \mu_{s \rest i}(\{ s \rest (i +1) \})\\ 
%                 %line 2
%                 & = \left( \prod_{0 \leq i \leq m - 1} \mu_{s \rest i} ( \{ s \rest (i + 1) \} ) \right)  \left( \prod_{m \leq i \leq m + n -1} \mu_{s \rest i}(\{ s \rest (i +1) \}) \right) \\
%                 % line 3
%                 & = \left ( \prod_{0 \leq i \leq m -1} \mu_{t \rest i} (\{ t \rest ( i+ 1) \}) \right)   \left( \prod_{0 \leq j \leq n -1} \mu_{s \rest (m + j)} (\{ s \rest (m + j + 1) \}) \right) \\
%                 % line 4
%                 & = \Xi_{m}^{\bar{\mu}}(\{ t \})  \Xi_{n}^{\bar{\mu} \rest T_{\geq t}}(\{ s \}). 
%             \end{split}
%         \end{equation*} 

%     Thus, $\Xi_{m+n}^{\bar{\mu}}(\{ t \}) = \Xi^{\bar{\mu} \rest T_{\geq t}}_{n}(\{ s \})  \Xi_{m}^{\bar{\mu}}(\{ t \})$.
% \end{PROOF}

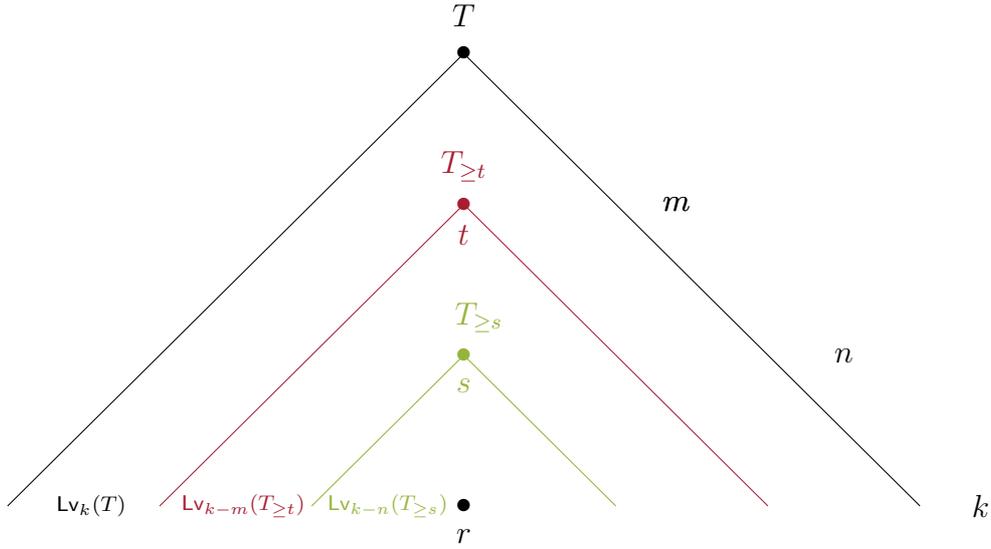
\begin{figure}
        \centering
        \begin{tikzpicture}
            % nodos

            \node at (6, 6) {$\bullet$}; 
            \node[redun] at (6, 4) {$\bullet$}; 
            \node[greenun] at (6, 2) {$\bullet$}; 
            \node at (6, 0) {$\bullet$}; 

            %triangulo grande 

            \draw (6, 6) -- (0, 0);
            \draw (6, 6) -- (12, 0);

            % triangulo mediano 

            \draw[redun] (6, 4) -- (2, 0);
            \draw[redun] (6, 4) -- (10, 0);

            % triangulo pequeño 

            \draw[greenun] (6, 2) -- (4, 0);
            \draw[greenun] (6, 2) -- (8, 0);

            % alturas 

            \node at (8.8, 4) {$m$};
            \node at (11, 2) {$n$};
            \node at (12.8, 0) {$k$};

            % nodos 

            \node at (6, 6.5) {$T$};
            \node[redun] at (6, 4.5){$T_{\geq t}$};
            \node[greenun] at (6.2, 2.5) {$T_{\geq s}$};
            \node at (8.8, 4) {$m$};         

            % niveles 

            \node at (1.1, 0) {\tiny $\Lv_{k}(T)$};
            
            \node[redun] at (3.1, 0) {\tiny $\Lv_{k-m}(T_{\geq t})$};
            
            \node[greenun] at (5, 0) {\tiny $\Lv_{k-n}(T_{\geq s})$}; 

            % funciones 

            \node at (6, -0.4) {$r$};
            \node[greenun] at (6, 1.6) {$s$};
            \node[redun] at (6, 3.6) {$t$};
        \end{tikzpicture}
        
        \caption{The situation in \autoref{p64}.}
        \label{f40}
    \end{figure}

The relative expected value can be calculated as a composition of relative expected values at intermediate levels, as follows (see \autoref{f40}).
%Now, we can show that, to calculate a relative expected value, we can do \emph{intermediate steps} (see \autoref{f40}), that is: 

\begin{theorem}\label{p64}
    Let $m\leq n\leq k<\omega$, $t \in \Lv_{m}(T),$ and assume that $T_{\geq t}\cap \Fr_k(T) = T_{\geq t}\cap \Lv_k(T)$. If $X$ is a random variable on $\Fr_{k}(T)$, then  
    $$\E_{k}[X \colon r \rest m = t] = \E_{n}\big[\E_{k}[ X \colon r \rest n = s ] \colon s \rest m = t \big].$$ 
\end{theorem}

\begin{PROOF}[\textbf{Proof}]{\ref{p64}}
    Let $E \coloneqq \E_{n}\big[ \E_{k}[X \colon r  \rest \ell = s ] \colon s  \rest m = t\big]$. By~\autoref{p60}, we have that: 
    \begin{equation*}
        \begin{split}
            E &\coloneqq    \sum_{s \in \Lv_{n-m}(T_{\geq t})} \E_{k}[X \colon r \rest n = s]  \Xi^{\bar{\mu}^{T_{\geq t}}}(\{ s \})  = \sum_{s \in \Lv_{n-m}(T_{\geq t})} \E[X \rest \Lv_{k-n}(T_{\geq s})]  \Xi^{\bar{\mu}^{T_{\geq t}}}(\{ s \}) \\
            % line 3
            & = \sum_{s \in \Lv_{n-m}( T_{\geq t})}  \left( \sum_{r \in \Lv_{k-n}(T_{\geq s})} X(r)  \Xi^{\bar{\mu}^{T_{\geq s}}}(\{ r \}) \right) \Xi^{\bar{\mu}^{T_{\geq t}}}(\{ s \})\\ 
            & = \sum_{s \in \Lv_{n-m}(T_{\geq t})} \sum_{r \in \Lv_{k-n}(T_{\geq s})} X(r)  \Xi^{\bar{\mu}^{T_{\geq t}}}(\{ r \})\\
            & =  \sum_{r \in \Lv_{k-m}(T_{\geq t})} X(r)  \Xi^{\bar{\mu}^{T_{\geq t}}}(\{ r \}) = \E[X \rest \Lv_{k-m}(T_{\geq t})]= \E_{k}[X \colon r \rest m = t],
        \end{split}
    \end{equation*}
    % \begin{equation*}
    %     \begin{split}
    %         \E_{\ell}[\, \E_{j}[X \colon r  \rest \ell = s ] \colon s  \rest m = t] & =  \!  \sum_{s \in \Lv_{n}(T_{\geq t})} \E_{j}[X \colon r \rest \ell = s] \cdot \mu_{T_{\geq t}}^{n}(s) \\
    %         % line 2
    %         & = \sum_{s \in \Lv_{n}(T_{\geq t})} \E_{\Lv_{k}(T_{\geq s})}[X \rest \Lv_{k}(T_{\geq s})] \cdot  \mu_{T_{\geq t}}^{n}(s) \\
    %         % line 3
    %         & = \sum_{s \in \Lv_{n}( T_{\geq t})}  \left( \sum_{r \in \Lv_{k}(T_{\geq s})} X(r) \cdot \Xi^{\bar{\mu} \rest T_{\geq s}}_{n-k}(r) \right) \mu_{T_{\geq t}}^{n}(s)\\
    %         % line 4
    %     \end{split}
    % \end{equation*}
    % \begin{equation*}
    %     \begin{split}
    %     \ \ \ \ \ \ \ \ \ \ \ \ \ \ \ \ \ \ \ \ \   & = \sum_{s \in T_{\geq t}} \left( \sum_{r \in T_{\geq s}} X(r) \cdot \Xi^{\bar{\mu} \rest T_{\geq t}}_{k}(r) \right)\\
    %     & =  \sum_{r \in \Lv_{k}(T_{\geq t})} X(r) \cdot \Xi^{\bar{\mu} \rest T_{\geq t}}_{k}(r)\\
    %     % line 5
    %     & = \E_{\Lv_{k}(T_{\geq t})}[X \rest \Lv_{k}(T_{\geq t})]\\
    %     & = \E_{m +k}[X \colon r \rest m = t],
    %     \end{split}
    % \end{equation*}
    where $ \Xi^{\bar{\mu}^{T_{\geq t}}}(\{ r \}) = \Xi^{\bar{\mu}^{T_{\geq s}}}(\{ r \}) \cdot \Xi^{\bar{\mu}^{T_{\geq t}}}(\{ s \})$ by \autoref{p62}.  
\end{PROOF}

Finally, as a consequence (when $m=0$), we can express the expected value of $X$ in terms of the relative expected value:

\begin{corollary}\label{p66}
    If $n\leq k<\alt(T)$, $\Lv_k(T)$ is a front of $T$ and $X$ is a random variable on $\Lv_{k}(T),$ then $\E[X] = \E_{n}\big[ \E_{k}[X \colon r \rest n = s]\big].$
\end{corollary}

% \begin{PROOF}[\textbf{Proof}]{\ref{p66}}
%     First, notice that $\E_{n}[X] = \E_{\Lv_{n}(T_{ \geq \langle \, \rangle)}}[X \rest T_{\geq \langle \, \rangle}] = \E_{n}[X \colon r \rest 0 = \langle \, \rangle].$ Now,  by~\autoref{p64}, we get: 
%     \begin{equation*}
%         \begin{split}
%             \E_{m}[X] & = \E_{m}[X \colon r \rest 0 = \langle \, \rangle]  = \E_{n}[\, \E_{m}[ X \colon r \rest n = s ] \colon s \rest 0 = \langle \, \rangle] = \E_{n}[\, \E_{m}[X \colon r \rest n = s]].
%         \end{split}\qedhere
%     \end{equation*}
% \end{PROOF}

So far we discussed the relative expected values of random variables on $\Fr_n(T)$. However, we can extend this notion to random variables on any front. We first fix some terminology.

\begin{definition}
    Let $t\in T$ and $A,A'$ fronts of $T$.
    \begin{enumerate}[label = \normalfont (\arabic*)]
        \item The node \emph{$t$ is below $A$} if $t\leq s$ for some $s\in A$.
        \item The front \emph{$A$ is below $A'$} if any node of $A$ is below $A'$. 
    \end{enumerate}
\end{definition}

For instance, $\Fr_m(T)$ is a front below $\Fr_n(T)$ whenever $m\leq n<\alt(T)$.

\begin{definition}
    Let $A$ be a front of $T$, $t\in T$ below $A$, and let $X$ be a random variable on $\la A,\Xi\frestr \pts(A)\ra$ (recall from \autoref{b035-3} that $\Xi\frestr\pts(A)$ is a probability measure). Define the \emph{relative expected value of $X$ with respect to $t$} as 
    \[\E_{A}[X\colon s\frestr |t| = t] \coloneqq E_{A\cap T_{\geq t}}[X\frestr A\cap T_{\geq t}].\]
    Notice that $A\cap T_{\geq t}$ is a front of $T_{\geq t}$. Also, $X\frestr A\cap T_{\geq t}$ is interpreted as a random variable on the probability space $\la  A\cap T_{\geq t}, \Xi^{\bar\mu^{T_{\geq t}}}\frestr\pts(A\cap T_{\geq t})\ra$. 
\end{definition}

We can generalize \autoref{p64} for fronts. We omit the proof, as it is very similar.

\begin{theorem}
    Assume that $A,A'$ are fronts, $t\in T$ is below $A$, and $A$ is below $A'$. If $X$ is a random variable on $A'$, then
    \[\E_{A'}[X\colon r\frestr |t| = t] = \E_{A}\big[ \E_{A'}[X\colon r\frestr |s| = s] \colon s\frestr |t| = t \big]. \]
    In particular, $\E_{A'}[X] = \E_A\big[ \E_{A'}[X \colon r\frestr|s| = s]\big]$.
\end{theorem}

%%%%%%%%%%%%%%%%%%%%%%%%%%%%%%%%%%%%%%%%%%%%%%%%%%%%%%%%%%%%%%%%%%%%%%%%%%%%%%%%%%%%%%%%%%%%%%%%%%%%%%%%%%%%%%%%%%%%%%%%%%%%%%%%%

\section{Bounding cumulative dependent Bernoulli distributions}\label{5.1}

As mentioned in the introduction it is well known that, by adding finite many independent and identically distributed random variables with Bernoulli distribution, we obtain a random variable with the binomial distribution. However, there are cases where we must deal with dependent random variables with Bernoulli distribution which may not be identically distributed.\footnote{See, for example, the proof of \cite[Main Lemma 4.3.17]{uribethesis} and the proof of \cite[Main Lemma~7.17]{CMU}.} 
In the following theorem, we show how the cumulative distribution of the number of successes of these random variables can be bounded by the cumulative distribution of the binomial distribution, which proves \autoref{thmA}. Here, success corresponds to $0$ and failure to $1$.

For a natural number $n$, $n$-many dependent trials with Bernoulli distribution can be understood as a probability tree $\la T,\bar\mu\ra$ where $T= {}^{\leq n}2$ is the \emph{complete binary tree of height $n+1$}, i.e.\ $t\in T$ iff $t$ is a sequence of length ${\leq}\, n$ composed by $0$'s and $1$'s (including the empty sequence). Any $t\in T$ of length $k\leq n$ represents a sequence of success and failures of the first $k$ Bernoulli tests and, whenever $k<n$, $\suc_T(t)$ is Bernoulli distributed with probability of success $\mu_t(\{\concat{t}{\la 0\ra}\})$, which clearly depends on $t$, i.e.\ on the previous trials. 

The random variable expressing the success at the $k$-th trial for $k<n$ is $X_k\colon \Lv_{k+1}(T)\to \{0,1\}$, defined by $X_k(t')\coloneqq t'(k)$, so success is attained when $X_k(t')=0$, i.e.\ $t'(k)=0$. Therefore, the probability of success is 
\[\Pr[X_k=0]=\sum_{t\in \Lv_k(T)}\Xi^{\bar\mu}(\{t\})\mu_t(\{\concat{t}{\la0\ra}\}) = \sum_{t\in \Lv_k(T)} \Xi^{\bar\mu}(\{\concat{t}{\la0\ra}\}).\]

\begin{theorem}\label{p70}
    Let $n < \omega$, $T  \coloneqq {}^{\leq n} 2$, and assume that $\langle T, \bar{\mu} \rangle$ is a probability tree. Define $Y \colon \Lv_{n}(T) \to \R$ as the random variable measuring the number of successes after $n$ trials, that is, for any $t \in \Lv_{n}(T),$ 
    $$ Y(t) \coloneqq \vert \{ k < n \colon t(n) = 0 \} \vert.$$
    
    Assume that there exists some $p \in [0, 1]$ such that, for any $t \in T\menos\max(T)$, $p \leq p_{t} \coloneqq \mu_{t}(\{ \concat{t}{\langle 0 \rangle} \})$. Then, for all $z \in \R,$ 
    $$\Pr[Y \leq z] \leq \Pr[\B_{n, p} \leq z],$$ 
    where $\B_{n, p}$ denotes the binomail distribution of $n$ trials, each with probability of success $p$. 
    
    %where $\Omega_{n^{\ast}} \coloneqq \{ i < \omega  \colon i \leq n^\ast{} \},$ $\calA_{n^{\ast}} \coloneqq \calP(\Omega_{n^{\ast}})$ and  $\mu_{n^{\ast}} \colon  \calA_{n^{\ast}} \to \R$ such that, for any $i \leq n^{\ast},$ $ \mu_{n^{\ast}}(i) \coloneqq \binom{n^{\ast}}{i} p^{i} (1 -p)^{n - i},$ and $\B_{n, p} \colon \Omega_{n} \to \R$ the identity map. This corresponds to the binomial distribution with parameters $n^{\ast}$ and $p$. 
\end{theorem}

\begin{PROOF}{\ref{p70}}
    %For any $n < n^{\ast},$ define the random variable $X_{n}$ on $\Lv_{n^{\ast}}(T)$ such that, for every $t \in \Lv_{n^{\ast}}(T),$ $X_{n}(t) \coloneqq 1 - t(n)$. Thereby, clearly for any $t \in T$ with $\suc_{T}(t) \neq \emptyset$,  $ X_{\mathrm{\vert t \vert}} \rest  \suc_{T}(t) \sim \Bernoulli(p_{t})$ and also, we have that $Y = \sum_{n < n^{\ast}} X_{n}$, which is the sum of random variables that we referred in the comment previous  to the theorem.
    
    For any $x \in [0, 1]$ and $d \in \{ 0, 1 \},$ define:
    $$I_{x}^{d} \coloneqq 
    \begin{cases}
            [0, x), &  d = 0, \\[1ex]
            [x, 1], & d = 1.
    \end{cases}$$
   For $t \in \Lv_{n}(T),$ let $C_{t}^{\bullet} \coloneqq \prod_{k < n} I_{p_{t \rest n}}^{t(n)},$ and notice that its volume is $\vol(C_{t}^{\bullet}) = \Xi^{\bar{\mu}}(\{ t \}).$  It is easy to show that $\{ C_{t}^{\bullet} \colon t \in \Lv_{n}(T) \}$ is a partition of the $n$-dimensional unitary cube ${}^n [0, 1]$. 
   % it is clear that it is a collection of pairwise disjoint sets, so let $x = \langle x_{n} \colon n < n^{\ast} \rangle \in [0, 1]^{n^{\ast}}$. Define, by recursion on $m < n^{\ast}$ a sequence $s_{x} \in \Lv_{n^{\ast}}(T),$  as follows: 
   %  $$s_{x}(0) \coloneqq \left\{ \begin{array}{lcc}
   %          0, &   \rm{if}  & x_{0} \in [0, p_{\langle \, \rangle}), \\
   %           \\ 1, &  \rm{if}  & x_{0} \in [ p_{\langle \, \rangle}, 1].
   %           \end{array}
   % \right.$$

   % Assume that we have constructed $s_{x}(m)$ $m + 1 < n^{\ast}.$ Then: 

   %  $$s_{x}(m+1) \coloneqq \left\{ \begin{array}{lcc}
   %          0, &   \rm{if}  & x_{m+1} \in [0, \, p_{s_{x} \rest m+1}), \\
   %           \\ 1, &  \rm{if}  & x_{m} \in [ p_{s_{x} \rest m+1},  1].
   %           \end{array}
   % \right.$$
   
   % It is clear that $s_{x} \in \Lv_{n^{\ast}}(T)$. Now, let $n < n^{\ast}$ and consider two possible cases: 

   % \begin{enumerate}
   %     \item[$\bullet$] When $x_{n} \in [0, p_{s_{x} \rest n})$: in this case, by definition of $s_{x},$ we have that $s_{x}(n) = 0$ and therefore, $ x_{n} \in [0, p_{s_{x} \rest n}) = I_{p_{s_{x} \rest n}}^{s_{x}(n)}$.

   %     \item[$\bullet$] When $x_{n} \in [p_{s_{x} \rest n}, 1]$: in this case, $s_{x}(n) = 1,$ hence $x_{n} \in [p_{s_{x} \rest n}, 1] = I_{p_{s_{x} \rest n}}^{s_{x}(n)}$. 
   % \end{enumerate}

   %  Therefore, in any case, $x_{n} \in I_{p_{s_{x} \rest n}}^{s_{x}}$. Thus, $x \in C_{s_{x}}^{\bullet}$.  

   For $z\in\R$ let $C^\bullet(z) \coloneqq \bigcup\set{C^\bullet_t}{Y(t)\leq z,\ t\in\Lv_n(T)}$. Thus
   \begin{equation}\label{p70-3}
        \Pr[Y\leq z] = \sum\set{\vol(C^\bullet_t)}{Y(t)\leq z,\ t\in \Lv_n(T)} = \vol(C^\bullet(z)).
   \end{equation}

   %  Now, let $z \in \R.$ By the definition of $C_{t}^{\bullet}$ and the constructions of $s_{x},$ we have that, for any $x = \langle x_{n} \colon n < n^{\ast} \rangle \in [0, 1]^{n^{\ast}}$, 
   %  $ \vert \{ n < n^{\ast} \colon x_{n} \leq p_{s_{x} \rest n} \} \vert \leq z \Leftrightarrow \sum_{n < n^{\ast}} X_{n}(s_{x}) \leq z,$ 
   %  whence follows that:
   %  \begin{equation}\label{p70-3}
   %     \text{$ \bigcup  \left \{ C_{t}^{\bullet} \colon t \in \Lv_{n^{\ast}}(T), \sum_{n < n^{\ast}} X_{n}(t) \leq  z \right \} =  \left \{ x \in [0, 1]^{n^{\ast}} \colon \vert \{ n < n^{\ast} \colon x_{n} \leq p_{s_{x} \rest n} \} \vert \leq  z \right \}$}.
   % \end{equation}

   % On the other hand,  
   % \begin{equation}\label{p70-5}
   %     \begin{split}
   %          % line 1
   %         \Xi_{n^{\ast}}^{\bar{\mu}} [Y \leq z]_{} & = \sum \left \{ \Xi_{n^{\ast}}^{\bar{\mu}}(t) \colon t \in \Lv_{n^{\ast}}(T) \conj \sum_{n < n^{\ast}} X_{n}(t) \leq z \right \} \\
   %         % % line 2
   %         % & = \sum \{ \vol(C_{t}^{\bullet}) \colon t \in \Lv_{n^{\ast}}(T) \conj \sum_{n < n^{\ast}} X_{n}(t) \leq z \}\\
   %         % line 3
   %         & = \vol \left( \bigcup \{  C_{t}^{\bullet} \colon t \in \Lv_{n^{\ast}}(T) \conj  \sum_{n < n^{\ast}} X_{n}(t) \leq z \} \right). 
   %     \end{split}
   % \end{equation}

    On the other hand, define the polyhedron $C_{t} \coloneqq \prod_{k < n} I_{p}^{t(n)}$ for $t\in \Lv_n(T)$. We can use this to express the cumulative binomial distribution because $\set{C_t}{t\in T}$ is a partition of ${}^n [0,1]$ and, by setting $C(z)\coloneqq \bigcup\set{C^\bullet_t}{Y(t)\leq z,\ t\in\Lv_n(T)}$, we obtain
    \begin{equation}\label{p70-e7}
        \Pr[\B_{n,p}\leq z] = \sum\set{\vol(C_t)}{Y(t)\leq z,\ t\in T} = \vol(C(z)).
    \end{equation}
    
    % Now, we are going to define another cube similar to $C_{t}^{\bullet}$ but, in order to be able to compare with a binomial distribution, we are going to define it with constant probability given by $p$:  for any $t \in \Lv_{n^{\ast}}(T),$ define $$ C_{t} \coloneqq \prod_{n < n^{\ast}} I_{p}^{t(n)}.$$  
    
    % In an analogous way for $C_{t}^{\bullet}$, we have that: 
    % \begin{equation}\label{p70-e9}
    %     \text{$ \bigcup \left \{ C_{t} \colon t \in  \Lv_{n^{\ast}}(T),  \sum_{n < n^{\ast}} X_{n}(t) \leq z \right \} = \{ x \in [0, 1]^{n^{\ast}} \colon \vert \{ n < n^{\ast} \colon x_{n} \leq p \} \vert \leq z \}$.}
    % \end{equation} 

    % Also, 
    % \begin{equation}\label{p70-e7}
    %     \text{$\mu_{n^{\ast}}[\B_{n^{\ast}, p} \leq z] = \vol \left( \bigcup \left \{ C_{t} \colon t \in  \Lv_{n^{\ast}}(T),  \sum_{n < n^{\ast}} t(n) \leq z \right \} \right)$.}
    % \end{equation}

    To conclude the proof, thanks to \eqref{p70-3} and~\eqref{p70-e7}, it is enough to show that $C^\bullet(z)\subseteq C(z)$ for all $z\in\R$. Let $\bar x = \seq{x_k}{k<n}\in C^\bullet(Z)$, so there is a unique $t\in\Lv_n(T)$ such that $\bar x\in C^\bullet_t$ and $Y(t) \leq z$. Likewise, there is a unique $t'\in\Lv_n(T)$ such that $\bar x\in C_{t'}$. For $k<n$, if $t'(k)=0$ then $x_k\in I^0_p =[0,p)\subseteq [0,p_{t\frestr k})=I^0_{p_{t\frestr k}}$ (the contention because $p_s\leq p$ for all $s\in T\menos\max(T)$), so $t(k)=0$. This shows that $Y(t')\leq Y(t)$, so $Y(t')\leq z$ and we can conclude that $\bar x\in C(z)$.
    %
   %  Since for any $t \in T, p \leq p_{t},$ we have that $\{ n < n^{\ast} \colon n_{n} \leq p \} \subseteq \{ n < n^{\ast} \colon x_{n} \leq p_{s_{x} \rest n}  \}$ and therefore, 
   % \begin{equation}\label{p70-e10}
   %     \text{$\{ x \in [0, 1]^{n^{\ast}} \colon \vert \{n < n^{\ast} \colon  n_{n} \leq p_{s_{x} \rest n} \} \vert \leq z \} \subseteq \{ x \in [0, 1]^{n^{\ast}} \colon \vert \{ n < n^{\ast} \colon x_{n} \leq p \} \vert \leq z \}$.}
   % \end{equation}
   %
   % Finally, by (\ref{p70-5}), (\ref{p70-3}), (\ref{p70-e10}), (\ref{p70-e9}) and (\ref{p70-e7}), in that order, we have:
   % \begin{equation*}
   %     \begin{split}
   %         \Xi_{n}^{\bar{\mu}} [Y \leq z] & = \vol \left( \bigcup \left \{  C_{t}^{\bullet} \colon t \in \Lv_{n^{\ast}}(T),  \sum_{n < n^{\ast}} X_{n}(t) \leq z \right \} \right) \\ 
   %         % line 2
   %         & = \vol( \left \{ x \in [0, 1]^{n^{\ast}} \colon \vert \{ n < n^{\ast} \colon x_{n} \leq p_{s_{x} \rest n} \} \vert \leq  z \right \} )\\
   %         % line 3
   %         & \leq \vol (\{ x \in [0, 1]^{n^{\ast}} \colon \vert \{ n < n^{\ast} \colon x_{n} \leq p \} \vert \leq z \})\\
   %         % line 4
   %         & = \vol \left(  \bigcup \left \{ C_{t} \colon t \in  \Lv_{n^{\ast}}(T) \conj  \sum_{n < n^{\ast}} t (n) \leq z \right \} \right)\\ 
   %         & = \mu_{n^{\ast}}[\B_{n^{\ast}, p} \leq z].  \qedhere
   %     \end{split}
   % \end{equation*}
\end{PROOF}

\section{Probability trees and the real line}\label{5.2}

In this section, we explore the connection between probability trees and the real line. To this end, we start with an alternative proof of~\autoref{ppp10}, specifically, given a probability tree $\langle T, \bar{\mu} \rangle$, we construct a measure $\lambda^{\bar{\mu}} \in \BP$ such that $\lambda^{\bar{\mu}}([t]_{T}) = \Xi^{\bar{\mu}}({ t })$ for any $t \in T$. This construction also shows a connection between the measure space $\la[T],\Bwf_T,\lambda^{\bar\mu}\ra$ and the Lebesgue measure space of $[0,1]$, reflected by a map that, under certain conditions, preserves measure (see~\autoref{b0100} and~\autoref{b095}).

For all this section, we assume that $\la T,\bar\mu\ra$ is a probability tree. Since $T$ is countable, $\Bwf_T$ is the $\sigma$-algebra generated by $\set{[t]}{t\in T}$.

\begin{definition}\label{b037}
    It is easy to show that $T$ is isomorphic with a tree $T_*$ such that, for any $t\in T_*\menos\max T_*$, there is some $\alpha_t\leq\omega$ such that $\suc_{T_*}(t) = \set{\concat{t}{\la k\ra}}{k<\alpha_t}$. 
    Such a $T^*$ is called a \emph{representation of\/ $T$}.
\end{definition}

To define the measure $\lambda^{\bar{\mu}}$, we fix a representation $T_*$ of $T$ as in \autoref{b037} and, without loss of generality, assume $T=T_*$. Our construction is motivated by known connections between the Cantor space, the Baire space, and $[0,1]$ as in e.g.~\cite[Ch.~VII, \S3]{LevyBasic}.
 
\begin{definition}\label{b040}
  Define a sequence $\bar I^{\bar\mu} \coloneqq \seq{I_t}{t\in T}$ of closed intervals $I_t= I^{\bar\mu}_t=[a_t,b_t]$ by recursion on $n=\hgt_T(t)$ as follows.
  \begin{itemize}
    \item $I_{\la\ \ra}  \coloneqq  [0,1]$, that is,  $a_{\la\ \ra} \coloneqq 0$ and $b_{\la\ \ra} \coloneqq 1$.
    
    \item Having defined the interval $I_t$, when $t\notin \max T$ let $\seq{I_{\concat{t}{\la k\ra}}}{k<\alpha_t}$ be the collection of consecutive closed intervals contained in $I_t$ where each $I_{\concat{t}{\la k\ra}}$ has length $\Lb(I_{t}) \mu_{t}(\{ \concat{t}{\langle k \rangle} \})$, that is:   $a_{\concat{t}{\la 0\ra}} \coloneqq a_t$, $b_{\concat{t}{\la k\ra}}\coloneqq a_{\concat{t}{\la k\ra}}+\Lb(I_t)\mu_t(\{\concat{t}{\la k\ra}\})$ whenever $k<\alpha_t$, and $a_{\concat{t}{\la k+1\ra}} \coloneqq b_{\concat{t}{\la k\ra}}$ whenever $k + 1 < \alpha_{t}$.
  \end{itemize}
  
  Denote $Q_{\bar\mu}  \coloneqq  \set{a_t, b_t}{t\in T}$. %i.e.\ $Q_{\bar\mu}$ is the set of endpoints of all the intervals $I_t$ with $t\in T$. 
  Then $|Q_{\bar\mu}|\leq 2|T|\leq \aleph_0$, so $Q_{\bar\mu}$ is countable.
  % Note that $I_t\subseteq I_s$ when $s\subseteq t$ in $T$, so $a_s\leq a_t\leq b_t\leq a_s$. Then, for $x\in[T]$, $\seq{a_{x\frestr n}}{n<\omega}$ is a monotone increasing sequence and $\seq{b_{x\frestr n}}{n<\omega}$ is a monotone decreasing sequence in $[0,1]$, so their limits exist. This allows to define two functions $f^-_{\bar\mu}\colon [T]\to [0,1]$ and $f^+_{\bar\mu}\colon [T]\to [0,1]$ by $\displaystyle f^-_{\bar\mu}(x) \coloneqq \lim_{n\to\infty} a_{x\frestr n}$ and $\displaystyle f^+_{\bar\mu}(x) \coloneqq \lim_{n\to\infty} b_{x\frestr n}$. Clearly, $f^-_{\bar\mu}(x) \leq f^+_{\bar\mu}(x)$.
  % For each $x\in[T]$, define $I^*_x \coloneqq  \bigcap_{n<\omega}I_{x\frestr n}$ and, for $A\subseteq[T]$, define $I^*(A) \coloneqq \bigcup_{x\in A}I^*_x$.
\end{definition}

Notice that $I_{\concat{t}{\la k\ra}}$ is just one point (that is, $a_t=b_t$) iff $\Lb(I_t)\mu_t(\{\concat{t}{\la k\ra}\}) = 0$. 

Let us look at many basic properties of the objects defined in~\autoref{b040}.  

\begin{lemma}\label{b045}
    Let $s,t\in T$.
    \begin{enumerate}[label=\normalfont(\alph*)]
    \item\label{b045-2}  If $t\notin\max T$ then, for any $ k < \alpha_{t}$,  $$a_{\concat{t}{\langle k \rangle}} = a_{t} + \Lb(I_{t})  \sum_{j < k} \mu_{t}(\concat{t}{\langle j \rangle}) \text{ and } 
    %
    %\item\label{b045-1} For any $k < \alpha_{t}$, 
    b_{\concat{t}{\langle k \rangle}} = a_{t} + \Lb(I_{t})  \sum_{j < k + 1} \mu_{t}(\concat{t}{\langle j \rangle}).$$ 

     \item\label{b045a} $\Lb(I_t) = b_t-a_t = \Xi^{\bar{\mu}}(\{t\})$.
     
     \item\label{b045c} If $t\notin \max T$ then $I_t\menos \{b_t\}\subseteq \bigcup_{k<\alpha_t}I_{\concat{t}{\la k\ra}} \subseteq I_t$. Moreover, $\bigcup_{k<\alpha_t}I_{\concat{t}{\la k\ra}} = I_t$  holds whenever $|\suc_T(t)|<\aleph_0$. As a consequence, $\bigcup_{k<\alpha_t}I_{\concat{t}{\la k\ra}}\menos Q_{\bar\mu} = I_t\menos Q_{\bar\mu}$.
    
     \item\label{b045gg} If $s\perp t$ then $I_s\cap I_t \subseteq Q_{\bar\mu}$ and $|I_s\cap I_t|\leq 1$. As a consequence, when $s \perp t$, $(I_s \menos Q_{\bar\mu})\cap(I_t\menos Q_{\bar\mu}) = \emptyset$.  
  \end{enumerate}
\end{lemma}

\begin{PROOF}{\ref{b045}} 
    We have that \ref{b045-2} is clear from the definition of $I_{\concat{t}{\langle k \rangle}}$. 
    
    \ref{b045a}: Proceed by induction on $n=|t|$. When $n=0$, $\Xi^{\bar{\mu}}(\{\la\ \ra\}) = 1 =\Lb(I_{\la\ \ra})$. Assume $n+1<\alt(T)$ and that the claim is true for all $s\in \Lv_n(T)$. Let $t\in \Lv_{n+1}(T)$, that is, $t=\concat{s}{\la k\ra}$ where $s \coloneqq t\frestr n\in \Lv_n(T)$ and $k \coloneqq t(n)<\alpha_s$. By induction hypothesis, $\Lb(I_s) = \Xi^{\bar{\mu}}(\{s\})$, so $\Lb(I_t) = \Lb(I_s)\mu_s(\{t\}) = \Xi^{\bar{\mu}}(\{s\})\mu_s(\{t\})= \Xi^{\bar{\mu}}(\{t\})$, where the last equality holds by~\autoref{pp13}. 
    
    \ref{b045c}: Observe that $\seq{b_{\concat{t}{\la k\ra}}}{k<\alpha_t}$ is a monotone increasing sequence. Even more, by~\ref{b045-2}, if $\alpha_t=\omega$ then $\lim_{k\to\infty}\sum_{j\leq k}\mu_t(\{\concat{t}{\la j\ra}\}) = \sum_{j<\omega}\mu_t(\{\concat{t}{\la j\ra}\})=1$, so $\lim_{k\to\infty}b_{\concat{t}{\la k\ra}} = a_t + \Lb(I_t) = b_t$. Then, $\seq{I_{\concat{t}{\la j\ra}}}{j<\omega}$ covers $I_t\menos\{b_t\}$. 
    
    To prove the ``moreover'' in the statement, assume that $\vert \suc_{T}(t) \vert < \aleph_{0}$. Hence, $\alpha_t<\omega$ and therefore, $\sum_{j\leq \alpha_t-1}\Lb(I_{\concat{t}{\la j\ra}}) = 1$, so $b_{\concat{t}{\la \alpha_t-1\ra}} = a_t + \Lb(I_t) = b_t$. As a consequence, $\seq{I_{\concat{t}{\la j\ra}}}{j<\alpha_t}$ covers $I_t$.
    
    \ref{b045gg}: If $s\perp t$, then there is some $n<\omega$ such that $s\frestr n = t\frestr n$ and $s(n)\neq t(n)$. Without loss of generality, we assume that $s(n)<t(n)$. Then, $b_{s\frestr(n+1)} \leq a_{t\frestr(n+1)}$, so $|I_{s\frestr(n+1)}\cap I_{t\frestr(n+1)}|\leq 1$ and, in case they intersect, the unique point in the intersection must be in $Q_{\bar\mu}$. Hence, the result follows because $I_s\cap I_t\subseteq I_{s\frestr(n+1)}\cap I_{t\frestr(n+1)}$. This easily implies that $(I_s \menos Q_{\bar\mu})\cap(I_t\menos Q_{\bar\mu}) = \emptyset$.
\end{PROOF}

Notice that, for $s, t \in T$,  $I_t\subseteq I_s$ when $s\subseteq t$ in $T$, so $a_s\leq a_t\leq b_t\leq a_s$. Then, for $x\in[T]$, $\seq{a_{x\frestr n}}{n<\omega}$ is a monotone increasing sequence and $\seq{b_{x\frestr n}}{n<\omega}$ is a monotone decreasing sequence (letting $x\frestr n=x$ in the case that $x$ has finite length and $m\geq |x|$) in $[0,1]$, so their limits exist. This allows us to introduce the following definition. 

\begin{definition}
    Define $f^-_{\bar\mu}\colon [T]\to [0,1]$ and $f^+_{\bar\mu}\colon [T]\to [0,1]$ by $\displaystyle f^-_{\bar\mu}(x) \coloneqq \lim_{n\to\infty} a_{x\frestr n}$ and $\displaystyle f^+_{\bar\mu}(x) \coloneqq \lim_{n\to\infty} b_{x\frestr n}$, respectively. Furthermore, for each $x\in[T]$, define $I^*_x \coloneqq  \bigcap_{n<\omega}I_{x\frestr n}$ and, for $A\subseteq[T]$, set $I^*(A) \coloneqq \bigcup_{x\in A}I^*_x$. 
\end{definition}

Clearly, $f^-_{\bar\mu}(x) \leq f^+_{\bar\mu}(x)$. Next, let us look at some basic properties of $I_{x}^{\ast}$ and $I^{\ast}([t])$.

\begin{lemma}\label{b045p}
    Let $t \in T$. Then:
    \begin{enumerate}[label=\normalfont(\alph*)]
        \item\label{b045b} For $x\in[T]$, $I^*_x=[f^-_{\bar\mu}(x),f^+_{\bar\mu}(x)]$ and 
        $$\Lb(I^*_x) = \lim_{n\to\infty} \Xi^{\bar{\mu}}(\{x\frestr n\})= \prod_{n<|x|}\mu_{x\frestr n}(\{x\frestr (n+1)\}).$$
     
        \item\label{b045f} $I_t\menos \{b_t\}\subseteq I^*([t])\subseteq I_t$. As a consequence,  $I^*([t])\menos Q_{\bar\mu} = I_t\menos Q_{\bar\mu}$.
    \end{enumerate}
\end{lemma}

\begin{PROOF}{\ref{b045p}}
    \ref{b045b}: Since $\seq{I_{x\frestr n}}{n<\omega}$ is a decreasing sequence of  closed intervals, 
    $$I^*_x = \bigcap_{n<\omega}I_{x\frestr n} = \bigcap_{n<\omega}[a_{x\frestr n}, b_{x \rest n}] = [f^-_{\bar\mu}(x),f^+_{\bar\mu}(x)]
    $$
    and $\Lb(I^*_x) = \lim_{n\to \infty}\Lb(I_{x\frestr n})$. The rest follows by~\autoref{b045}~\ref{b045a}. 
    
    \ref{b045f}: Assume that $y\in I^*([t])$, so $y\in I^*_x = \bigcap_{n<\omega}I_{x\frestr n}$ for some $x\in [t]$. Since $x\frestr |t| =t$, we get $y\in I_{x\frestr |t|} = I_t$. 
    
    Now assume that $y\in I_t\menos\{b_t\}$. By recursion on $n\geq |t|$, use \autoref{b045}~\ref{b045c} to find $t_n\in \Lv_n(T)$ such that $y\in I_{t_n}\menos Q_{\bar\mu}$ and $t_{n+1}\supseteq t_n\supseteq t$. This recursion can end when reaching a $t_{n_*}\in \max T$, in which case we set $x\coloneqq t_{n_*}$. Otherwise, 
    %Let $t_{|t|} \coloneqq t$; assume we have $t_n$, that is, $y\in I_{t_n}\menos Q_{\bar\mu}$. Then, by~\autoref{b045}~\ref{b045e}, $y\in I_{t_{n+1}}\menos Q_{\bar\mu}$ for some $t_{n+1}\in \suc_T(t_n)$.
    set $x \coloneqq \bigcup_{n\geq |t|}t_n$. In both cases, $x\in[t]$. Thus, $y\in \bigcap_{n<\omega}I_{x\frestr n} = I^*_x$, so $y\in I^*([t])$ because $x\in [t]$.
\end{PROOF}

As a consequence of~\autoref{b045} and \autoref{b045p}:

\begin{corollary}\label{b045-c1}
    $\{ I_{x}^{\ast} \menos Q_{\bar\mu{}} \colon x \in [T] \}$ is a partition of $[0,1] \menos Q_{\bar{\mu}}$. In particular, for any $y \in [0, 1] \menos Q_{\bar{\mu}}$, there exists an unique $x \in [T]$ such that $y \in I_{x}^{\ast}$. 
\end{corollary}

% \begin{PROOF}{\ref{b045-c1}}
%     On the one hand, we have that if $x\neq y$ in $[T]$ then $(I^*_x \menos Q_{\bar\mu})\cap(I^*_y\menos Q_{\bar\mu}) = \emptyset$, as a consequence of~\autoref{b045}~\ref{b045gg}. Now, let us show that $\{ I_{x}^{\ast} \menos Q_{\bar{\mu}} \colon x \in [T] \}$ covers $[0, 1] \menos Q_{\bar{\mu}}$. By~\autoref{b045}~\ref{b045f}, $$[0,1]\menos Q_{\bar\mu} = I_{\la\ \ra}\menos Q_{\bar\mu} = I^*([\la\ \ra])\menos Q_{\bar\mu} = I^*([T])\menos Q_{\bar\mu} = \bigcup_{x\in[T]}I^*_x\menos Q_{\bar\mu},$$ as is required.  
% \end{PROOF}

We can compute the functions $f_{\bar{\mu}}^{-}$ and $f_{\bar{\mu}}^{+}$ for the probability trees from~\autoref{b036}.  

\begin{example}\label{b042}
    \
    \begin{enumerate}[label=\normalfont(\arabic*)]
        \item\label{b042-1} For $t\in T={}^{<\omega}2$, $\la I_{\concat{t}{\la0\ra}},I_{\concat{t}{\la1\ra}}\ra$ results from splitting $I_t$ in half. It can be proved that $a_t = \sum_{i<|t|} \frac{t(i)}{2^{i+1}}$, $b_t = a_t + 2^{-|t|}$ and, for $x\in \cantor$, $f^-_{\bar\mu}(x) = f^+_{\bar\mu}(x) = \sum_{i<\omega}\frac{x(i)}{2^{i+1}}$, hence $I^*_x$ is a singleton.
        
        \item\label{b042-2} For $t\in T = {}^{<\omega}\omega$, $I_{\concat{t}{\la 0\ra}}$ is the first half of $I_t$, $I_{\concat{t}{\la 1\ra}}$ is the first half of $[a_{\concat{t}{\la 0\ra}},b_t]$, $I_{\concat{t}{\la 2\ra}}$ is the first half of $[a_{\concat{t}{\la 1\ra}},b_t]$ and, in general, $I_{\concat{t}{\la \ell+1\ra}}$ is the first half of $[a_{\concat{t}{\la \ell\ra}},b_t]$. Therefore, $\bigcup_{\ell<\omega}I_{\concat{t}{\la \ell\ra}} = I_t\menos\{b_t\}$. 
        
        It can be proved by induction on $n=|t|$ that $\Lb(I_t)=2^{-|t|+\sum_{k<|t|}t(k)}\leq 2^{-|t|}$. As a consequence, $I^*_x$ is a singleton for all $x\in\baire$.
        
        \item\label{b042-3} For $t\in T = {}^{<\omega}\omega$, 
        \[I_{\concat{t}{\la\ell\ra}} = \left\{
        \begin{array}{ll}
           \{a_t\} & \text{ if $\ell<5$,}\\
           I_t & \text{ if $\ell=5$,}\\
           \{b_t\} & \text{ if $\ell>5$.}
        \end{array}
        \right.\]
        As a consequence, $I_t = [0,1]$ when $t\in {}^{<\omega}\omega$ is the constant sequence of $5$. Otherwise, let $n<|t|$ be the minimum such that $t(n)\neq5$. If $t(n)<5$ then $I_t =\{0\}$, and if $t(n)>5$ then $I_t=\{1\}$. Therefore, $I^*_x =[0,1]$ when $x\in\baire$ is the constant sequence of $5$. Otherwise, letting $n<\omega$ be the minimum such that $x(n)\neq 5$, if $x(n)<5$ then $I^*_x = \{0\}$ and, if $x(n)>5$, then $I^*_x = \{1\}$.
    \end{enumerate}
    
    In~\ref{b042-1} and~\ref{b042-2}, $Q_{\bar\mu}$ is dense in $[0,1]$, while in~\ref{b042-3} it is just $\{0,1\}$.
\end{example}

Next, we analyze the conditions under which $a_{t}$ and $b_{t}$ take extreme values, that is, when $a_{t} = 0$ and $b_{t} = 1$. For this, we use the \emph{lexicographic order}: for $s,t\in {}^{<\omega}\omega$, $s\llex t$ iff there is some $n<\omega$ such that $s\frestr n = t\frestr n$ and $s(n)<t(n)$. Notice that $\llex$ is a partial order on ${}^{<\omega}\omega$ but not necessarily linear, e.g.\ comparable nodes in ${}^{<\omega}\omega$ are not $\llex$-comparable. However, $\llex$ is linear at any level of ${}^{<\omega}\omega$.

\begin{lemma}\label{b046}
    Let $t \in T$. Then: 
    \begin{enumerate}[label=\normalfont(\alph*)]
        \item\label{b046.0} If $s\in T$ and $s\llex t$ then $b_s\leq a_t$. 
        \item\label{b046.1} $a_t=0$ iff\/ $\Xi^{\bar\mu}(\{s\})=0$ for all $s\in T$ satisfying $s\llex t$.

        \item\label{b046.2} $b_t=1$ iff\/ $\Xi^{\bar\mu}(\{s\})=0$ for all $s\in T$ satisfying $t\llex s$.
    \end{enumerate}
\end{lemma}

\begin{PROOF}{\ref{b046}}
   \ref{b046.0}: If $s\llex t$ then there is some $n<\omega$ such that $s\frestr n = t\frestr n$ and $s(n)<t(n)$. This implies that $b_{s\frestr (n+1)} \leq a_{t\frestr (n+1)}$. On the other hand, $b_s\leq b_{s\frestr (n+1)}$ and $a_{t\frestr (n+1)} \leq a_t$, so $b_s\leq a_t$. 
   
   Next, we only show~\ref{b046.1}, as~\ref{b046.2} is similar. If $a_t=0$, $s\in T$ and $s\llex t$, then $b_s\leq a_t =0$ by~\ref{b046.0}, so $a_s=b_s=0$, i.e.\ $\Xi^{\bar\mu}(\{s\})=b_s-a_s=0$ by \autoref{b045}~\ref{b045a}.

   Conversely, assume that $\Xi^{\bar\mu}(\{s\})=0$ for all $s\in T$ such that $s\llex t$. We show by induction on $n\leq |t|$ that $a_{t\frestr n} = 0$. This is clear for $n=0$. Now assume that $n<|t|$ and $a_{t\frestr n}=0$. For every $k<t(n)$, $\Lb(I_{\concat{t\frestr n}{\la k\ra}})=\Xi^{\bar\mu}(\{\concat{t\frestr n}{\la k\ra}\})=0$ because $\concat{t\frestr n}{\la k\ra}\llex t$, so $a_{\concat{t\frestr n}{\la k\ra}}=b_{\concat{t\frestr n}{\la k\ra}}=0$. By induction on $k\leq t(n)$, it is easy to show that $a_{\concat{t\frestr n}{\la k\ra}}=a_{t\frestr n}=0$, so $a_{t\frestr(n+1)}=0$. 
    %    
    % We prove~\ref{b046.1} and~\ref{b046.2} simultaneously by induction on $n < \alt(T)$. The base case \red{$n = 1$} is trivial. Now, assume that we have~\ref{b046.1} and~\ref{b046.2} for any $s \in \Lv_{n}(T)$ and consider $t \in \Lv_{n+1}(T)$. Then, $t = \concat{s}{\langle k_{t} \rangle}$, where $s = t \rest n \in \Lv_{n}(T)$ and $k_{t} = t(n) < \alpha_{s}$. 
    %
    % First, assume that $a_{t} = 0$. This implies that $a_{s} = 0$. By induction hypothesis, for any $0 < j < n$, $s(j) = 0$. Now, if $k_{t} = t(n) > 0$, then $a_{\concat{s}{\langle k_{t} - 1 \rangle}} < a_{t} = 0$, which is not possible. Therefore, $k_{t} = 0$ and thus, for any $0 < j < n + 1$, $t(j) = 0$. To show the converse, suppose that for any $0 < j < n + 1$, $t(j) = 0$. In particular, $s(j) = 0$ whenever $0 < j < n$, and $k_{t} = 0$. Therefore, by the induction hypothesis, $a_{s} = 0$. Finally, $a_{t} = a_{\concat{s}{\langle 0 \rangle}} = a_{s} = 0$.  
    %
    % Second, assume that $b_{t} = 1$. This implies that $b_{s} = 1$ and by the induction hypothesis, for any $0 < j < n$, $s(j) = \alpha_{s \rest (j - 1)}$. Also, $k_{t} = \alpha_{t \rest (n - 1)} - 1$ because otherwise, we have $1 = b_{t} < b_{\concat{s}{\langle k_{t}} \rangle} = b_{\concat{s}{\langle \alpha_{t \rest (n - 1)} - 1} \rangle}$, which is not possible. Finally, the converse follows from the induction hypothesis and~\autoref{b045}~\ref{b045c} because $t$ is a finitely splitting node. 
\end{PROOF}

As a consequence, any element in $Q_{\bar{\mu}}\menos\{1\}$ has the form $a_{t}$ for some $t \in T$.

\begin{corollary}\label{b047}
    $Q_{\bar\mu} = \set{a_t}{t\in T}\cup\{1\}$.
\end{corollary}

\begin{PROOF}{\ref{b047}}
    Let $b_{t} \in Q_{\bar{\mu}}$ such that $b_{t} < 1$. By~\autoref{b046}~\ref{b046.2}, we can find some $s\in T$ such that $t\llex s$ and $\Xi^{\bar\mu}(\{s\})>0$, i.e.\ $b_s-a_s>0$ by \autoref{b045}~\ref{b045a}. Then, there is some $m<|t|$ such that $t\frestr m = s\frestr m$ and $t(m)<s(m)$. We show by induction on $m+1\leq n \leq |t|$ that $b_{t\frestr n}$ has the form $a_{t'}$ for some $t'\in T$. In the case $n=m+1$, since $t(m)<s(m)<\alpha_{t\frestr m}$, we get $t(m)+1<\alpha_{t\frestr m}$, so $b_{t\frestr (m+1)} = a_{\concat{t\frestr m}{\la t(m)+1\ra}}$. Now assume that $m+1\leq n<|t|$ and $b_{t\frestr n}=a_{t'}$ for some $t'\in T$. If $t(n)+1<\alpha_{t\frestr n}$ then $b_{t\frestr (n+1)} = a_{\concat{t\frestr n}{\la t(n)+1\ra}}$, otherwise, if $t(n)+1=\alpha_t$ then $b_{t\frestr (n+1)} = b_{t\frestr n} = a_{t'}$. 
    %
    %$n < \vert t \vert$ such that $t(n) + 1 < \alpha_{s}$, where $s \coloneqq t \rest (n - 1)$. Since $t(n) + 1 < \alpha_{s}$, $b_{t} = b_{\concat{s}{\langle t(n) \rangle}} = a_{\concat{s}{\langle t(n) + 1} \rangle}.$ 
    % Let $t \in T$ and assume that $b_{t} < 1$. If $\Xi_{\vert t \vert}^{\bar{\mu}}(\{ t \}) = 0$, \red{the result is trivial} \Diego{no porque aqu\'i no terminar\'ia la prueba, lo ``trivial" es que $b_t=a_t$}, so assume that it is positive. By~\autoref{b046}~\ref{b046.2}, we can find a $k < \vert t \vert$ such that $\red{t(k) \coloneqq j} < \alpha_{s} -1$ (\Diego{al rev\'es, aunque ese $j$ nunca se usa}), where \red{$s \coloneqq t \rest (k-1) \in T$} \Diego{no es obvio?}. Since $j < \alpha_{s} - 1$, by~\autoref{b040} we have that $b_{t} = b_{\concat{s}{\langle k \rangle}} = a_{\concat{s}{\langle k + 1 \rangle}}$. %\in \{ a_{r} \colon r \in T \}$. 
\end{PROOF}

The construction of $\seq{I_t}{t\in T}$ yields the following important connection between $[T]$ and $[0,1]$.

\begin{definition}\label{b050-1}
    Define $g_{\bar\mu}\colon [0,1]\menos Q_{\bar\mu}\to [T]$ such that, for any $y \in [0,1]\menos Q_{\bar\mu}$, $g_{\bar\mu}(y)$ is the unique $x\in[T]$ such that $y\in I^*_x$.
\end{definition}

Notice that $g_{\bar{\mu}}$ is well-defined by virtue of~\autoref{b045-c1}. Now,  we will use it to connect $\calB_{T}$ with $\calB([0, 1])$. 

\begin{lemma}\label{b050}
    The map $g_{\bar\mu}$ is continuous. Moreover:
    \begin{enumerate}[label=\normalfont(\alph*)]    
        \item\label{b050b} For any $A\subseteq[T]$, $g_{\bar\mu}^{-1}[A] = I^*(A)\menos Q_{\bar\mu}$, and 
    
        \item\label{b050c} If $A\in\Bwf_T$ then $I^*(A)\in\Bwf([0,1])$.
    \end{enumerate}
\end{lemma}

\begin{PROOF}{\ref{b050}}    
The continuity of $g_{\bar\mu}$ follows from~\ref{b050b}: by~\autoref{b045}~\ref{b045f}, for any $t\in T$, $g^{-1}_{\bar\mu}\big[[t]\big] = I^*([t])\menos Q_{\bar\mu} = I_t\menos Q_{\bar\mu}$, which is an open set in $[0,1]\menos Q_{\bar\mu}$. 

    \ref{b050b}: Let $A\subseteq [T]$. On the one hand, if $y\in I^*(A)\menos Q_{\bar\mu}$ then $y\in I^*_x\menos Q_{\bar\mu}$ for some $x\in A$. By~\autoref{b050-1}, $x = g_{\bar\mu}(y)$, so $y\in g^{-1}_{\bar\mu}[A]$. On the other hand, if $y\in g^{-1}_{\bar\mu}[A]$ then $y\in [0,1]\menos Q_{\bar\mu}$ and $x \coloneqq   g_{\bar\mu}(y)\in A$, so $y\in I^*_x$ and, thus, $y\in I^*(A)\menos Q_{\bar\mu}$.
    
    \ref{b050c}: By \autoref{b017}, $g_{\bar\mu}$ is a Borel map. Then, for any $A\in\Bwf_T$, $I^*(A)\menos Q_{\bar\mu} \in \Bwf([0,1]\menos Q_{\bar\mu})$. Also notice that, by \autoref{b016}, $\Bwf([0,1]\menos Q_{\bar\mu}) = \Bwf([0,1])|_{[0,1]\menos Q_{\bar\mu}}$, so there is some $B\in \Bwf([0,1])$ such that $I^*(A)\menos Q_{\bar\mu} = B\menos Q_{\bar\mu}$. But $Q_{\bar\mu}\in\Bwf([0,1])$ because it is countable, so $B\menos Q_{\bar\mu}\in \Bwf([0,1])$, that is, $I^*(A)\menos Q_{\bar\mu}\in\Bwf([0,1])$. On the other hand, $I^*(A)\cap Q_{\bar\mu} \in \Bwf([0,1])$ because it is countable, so $I^*(A) =  [I^*(A)\menos Q_{\bar\mu}] \cup [I^*(A)\cap Q_{\bar\mu}]$ is Borel in $[0,1]$. 
    %
    %Let $\Bwf' \coloneqq \set{B\subseteq [T]}{I^*(B)\in \Bwf([0,1])}$. We show that this is a $\sigma$-algebra containing the basic clopen subsets of $[T]$. If $t\in T$, by \autoref{b045}~\ref{b045f}, $I^*([t])\menos Q_{\bar\mu} = I_t\menos Q_{\bar\mu}$, which is Borel in $[0,1]$ (recall that any countable subset of $[0,1]$ is Borel in $[0,1]$). On the other hand, $I^*([t])\cap Q_{\bar\mu}\in \Bwf([0,1])$ because it is countable, so $I^*([t]) =  (I_t\menos Q_{\bar\mu}) \cup (I^*([t])\cap Q_{\bar\mu})$ is in $\Bwf([0,1])$. Therefore, $[t]\in \Bwf'$.
    % 
    % Let $\set{A_n}{n<\omega}\subseteq \Bwf'$ and $A \coloneqq \bigcup_{n<\omega} A_n$. It is not hard to show that $I^*(A) = \bigcup_{n<\omega}I^*(A_n)$, so $A_n\in\Bwf'$ means that $I^*(A_n)\in\Bwf([0,1])$ for all $n<\omega$, hence $I^*(A)\in \Bwf([0,1])$, i.e.\ $A\in \Bwf'$. 
    % 
    % Now assume that $B\in\Bwf'$, i.e.\ $B\subseteq[T]$ and $I^*(B)\in\Bwf([0,1])$. By \autoref{b045}~\ref{b045h} and~\ref{b045i},  % $$([0,1]\menos Q_{\bar\mu})\menos I^*(A) = \bigcup_{x\in [T] \menos A} I^*_x\menos Q_{\bar\mu} = I^*([T]\menos A)\menos Q_{\bar\mu},$$ 
    % and both $[0,1]\menos Q_{\bar\mu}$ and $I^*(A)$ are in $\Bwf([0,1])$, so $I^*([T]\menos A)\menos Q_{\bar\mu}$ is in $\Bwf([0,1])$. On the other hand, $I^*([T]\menos A)\cap Q_{\bar\mu}\in \Bwf([0,1])$, so $I^*([T]\menos A)\in \Bwf([0,1])$, i.e.\ $[T]\menos A\in \Bwf'$.
    % 
    % This shows that $\Bwf'$ is a $\sigma$-algebra over $[T]$ containing $\set{[t]}{t\in T}$. Since the latter generates $\Bwf_T$, we conclude that $\Bwf_T\subseteq\Bwf'$. 
\end{PROOF}

We are ready to introduce a more precise definition of $\lambda^{\bar{\mu}}$.

\begin{definition}\label{b055}
  Define $\lambda^{\bar{\mu}}\colon \Bwf_T\to[0,1]$ by $\lambda^{\bar{\mu}}(A) \coloneqq  \Lb(I^*(A))$, which is well-defined by \autoref{b050}~\ref{b050c}. 
\end{definition}

Although this definition uses the representation of $T$ that we fixed, we will show in \autoref{b078} that $\lambda^{\bar{\mu}}$ does not depend on the representation of $T$.

\begin{theorem}\label{b060}
  The map $\lambda^{\bar{\mu}}$ from \autoref{b055} is a probability measure on $\Bwf_T$ such that $ \lambda^{\bar{\mu}}([t]) = \Xi^{\bar{\mu}}(\{t\})$ for all $t \in T$.
\end{theorem}

\begin{PROOF}{\ref{b060}}
    By \autoref{b028}, $\langle [T], \calA', \lambda' \rangle$ is a measure space, where $\Awf' \coloneqq  g^{\to}_{\bar\mu}(\Bwf([0,1]\menos Q_{\bar\mu}))$ and  $\lambda'(A) \coloneqq  \Lb(g^{-1}_{\bar\mu}[A])$ for any $A\in \Awf'$. Since $g_{\bar\mu}$ is a Borel function, $\Bwf_T \subseteq \Awf'$. Then, for any $A\in\Bwf_T$, since $Q_{\bar\mu}$ is countable,
    $$\lambda'(A) = \Lb(g^{-1}_{\bar\mu}[A]) = \Lb(I^*(A)\menos Q_{\bar\mu}) = \Lb(I^*(A)) = \lambda^{\bar{\mu}}(A).
    $$  
    As a consequence, $\lambda^{\bar{\mu}} =\lambda'\frestr\Bwf_T$, which is a measure. On the other hand, by~\autoref{b045}, for any $t \in T$, 
    $\lambda^{\bar{\mu}}([t]) = \lambda'([t]) = \Lb(I^*([t])\menos Q_{\bar\mu}) = \Lb(I_t\menos Q_{\bar\mu}) = \Lb(I_t) = \Xi^{\bar{\mu}}(\{t\}),
    $
    as required. In particular, 
    $\lambda^{\bar{\mu}}([T]) = \lambda^{\bar{\mu}}([\la\ \ra]) = \Lb(I_{\la\ \ra}) = 1$. Thus, $\lambda^{\bar{\mu}}$ is a probability measure. 
\end{PROOF}

\begin{remark}\label{b065}
   In the previous proof, we do not always have $\Awf'=\Bwf_T$, neither the equivalent formulation ``$I^*(A)\in\Bwf([0,1])$ implies $A\in\Bwf_T$''. For example, when $\Xi^{\bar{\mu}}(\{t\}) = 0$ for some $t\in T$ such that $[t]$ is uncountable, $\lambda^{\bar{\mu}}([t])=0$ and $I^*([t]) = I_t$ is just a singleton, and then $I^*(A) = I_t\in \Bwf([0,1])$ for every non-empty $A\subseteq [t]$, and there are many such $A$ that are not Borel in $[T]$. This is not much of a problem anyway because the converse holds for the completions of $\la[T],\Bwf_T,\lambda^{\bar{\mu}}\ra$ and $\la [0,1], \Bwf([0,1]),\Lb\ra$. %(see~\autoref{b100}).
   For more details, see \autoref{b0100}, \autoref{b0103}, and \autoref{b101}.
   
   Although it is possible that $\lambda^{\bar{\mu}}([t])=0$ for some $t\in T$, this is a very undesirable situation that is typically avoided. When $\lambda^{\bar{\mu}}([t])>0$ for all $t\in T$ and all singletons in $[T]$ have measure zero (which implies $[T]=\lim T$), the equivalence discussed above will hold (see~\autoref{b095}).
\end{remark}

% \begin{theorem}\label{b075}
%   The map $\lambda^{\bar{\mu}}$ is a measure on $\Bwf_T$ such that $\lambda^{\bar{\mu}}([T]) = \Xi_{\vert t \vert}^{\bar{\mu}}(\{ t \})$ whenever $t \in T$. 
% \end{theorem}

% \begin{PROOF}{\ref{b075}}
%     Assume that $\lambda'$ is a measure on $\Bwf_T$ such that $\lambda'([t]) = \Xi_{|t|}^{\bar{\mu}}(\{t\})$ for all $t\in T$. Define 
%     \[\Fwf \coloneqq \bigset{\bigcup_{t\in F}[t]}{F\subseteq \Lv_n(T),\ n<\omega}.\]
%     This is an algebra of sets over $[T]$ and every $C\in\Fwf$ is clopen in $[T]$, so $\Fwf\subseteq \Bwf_T$. It is clear that $[t]\in \Fwf$ for all $t\in T$, so $\sigma(\Fwf) = \Bwf_T$. On the other hand, if $n<\omega$ and $F\subseteq \Lv_n(T)$, then
%     \[\lambda'\left(\bigcup_{t\in F}[t]\right) = \sum_{t\in F}\lambda'([t]) = \sum_{t\in F} \Xi_{\vert t \vert}^{\bar{\mu}}(\{ t \}) =  \sum_{t\in F} \lambda^{\bar{\mu}}([t]) = \lambda^{\bar{\mu}}\left(\bigcup_{t\in F}[t]\right).\]
%     Therefore, $\nu \coloneqq \lambda'\frestr\Fwf = \lambda^{\bar{\mu}}\frestr\Fwf$, which is a measure on $\Fwf$. Both $\lambda'$ and $\lambda^{\bar{\mu}}$ are probability measures on $\sigma(\Fwf)$ extending $\nu$, so $\lambda'= \lambda^{\bar{\mu}}$ by~\autoref{b070}.
% \end{PROOF}

The uniqueness of $\lambda^{\bar\mu}$ in~\autoref{ppp10} simply follows by \autoref{b070}. Therefore: %the construction of $\lambda^{\bar{\mu}}$ does not depend on the representation of $T$. 

\begin{corollary}\label{b078}
  The measure $\lambda^{\bar{\mu}}$ does not depend on the representation of $T$.
\end{corollary}

The construction of $\lambda^{\bar\mu}$ uses the Lebesgue measure on $[0,1]$ and, unlike the first construction in \autoref{ppp10}, does not rely on \autoref{b035-3}, neither on its consequences. Using the second construction, we can prove \autoref{b035-3} more directly. 

\begin{PROOF}[Second proof of \autoref{b035-3}]{\ref{b035-3}}\label{proofb035-3}
    Let $\Xi\in\IP$, $T\coloneqq T_\Xi$ and let $A$ be a front of $T$. Since $\Pi$ is surjective, there is some $\bar\mu\in\TP$ such that $\Xi=\Xi^{\bar\mu}$. Hence, $\Xi^{\lambda^{\bar\mu}} = \Xi^{\bar\mu}$ by \autoref{b060}. 

    Let $t\in T$ below $A$ and $A_{\geq t}\coloneqq \set{s\in A}{t\leq s}$. It is enough to show that $\Xi(\{t\}) = \Xi(A_{\geq t})$. Since $A$ is a front and $t$ is below $A$, we obtain that $[t] = \bigcup_{s\in A_{\geq t}}[s]$ is a disjoint union, so
    \[\Xi(\{t\}) = \lambda^{\bar\mu}([t]) = \sum_{s\in A_{\geq t}}\lambda^{\bar\mu}([s]) = \sum_{s\in A_{\geq t}} \Xi(\{s\}) = \Xi(A_{\geq t}).\qedhere\]
\end{PROOF}

The function $g_{\bar\mu}$ has more properties than it appears to have. Under certain conditions, it is a topological embedding into $[T]$. To construct the inverse function, we look at the measure zero points of $[T]$.

\begin{lemma}\label{b077}
    Let $S$ be a subtree of $T$ with $\max S\subseteq \max T$, and let $x \in [T]$. Then:      
    \begin{enumerate}[label=\normalfont(\alph*)]
     \item\label{b077a} $\displaystyle \lambda^{\bar{\mu}}([S]) = \inf_{n<\alt(T)} \Xi^{\bar{\mu}}(\Fr_n(S))$.
     
     \item\label{b077b} $\lambda^{\bar{\mu}}(\{x\}) = \Lb(I^*_x)$.
     
     \item\label{b077c} $\lambda^{\bar{\mu}}(\{x\}) =0$ iff $I^*_x$ is a singleton.
  \end{enumerate}
\end{lemma}

\begin{PROOF}{\ref{b077}}
    \ref{b077a}: For any $n < \alt(T)$, consider $C_{n} \coloneqq \bigcup \{ [t]_{T} \colon t \in \Fr_{n}(S) \}$, which is a pairwise disjoint union. Then $\langle C_{n} \colon n < \alt(T) \rangle$ is a decreasing sequence of clopen sets in $[T]$ whose intersection is $[S]$. Furthermore, by~\autoref{b060}, it is clear that, for any $n < \alt(T)$, $\lambda^{\bar{\mu}}(C_{n}) = \Xi^{\bar{\mu}}(\Fr_{n}(S))$. As a consequence, 
    $$ \lambda^{\bar{\mu}} ( \left [S] \right) =  \lambda^{\bar{\mu}} \left( \bigcap_{n < \alt(T)} C_{n} \right) =  \inf_{n < \alt(T)} \lambda^{\bar{\mu}}(C_{n}) = \inf_{n < \alt(T)} \Xi^{\bar{\mu}}(\Fr_{n}(S)).
    $$ 

    \ref{b077b}: By the definition of $\lambda^{\bar{\mu}}$, 
    $ \lambda^{\bar{\mu}}(\{ x \}) = \Lb(I^{\ast}(\{ x \})) = \Lb(I_{x}^{\ast}). $

    \ref{b077c}: On the one hand, if $I_{x}^{\ast}$ is a singleton, then by~\ref{b077b} $\lambda^{\bar{\mu}}(\{ x \}) = \Lb(I_{x}^{\ast}) = 0$. On the other hand, by~\autoref{b045p}~\ref{b045b}, $I_{x}^{\ast}$ is a closed interval, namely $[f_{\bar{\mu}}^{-}(x), f_{\bar{\mu}}^{+}(x)]$. Therefore, if $\lambda^{\bar{\mu}}(\{ x \}) = 0$, then $\Lb(I_{x}^{\ast}) = 0$, so $I_{x}^{\ast}$ must be a singleton. 
\end{PROOF}

Using~\autoref{b077}~\ref{b077c}, we can introduce a sort of inverse of $g_{\bar{\mu}}$ as follows (inverse in the sense of~\autoref{b087}~\ref{b087g}).

\begin{definition}\label{m45}
  \ 
    \begin{enumerate}[label=\normalfont (\arabic*)]
       \item For $\lambda\in\BP$, define $V^*_\lambda\coloneqq \set{x\in[T_\lambda]}{\lambda(\{x\})=0}$, the \emph{free part of $[T_\lambda]$.}
    
       \item Let $f_{\bar\mu}$ be the function with domain $V^*_{\bar\mu}\coloneqq V^*_{\lambda^{\bar\mu}}$ such that, for $x\in V^*_{\bar\mu}$, $f_{\bar\mu}(x)$ is the unique point in $I^*_x$. Notice that $V^*_{\bar\mu}$ could be empty.

       \item For $A \subseteq [T]$, define $G^{\ast}_{A} \coloneqq \bigcup_{x \in A \menos  V^*_{\bar\mu}} \left( f_{\bar{\mu}}^{-}(x), f_{\bar{\mu}}^{+}(x) \right)$ and $G^{\ast}_{\bar{\mu}} \coloneqq G^*_{[T]}$, which are open in $[0,1]$.

       \item Denote $N^*_{\bar\mu}\coloneqq N^*_{\lambda^{\bar\mu}}$ (see \autoref{pp134}), which is the largest open measure zero subset of $[T]$ (see \autoref{pp140}). Recall that $[T^+_{\bar\mu}] = [T]\menos N^*_{\bar\mu}$, so it has measure $1$.
	\end{enumerate}
\end{definition}

%Notice that, $N_{\bar{\mu}}^{\ast}$  is an open set of $\lambda^{\bar{\mu}}$- measure zero sets (because $T$ is countable), and for any $A\subseteq [T]$, $G_{A}^{\ast}$ is open in $[0, 1]$. Furthermore, by~\autoref{b077}~\ref{b077c}, $\lambda^{\bar{\mu}}([T^+_{\bar\mu}])=1$.  

\begin{lemma}\label{b087}
    \ 
    \begin{enumerate}[label=\normalfont(\alph*)]
    \item\label{b087a} $[T]\menos V^*_{\bar\mu}$ is countable. In particular, $V^*_{\bar\mu}\in \Bwf_T$.
    
    \item\label{b087b} For $y\in[0,1]\menos Q_{\bar\mu}$, if $g_{\bar\mu}(y)\in V^*_{\bar\mu}$ then $f_{\bar\mu}(g_{\bar\mu}(y))=y$.
    
    \item\label{b087c} If $x\in V^*_{\bar\mu}$ and $f_{\bar\mu}(x)\notin Q_{\bar\mu}$ then $g_{\bar\mu}(f_{\bar\mu}(x))=x$.
    
    \item\label{b087d} $f_{\bar\mu}$ is continuous and $(Q_{\bar\mu}\cup \ran f_{\bar\mu})\cap G^*_{\bar\mu} = \emptyset$.
    
    \item\label{b087e} $\ran g_{\bar\mu} = [T]\menos f^{-1}_{\bar\mu}[Q_{\bar\mu}]$. In particular, $\ran g_{\bar\mu}\in \Bwf_T$.
    
    \item\label{b087f} $f_{\bar\mu}[\ran g_{\bar\mu}] = \ran f_{\bar\mu} \menos Q_{\bar\mu}$.
    
    \item\label{b087g} $f_{\bar\mu}\frestr \ran g_{\bar\mu}$ is an homeomorphism from $V^*_{\bar\mu}\menos f^{-1}[Q_{\bar\mu}]$ onto $\ran f_{\bar\mu}\menos Q_{\bar\mu}$ with inverse $g_{\bar\mu}\frestr (\ran f_{\bar\mu}\menos Q_{\bar\mu})$.
    
    \item\label{b087h} $N^*_{\bar\mu}\subseteq f^{-1}_{\bar\mu}[Q_{\bar\mu}]$ and $\max T \subseteq N^*_{\bar\mu}\cup([T]\menos V^*_{\bar\mu})$, i.e.\ $V^*_{\bar\mu}\menos N^*_{\bar\mu}\subseteq\lim T$.
    
    \item\label{b087i} $f^{-1}_{\bar\mu}[Q_{\bar\mu}]\menos N^*_{\bar\mu}$ is countable and $\lambda^{\bar{\mu}}(f^{-1}_{\bar\mu}[Q_{\bar\mu}]) = 0$.

    \item\label{b087j} $f_{\bar\mu}[A] = I^*(A \cap  V^*_{\bar\mu})$ for all $A\subseteq [T]$. 
    \end{enumerate}
\end{lemma}

\begin{figure}[h]
    \centering
    \begin{tikzpicture}[scale=0.86]
    % Cuadrado izquierdo
    \draw[thick] (0, 0) rectangle (5, 5);
    % Interseccion homeo cuadro derecho
    \fill[redun!15] (10, 1) rectangle (15, 3.5);
    
    % Cuadrado en la parte derecha del rectángulo del cuadrado izquierdo

    % Colocar el texto dentro del cuadrado pequeño
    \node at (4.1, 4.25) {\footnotesize \textbf{$N_{\bar{\mu}}^{\ast}$}};
    
    % Cuadrado derecho
    \draw[thick] (10, 0) rectangle (15, 5);
    
    % interseccion homeomorphismo cuadro derecho
    \fill[greenun!15] (0, 1) rectangle (5, 3.5);
    
    % Subcuadrado con solo el perímetro coloreado (azul) en el cuadrado izquierdo
    \draw[thick, greenun] (0, 1) rectangle (5, 5);
    
    % Colocar el nombre [T] arriba, en el centro del cuadrado izquierdo y más arriba
    \node at (2.5, 5.6) {\textbf{$[T]$}};
    % Colocar el nombre [0, 1] arriba, en el centro del cuadrado derecho y más arriba
    \node at (12.5, 5.6) {\textbf{$[0, 1]$}};
    
    % Flecha curva que representa la función f_{\bar{\mu}}, en la parte superior
    \draw[->, thick] (5.5, 5.5) to[out=30, in=150] node[above] {\footnotesize $f_{\bar{\mu}}$} (9.5, 5.5);
    % Flecha curva que representa la función g_{\bar{\mu}}, en la parte inferior (simétrica)
    \draw[<-, thick] (5.5, -0.5) to[out=-30, in=-150] node[below] {\footnotesize $g_{\bar{\mu}}$} (9.5, -0.5);
    
    % define ran g
    \draw[thick] (-1, 0) -- (-1, 1.60);

    \draw[thick] (-1, 1.97) -- (-1, 3.5);

    \draw[thick] (-1, 0) -- (-.88, 0);

    \draw[thick] (-1, 3.5) -- (-0.88, 3.5);
    
    % Dividir la línea a la izquierda del cuadrado izquierdo
    \node at (-1, 1.75) {\footnotesize \textbf{$\mathrm{ran} g_{\bar{\mu}}$}}; % Texto en la línea izquierda

    % Definiendo dom g
    \draw[thick] (16, 0) -- (16, 1.60);
    \draw[thick] (16, 1.97) -- (16, 3.5);
    \node at (16, 1.75) {\footnotesize \textbf{$\mathrm{dom} g_{\bar{\mu}}$}};

    \draw[thick] (16, 0) -- (15.88, 0);
    \draw[thick] (16, 3.5) -- (15.88, 3.5);
    
    % Etiqueta de de f_{\bar{\mu}}[Q] 
    
    \node at (0.9, 4.25) {\textbf{ \footnotesize $f_{\bar{\mu}}^{-1}[Q_{\bar{\mu}}]$}};

    % Etiqueta de Q 
    
    \node at (10.5, 4.25) {\textbf{\footnotesize $Q_{\bar{\mu}}$}};

    % Etiqueta de G

    \node at (14.43, 0.5) {\textbf{\footnotesize $G_{\bar{\mu}}^{\ast}$}};

    % Delimitar la inversa de Q 
    
    \draw[thick] (0, 3.5) -- (5, 3.5);

    %\draw (3.3, 5) -- (3.3, 3.5);

     \draw (3.3, 5) to[out=270, in=180] (4, 3.5);

    % Delimitar Q

    \draw[thick] (10, 3.5) -- (15, 3.5);

    % Etiqueta de dom f

    \draw[thick, greenun] (6, 3.32) -- (6, 5);
    \draw[thick, greenun] (6, 1) -- (6, 2.7);
    \node[greenun] at (6, 3) {\textbf{\footnotesize $V^*_{\bar\mu}$}};

    \draw[thick, greenun] (6, 5) -- (5.88, 5);
    
    \draw[thick, greenun] (6, 1) -- (5.88, 1);
    % Etiqueta de ran f

    \draw[thick, redun] (9, 3.2) -- (9, 4.25);
    
    \draw[thick, redun] (9, 1) -- (9, 2.8);
    
    \node[redun] at (9, 3) {\footnotesize \textbf{$\ran f_{\bar{\mu}}$}};

    \draw[thick, redun] (9, 4.25) -- (9.12, 4.25);
    \draw[thick, redun] (9, 1) -- (9.12, 1);

    % definiendo ran f 

    \draw[thick, redun] (10, 1) rectangle (15, 5);    
\end{tikzpicture}

    \caption{Graphic situation of~\autoref{b087}: $N_{\bar{\mu}}^{\ast} \subseteq f_{\bar{\mu}}^{-1}[Q_{\bar{\mu}}]$,  $\ran f_{\bar{\mu}}$ may include some elements of $Q_{\bar{\mu}}$, and the shaded regions are homeomorphic via $f_{\bar{\mu}} \rest \ran g_{\bar{\mu}}$, whose inverse is $g_{\bar{\mu}} \rest \ran f_{\bar{\mu}}$.}
    \label{fig:lines-divided}
\end{figure}
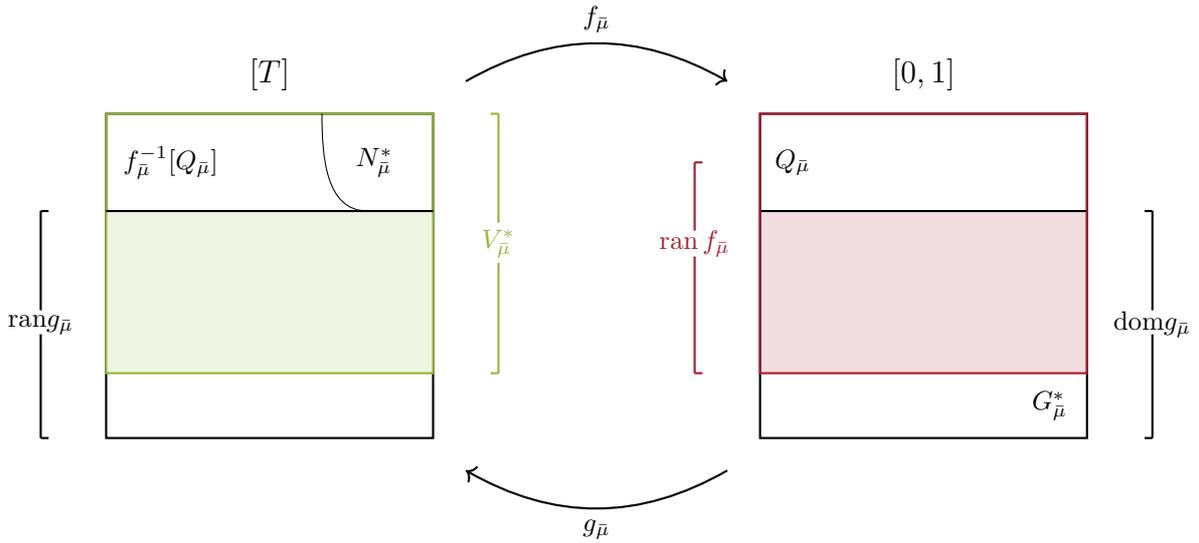

\begin{PROOF}{\ref{b087}}
    \ref{b087a}: This is a direct consequence of~\autoref{b029}.
    
    \ref{b087b}: By the definition of $g_{\bar\mu}$, $y\in I^*_{g_{\bar\mu}(y)}$ for $y\in [0,1]\menos Q_{\bar\mu}$, but $I^*_{g_{\bar\mu}(y)} =\{f_{\bar\mu}(g_{\bar\mu}(y))\}$ when $g_{\bar\mu}(y)\in V^*_{\bar\mu}$, so $f_{\bar\mu}(g_{\bar\mu}(y)) = y$.
    
    \ref{b087c}: If $x\in V^*_{\bar\mu}$ and $f_{\bar\mu}(x)\notin Q_{\bar\mu}$ then $f_{\bar\mu}(x)\in\dom g_{\bar\mu}$. Also $f_{\bar\mu}(x)\in I^*_x$, so $g_{\bar\mu}(f_{\bar\mu}(x))=x$.
    
    \ref{b087d}: Let $x\in V^*_{\bar\mu}$ and $\varp>0$. Since $\{f_{\bar\mu}(x)\}= I^*_x = \bigcap_{n<\omega}I_{x\frestr n} \subseteq(f_{\bar\mu}(x)-\varp,f_{\bar\mu}(x)+\varp)$, there is some $n<\omega$ such that $I_t\subseteq (f_{\bar\mu}(x)-\varp,f_{\bar\mu}(x)+\varp)$, where $t \coloneqq  x\frestr n$. Then, for any $x'\in [t]\cap  V^*_{\bar\mu}$, $f_{\bar\mu}(x')\in I^*([t]) \subseteq I_t$, so $|f_{\bar\mu}(x)-f_{\bar\mu}(x')|<2\varp$. This shows that $f_{\bar\mu}$ is continuous. 
    
    Notice that $\seq{(f_{\bar\mu}(x),f^+_{\bar\mu}(x))}{x\in [T]\menos V^*_{\bar\mu}}$ is a sequence of pairwise disjoint intervals that cannot contain points in $Q_{\bar\mu}$, and thus neither in $\ran f_{\bar\mu}$. Hence, $(Q_{\bar\mu}\cup \ran f_{\bar\mu})\cap G^*_{\bar\mu} = \emptyset$.
    
    \ref{b087e}: If $x\in\ran g_{\bar\mu}$ then $x = g_{\bar\mu}(y)$ for some $y\in[0,1]\menos Q_{\bar\mu}$. Either $x\notin  V^*_{\bar\mu}$ or $x\in V^*_{\bar\mu}$, and in the latter case, 
    $f_{\bar\mu}(x) = y \in [0,1]\menos Q_{\bar\mu}$ by~\ref{b087b}. Thus, $x\in [T]\menos f^{-1}_{\bar\mu}[Q_{\bar\mu}]$. On the other hand, assume that $x\in [T]\menos f^{-1}_{\bar\mu}[Q_{\bar\mu}]$. If $x\notin V^*_{\bar\mu}$ then $\Lb(I^*_x)>0$, so $I^*_x\nsubseteq Q_{\bar\mu}$, i.e.\ there is some $y\in I^*_x\menos Q_{\bar\mu}$. Thus, $x= g_{\bar\mu}(y) \in \ran g_{\bar\mu}$.
    
    In the case that $x\in  V^*_{\bar\mu}$, we have
    $y \coloneqq f_{\bar\mu}(x)\notin Q_{\bar\mu}$, so $x = g_{\bar\mu}(y) \in \ran g_{\bar\mu}$ by~\ref{b087c}.
    
    \ref{b087f}: If $y\in f_{\bar\mu}[\ran g_{\bar\mu}]$ then $y=f_{\bar\mu}(x)$ for some $x\in  V^*_{\bar\mu} \cap \ran g_{\bar\mu}$, so $x\notin f^{-1}_{\bar\mu}[Q_{\bar\mu}]$, i.e.\ $y= f_{\bar\mu}(x) \in \ran f_{\bar\mu}\menos Q_{\bar\mu}$.
    
    For the converse, if $y\in \ran f_{\bar\mu}\menos Q_{\bar\mu}$ then $y=f_{\bar\mu}(x)$ for some $x\in V^*_{\bar\mu}$ and $f_{\bar\mu}(x)\notin Q_{\bar\mu}$, i.e.\ $x\notin f^{-1}_{\bar\mu}[Q_{\bar\mu}]$, so $x\in \ran g_{\bar\mu}$ by~\ref{b087e}. Hence, $y\in f_{\bar\mu}[\ran g_{\bar\mu}]$.
    
    \ref{b087g}: Recall that $ V^*_{\bar\mu} \menos f^{-1}_{\bar\mu}[Q_{\bar\mu}] =  V^*_{\bar\mu}\cap \ran g_{\bar\mu}$ by~\ref{b087e}. 
    If $x\in  V^*_{\bar\mu}\menos f^{-1}_{\bar\mu}[Q_{\bar\mu}]$ then $f_{\bar\mu}(x)\notin Q_{\bar\mu}$, so $g_{\bar\mu}(f_{\bar\mu}(x))$ is defined and it is equal to $x$ by~\ref{b087c}. Conversely, if $y\in\ran f_{\bar\mu}\menos Q_{\bar\mu}$ then $y= f_{\bar\mu}(x)$ for some $x\in V^*_{\bar\mu}\cap \ran g_{\bar\mu}$ by~\ref{b087f}, so $x=g_{\bar\mu}(y')$ for some $y'\in[0,1]\menos Q_{\bar\mu}$. By~\ref{b087b}, $y'=f_{\bar\mu}(g_{\bar\mu}(y')) = f_{\bar\mu}(x) =y$, so $f_{\bar\mu}(g_{\bar\mu}(y))=y$. This shows that $f_{\bar\mu}\frestr \ran g_{\bar\mu}$ is a bijection from $ V^*_{\bar\mu}\menos f^{-1}_{\bar\mu}[Q_{\bar\mu}]$ onto $\ran f_{\bar\mu}\menos Q_{\bar\mu}$ with inverse $g_{\bar\mu}\frestr(\ran f_{\bar\mu}\menos Q_{\bar\mu})$. Since both $f_{\bar\mu}$ and $g_{\bar\mu}$ are continuous, we are done.
    
    \ref{b087h}: If $x\in N^*_{\bar\mu}$ then $x\in[t]$ for some $t\in T$ such that $\lambda^{\bar{\mu}}([t])=0$. Then $\lambda^{\bar{\mu}}(\{x\}) = 0$, so $x\in V^*_{\bar\mu}$. Moreover, $\{f_{\bar\mu}(x)\}=I^*_x = I_t$, so $f_{\bar\mu}(x)\in Q_{\bar\mu}$. 

    On the other hand, if $x\in\max T$ and $x\notin N^*_{\bar\mu}$, then $[x]=\{x\}$ and $\lambda^{\bar\mu}(\{x\})=\lambda^{\bar\mu}([x])>0$, so $x\notin  V^*_{\bar\mu}$. 

    \ref{b087i}: It is enough to show that, for $y\in Q_{\bar\mu}$, $|f^{-1}_{\bar\mu}[\{y\}]\menos N^*_{\bar{\mu}}|\leq 2$. Let $x\in f^{-1}_{\bar\mu}[\{y\}]\menos N^*_{\bar{\mu}}$, so $x\in \lim T_{\bar{\mu}}^+$ (by~\ref{b087h}) and $f_{\bar\mu}(x) = y\in Q_{\bar\mu}$. 
    By the definition of $f_{\bar\mu}$, $I^*_x = \{y\}$, so $a_{x\frestr n}\nearrow y$ and $b_{x\frestr n}\searrow y$. Since $y\in Q_{\bar\mu}$, we must have that, for some $m<\omega$, either $a_{x\frestr n} = y$ for all $n\geq m$, or $b_{x\frestr n} = y$ for all $n\geq m$. We show that at most only one $x \in f_{\bar{\mu}}^{-1}[\{ y \}] \menos N_{\bar{\mu}}^{\ast}$ satisfies that $a_{x\frestr n} = y$ for all but finitely many $n$. Likewise, there is at most one $x$ satisfying $b_{x\frestr n} = y$ for all but finitely many $n$, so $|f^{-1}_{\bar\mu}[\{y\}] \cap [T_{\bar{\mu}}^+]|\leq 2$. %, because we only consider two cases.
    
    Assume that $x,x'\in f^{-1}_{\bar\mu}[\{y\}]\menos N^*_{\bar\mu}$ satisfy that, for some $m'<\omega$, $a_{x\frestr n} = a_{x'\frestr n} = y$ for all $n\geq m'$. If $x\neq x'$ then let $n_0<\omega$ be such that $x\frestr n_0 = x\frestr n_0$ and $x(n_0)\neq x'(n_0)$. Without loss of generality, $x(n_0)<x'(n_0)$. Since $x,x'\in [T_{\bar{\mu}}^+]$, $x\frestr(n_0+1), x'\frestr(n_0+1) \in T_{\bar{\mu}}^+$, so $a_{x\frestr(n_0+1)}<b_{x\frestr(n_0+1)} \leq a_{x'\frestr(n_0+1)}$. Then, for some $n>\max\{m',n_0\}$, $y=a_{x\frestr n}< b_{x\frestr n}\leq a_{x'\frestr n}$, so $y<a_{x'\frestr n}$, a contradiction.
    
    %A similar argument shows that, in the second case, there is no other $x'\in \red{[T_{\bar{\mu}}']}$ satisfying that, for some $m'<\omega$, $b_{x'\frestr n} = y$ for all $n\geq m'$.
    
    Finally, since $f_{\bar\mu}^{-1}[Q_{\bar\mu}]\menos N^*_{\bar\mu}$ is countable and contained in $ V^*_{\bar\mu}$, $$\lambda^{\bar{\mu}}(f^{-1}_{\bar\mu}[Q_{\bar\mu}]) = \lambda^{\bar\mu}(f_{\bar\mu}^{-1}[Q_{\bar\mu}]\menos N^*_{\bar\mu}) + \lambda^{\bar\mu}(N^*_{\bar\mu})= 0.$$

    \ref{b087j}: By the definition of $f_{\bar{\mu}}$, we have   
    $ f_{\bar{\mu}}[A] = \bigcup_{x \in A \cap  V^*_{\bar\mu}} I_{x}^{\ast} = %f_{\bar{\mu}}[A \cap  V^*_{\bar\mu}] = 
    I^{\ast}(A \cap  V^*_{\bar\mu}).$
\end{PROOF}

Now we analyze the effect of applying the functions $f_{\bar{\mu}}$ and $g_{\bar{\mu}}$ to Borel sets and their respective measures.

\begin{theorem}\label{b0100}
    Let $A \subseteq [T]$ and $B \subseteq [0, 1]$. Then: 
    \begin{enumerate}[label=\normalfont(\alph*)]
        \item\label{b0100a} If $A \in \calB_{T}$, then $f_{\bar{\mu}}[A] \in \calB([0, 1])$. Furthermore, $\Lb(f_{\bar{\mu}}[A]) = \lambda^{\bar{\mu}}(A \cap  V^*_{\bar\mu})$ and 
        $$ \lambda^{\bar{\mu}}(A) = \Lb(f_{\bar{\mu}}[A])  + \Lb(G^*_A)=\Lb(f_{\bar{\mu}}[A])  + \sum_{x \in A \menos \dom(f_{\bar{\mu}})} \Lb(I_{x}^{\ast}).$$

        \item\label{b0100b} If $B \in \calB([0, 1])$ then  $f_{\bar{\mu}}^{-1}[B] \in \calB_{T}$ and $\lambda^{\bar{\mu}}(f_{\bar{\mu}}^{-1}[B])= \Lb \left( B \menos G^{\ast}_{\bar{\mu}} \right)$. As a consequence, $\Lb(B) = \lambda^{\bar{\mu}}(f_{\bar{\mu}}^{-1}[B]) + \Lb \left( B \cap G^{\ast}_{\bar{\mu}} \right). $

        \item\label{b0100c}  $A \in \calB_{T}$ iff $A \cap N_{\bar{\mu}}^{\ast} \in \calB_{T}$ and $f_{\bar{\mu}}[A] \in \calB([0, 1])$.  

        \item\label{b0100e} $B \in \calB([0, 1])$ iff $B \cap G^{\ast}_{\bar{\mu}} \in \calB([0, 1])$ and $f_{\bar{\mu}}^{-1}[B] \in \calB_{T}$.

        \item\label{b0100f} $B \in \calB([0, 1])$ iff $B \cap G^{\ast}_{\bar{\mu}} \in \calB([0, 1])$ and $g_{\bar{\mu}}[B] \in \calB_{T}$. In this case, $\Lb(B \menos G^{\ast}_{\bar{\mu}}) = \lambda^{\bar{\mu}}(g_{\bar\mu}[B] \cap  V^*_{\bar\mu})$.

        \item\label{b0100d}   $A \in \calB_{T}$ iff $A \cap N_{\bar{\mu}}^{\ast} \in \calB_{T}$ and $g_{\bar{\mu}}^{-1}[A] \in \calB([0, 1])$. In this case, $\Lb(g_{\bar{\mu}}^{-1}[A]) = \lambda^{\bar{\mu}}(A)$. 
    \end{enumerate}
\end{theorem}

\begin{PROOF}{\ref{b0100}}
    \ref{b0100a}: Let $A \in \calB_{T}$. For $x \in A$: if $x \in A \cap  V^*_{\bar\mu}$, then $I_{x}^{\ast} = \{ f_{\bar{\mu}}(x) \}$; and if $x \in A \menos  V^*_{\bar\mu}$, then $I_{x}^{\ast} = [f_{\bar{\mu}}^{-}, f_{\bar{\mu}}^{+}]$ is an uncountable interval.  Now, 
    $ f_{\bar{\mu}}[A] \in \calB([0, 1])
    $ by~\autoref{b087}~\ref{b087j}, ~\autoref{b050}~\ref{b050c} and because $ V^*_{\bar\mu} \in \calB_{T}$. As a consequence, $\Lb(f_{\bar{\mu}}[A]) = \Lb(I^{\ast}(A \cap  V^*_{\bar\mu})) = \lambda^{\bar{\mu}}(A \cap  V^*_{\bar\mu}).$ Finally, since we have $f_{\bar{\mu}}[A] \cup G_{A}^{\ast}\subseteq I^{\ast}(A) \subseteq f_{\bar{\mu}}[A] \cup G_{A}^{\ast}\cup Q_{\bar\mu}$, it follows that: 
    $$ \lambda^{\bar{\mu}}(A) = \Lb(I^{\ast}(A)) = \Lb \left( f_{\bar{\mu}}[A] \cup G_{A}^{\ast} \right) = \Lb(f_{\bar{\mu}}[A]) + \sum_{x \in A \menos  V^*_{\bar\mu}} \Lb(I_{x}^{\ast}).
    $$

    \ref{b0100b}: If $B \in \calB([0, 1])$ we have that $f_{\bar{\mu}}^{-1}[B] \in \calB_{T}$ because $f_{\bar{\mu}}$ is continuous. On the other hand, by~\ref{b0100a}, 
    $$ \lambda^{\bar{\mu}}(f_{\bar{\mu}}^{-1}[B]) = \lambda^{\bar{\mu}}(f_{\bar{\mu}}^{-1}[B] \cap  V^*_{\bar\mu}) = \Lb(B \cap \ran f_{\bar{\mu}}) = \Lb(B \menos G^{\ast}_{\bar{\mu}}).  
    $$

    \ref{b0100c}: If $A \in \calB_{T}$, then $N_{\bar{\mu}}^{\ast} \cap A \in \calB_{T}$ because $N_{\bar{\mu}}^{\ast}$ is open in $[T]$. Also, $f_{\bar{\mu}}[A] \in \calB([0, 1])$ follows by~\ref{b0100a}. To prove the converse, assume that $A \cap N_{\bar{\mu}}^{\ast} \in \calB_{T}$ and $f_{\bar{\mu}}[A] \in \calB_{T}$. Notice that we can write $A$ as a union of four sets:
    $$
    A = [A \cap N_{\bar{\mu}}^{\ast}] \cup [A \cap f_{\bar{\mu}}^{-1}[Q_{\bar\mu}] \menos N_{\bar{\mu}}^{\ast}] \cup [A \cap  V^*_{\bar\mu} \menos f_{\bar{\mu}}^{-1}[Q_{\bar{\mu}}]]\cup [A \cap [T] \menos  V^*_{\bar\mu}],
    $$
    so it is enough to show that these four sets are Borel in $\calB_{T}$. Since $f_{\bar{\mu}}^{-1}[Q_{\bar\mu}] \menos N_{\bar{\mu}}^{\ast}$ and $[T] \menos  V^*_{\bar\mu}$ are countable (see~\autoref{b087}),  and $A \cap N_{\bar{\mu}}^{\ast} \in \calB_{T}$ by hypothesis, it only remains to prove that $A \cap  V^*_{\bar\mu} \menos f_{\bar{\mu}}^{-1}[Q_{\bar{\mu}}] \in \calB_{T}$. 
    Since $f_{\bar{\mu}}$ is Borel, $f_{\bar{\mu}}[A], \dom g_{\bar{\mu}} \in \calB_{T}$, it follows that $
    A \cap  V^*_{\bar\mu} \menos f_{\bar{\mu}}^{-1}[Q_{\bar{\mu}}] = f_{\bar{\mu}}^{-1}[ f_{\bar{\mu}}[A] \cap \dom g_{\bar{\mu}}] \in \calB_{T}.
    $

    \ref{b0100e}: The implication from left to right follows by~\ref{b0100b} and because $G^{\ast}_{\bar{\mu}}$ is clearly a Borel set. To prove the converse, notice that we can write  $ B = [B \cap \ran f_{\bar{\mu}}] \cup [B \cap G_{\bar{\mu}}^{\ast}] \cup [B \cap Q_{\bar{\mu}}]$. Since $Q_{\bar{\mu}}$ is countable and by the hypothesis, it is enough to show that $B \cap \ran f_{\bar{\mu}}$ is Borel. This holds because, by hypothesis and by~\ref{b0100a}, $B \cap \ran f_{\bar{\mu}} = f_{\bar{\mu}}[f_{\bar{\mu}}^{-1}[B]] \in \calB([0, 1])$.

    \ref{b0100f}: Assume that $B \in \calB([0, 1])$. Then, $B \cap G^{\ast}_{\bar{\mu}} \in \calB([0, 1])$ because $G_{\bar{\mu}}^{\ast}$ is Borel. On the other hand, $g_{\bar{\mu}}[B] \in \calB_{T}$ because $
    g_{\bar{\mu}}[B] = [g_{\bar{\mu}}[B] \cap  V^*_{\bar\mu}] \cup [g_{\bar{\mu}}[B] \cap [T] \menos  V^*_{\bar\mu}]
    $, $[T] \menos  V^*_{\bar\mu}$ is countable, $f_{\bar{\mu}}^{-1}[B \cap \ran f_{\bar{\mu}}] = g_{\bar{\mu}}[B] \cap  V^*_{\bar\mu}$ (by~\autoref{b087}~\ref{b087g}), $f_{\bar{\mu}}$ is a Borel function, and $\ran f_{\bar{\mu}}$ is Borel in $[0, 1]$.   To prove the converse, notice that we can write $ B = [B \cap \ran f_{\bar{\mu}} \menos Q_{\bar{\mu}}] \cup [B \cap G_{\bar{\mu}}^{\ast}] \cup [B \cap Q_{\bar{\mu}}]$. Therefore, it is enough to show that $B \cap \ran f_{\bar{\mu}} \menos Q_{\bar{\mu}}$ is Borel in $\calB([0, 1])$. This holds because $B \cap \ran f_{\bar{\mu}} \menos Q_{\bar{\mu}} = g_{\bar{\mu}}^{-1}[ g_{\bar{\mu}}[B] \cap  V^*_{\bar\mu}]$, $g_{\bar{\mu}}$ is Borel, and $g_{\bar{\mu}}[B]$ and $ V^*_{\bar\mu}$ are Borel sets. Finally, $\Lb(B \menos G^{\ast}_{\bar{\mu}}) = \lambda^{\bar{\mu}}(g[B] \cap  V^*_{\bar\mu})$ follows by~\ref{b0100b} and~\autoref{b087}~\ref{b087g}.

     \ref{b0100d}: Assume that $A \in \calB_{T}$. Then $A \cap N_{\bar{\mu}}^{\ast} \in \calB_{T}$ because $N_{\bar{\mu}}^{\ast}$ is open, and $g_{\bar{\mu}}^{-1}[A] \in \calB([0, 1])$ follows by~\autoref{b050}. To show the converse, write $A$ as follows: 
    $$
    A = [A \cap N_{\bar{\mu}}^{\ast}] \cup [ A \cap f_{\bar{\mu}}^{-1}[Q_{\bar\mu}] \menos N_{\bar{\mu}}^{\ast}]\cup [A \cap \ran g_{\bar{\mu}}].
    $$
    Hence, it is enough to prove that the last set is Borel. Since $g_{\bar{\mu}}^{-1}[A] \in \calB([0, 1])$, by~\ref{b0100f}, $ A \cap \ran g_{\bar{\mu}} = g_{\bar{\mu}}[g_{\bar{\mu}}^{-1}[A]] \in \calB_{T}$. 
    
    Finally, regarding the measure, we have $\Lb(g_{\bar{\mu}}^{-1}[A]) = \Lb(I^{\ast}(A) \menos Q_{\bar{\mu}}) = \Lb(I^{\ast}(A)) = \lambda^{\bar{\mu}}(A)$. 
\end{PROOF}

\begin{remark}\label{b0103}
    The converse in~\autoref{b0100}~\ref{b0100a} does not hold when $N_{\bar{\mu}}^{\ast}$ is uncountable because it contains non-Borel sets that are mapped into the countable $Q_{\bar{\mu}}$. Similarly, the converse in~\autoref{b0100}~\ref{b0100b} is not true when $G^{\ast}_{\bar{\mu}} \neq \emptyset$, because it contains non-Borel subsets whose pre-images are empty. 
\end{remark}

We can also analyze the completion of $\la[T],\calB_T,\lambda^{\bar\mu}\ra$.

\begin{definition}\label{b100}
  Denote by $\la [T],\calL_{\bar\mu},\lambda^{\bar{\mu}}\ra$ the \emph{completion of $\la [T],\Bwf_T, \lambda^{\bar{\mu}}\ra$}.
\end{definition}

In the cases of the Cantor space and the Bairse space, we have: 

\begin{example}
    \
    \begin{enumerate}[label=(\arabic*)]
    \item $\la\cantor,\calL(\cantor),\Lb_2\ra$ is the completion of the measure on $\Bwf(\cantor)$ from the probability tree in \autoref{b036}~\ref{b036-1}. Here, $\calL(\cantor)$ is the \emph{Lebesgue $\sigma$-algebra on the Cantor space} and $\Lb_2$ is the \emph{Lebesgue measure on the Cantor space}. Note that $\Lb_2$ is free (see \autoref{b042}~\ref{b042-1}).
  
    \item $\la\baire,\calL(\baire),\Lb_\omega\ra$ is the completion of the measure on $\Bwf(\baire)$ from the probability tree in \autoref{b036}~\ref{b036-2}. Here, $\calL(\baire)$ is the \emph{Lebesgue $\sigma$-algebra on the Baire space} and $\Lb_\omega$ is the \emph{Lebesgue measure on the Baire space}. Note that $\Lb_\omega$ is free (see \autoref{b042}~\ref{b042-2}).
\end{enumerate}
\end{example}

Since every measurable set can be decomposed as a union of a Borel set and a null set, from~\autoref{b0100} we get: 

\begin{corollary}\label{b101}
    Let $A \subseteq [T]$ and $B \subseteq [0, 1]$. Then:
        \begin{enumerate}[label=\normalfont(\alph*)]
        \item  $A \in \calL_{\bar\mu}$ iff $f_{\bar{\mu}}[A] \in \calL([0, 1])$. In this case, $\Lb(f_{\bar{\mu}}[A]) = \lambda^{\bar{\mu}}(A \cap  V^*_{\bar\mu})$ and $\lambda^{\bar{\mu}}(A) = \Lb(f_{\bar{\mu}}[A]) + \Lb(G^*_A) = \Lb(f_{\bar{\mu}}[A]) +\sum_{x \in A \menos  V^*_{\bar\mu}} \Lb(I_{x}^{\ast})$. 

        \item $B \in \calL([0, 1])$ iff $f_{\bar{\mu}}^{-1}[B] \in \calL_{\bar\mu}$ and $B \cap G^{\ast}_{\bar{\mu}} \in \calL([0, 1])$. In this case, we have that $\lambda^{\bar{\mu}}(f_{\bar{\mu}}^{-1}[B]) = \Lb(B \menos G^{\ast}_{\bar{\mu}})$ and $\Lb(B) = \lambda^{\bar{\mu}}(f_{\bar{\mu}}^{-1}[B]) + \Lb(B \cap G^{\ast}_{\bar{\mu}})$. 

        \item $A \in \calL_{\bar\mu}$ iff $g_{\bar{\mu}}^{-1}[A] \in \calL([0, 1])$. In this case, $\Lb(g_{\bar{\mu}}^{-1}[A]) = \lambda^{\bar{\mu}}(A)$. 

        \item $B \in \calL([0, 1])$ iff $g_{\bar{\mu}}[B] \in \calL_{\bar\mu}$ and $B \cap G^{\ast}_{\bar{\mu}} \in \calL([0, 1])$. In this case, we have that $\Lb(B \menos G^{\ast}_{\bar{\mu}}) = \lambda^{\bar{\mu}}(g_{\bar{\mu}}[B] \cap  V^*_{\bar\mu})$.
    \end{enumerate}
\end{corollary}

Finally, let us consider the case in which $\lambda^{\bar{\mu}}$ is free, i.e.\ when every point in $[T]$ has measure zero. Thanks to \autoref{b077} and \autoref{b045p}~\ref{b045b}, we have the following characterization.

\begin{lemma}\label{b080}
  The following statements are equivalent.
  \begin{enumerate}[label= \normalfont(\roman*)]
     \item\label{b0801} $\lambda^{\bar{\mu}}$ is a free.

    \item\label{b0802} $\Lb(I^*_x) = 0$ for all $x\in [T]$, i.e.\ $ V^*_{\bar\mu}=[T]$.
     
     \item\label{b0803}$\displaystyle\lim_{n\to\infty}\Xi^{\bar\mu}(\{x\frestr n\}) = \prod_{n<|x|}\mu_{x\frestr n}(\{x\frestr(n+1)\}) = 0$ for all $x\in [T]$.
     
     \item\label{b0804} $I^*_x$ is a singleton for all $x\in [T]$.
  \end{enumerate}
\end{lemma}

% \begin{PROOF}{\ref{b080}}
%     We have~\ref{b0801}~$\Leftrightarrow$~\ref{b0802} and ~\ref{b0802}~$\Leftrightarrow$~\ref{b0804} as a consequence of~\autoref{b077}~\ref{b077b} and~\ref{b077c}. Also, ~\ref{b0802}~$\Leftrightarrow$~\ref{b0803} follows from~\autoref{b045p}~\ref{b045b}. 
% \end{PROOF}

As a direct consequence, when $\lambda^{\bar{\mu}}$ is free we get information about the structure of $T_{\bar{\mu}}^{+}$. 

\begin{corollary}\label{pp166}
    If $\lambda^{\bar{\mu}}$ is free then $T_{\bar{\mu}}^{+}$ is a perfect tree.
\end{corollary} 
\begin{PROOF}{\ref{pp166}}
    If $t\in T^+_{\bar\mu}$ then $\lambda^{\bar\mu}([t]_T\cap[T^+_{\bar\mu}])=\lambda^{\bar\mu}([t]_T)>0$, so $[t]_T\cap[T^+_{\bar\mu}]$ contains more than two points. This implies that there are two incompatible nodes in $T^+_{\bar\mu}$ above $t$.
\end{PROOF}

\begin{figure}[h]
    \centering
        \centering
    \begin{tikzpicture}[scale=0.8]
    % Cuadrado izquierdo
    \draw[thick] (0, 0) rectangle (5, 5);
    % Interseccion homeo cuadro derecho
    \fill[redun!15] (10, 0) rectangle (15, 3.5);
    
    % Cuadrado en la parte derecha del rectángulo del cuadrado izquierdo

    % Colocar el texto dentro del cuadrado pequeño
    \node at (4.1, 4.25) {\footnotesize \textbf{$N_{\bar{\mu}}^{\ast}$}};
    
    % Cuadrado derecho
    \draw[thick] (10, 0) rectangle (15, 5);
    
    % interseccion homeomorphismo cuadro derecho
    \fill[greenun!15] (0, 0) rectangle (5, 3.5);
    
    % Subcuadrado con solo el perímetro coloreado (azul) en el cuadrado izquierdo
    \draw[thick,greenun] (0, 0) rectangle (5, 5);

     \draw[thick, black] (0, 0) rectangle (5, 3.5);

    % Colocar el nombre [T] arriba, en el centro del cuadrado izquierdo y más arriba
    \node at (2.5, 5.6) {\textbf{$[T]$}};
    % Colocar el nombre [0, 1] arriba, en el centro del cuadrado derecho y más arriba
    \node at (12.5, 5.6) {\textbf{$[0, 1]$}};
    
    % Flecha curva que representa la función f_{\bar{\mu}}, en la parte superior
    \draw[->, thick] (5.5, 5.5) to[out=30, in=150] node[above] {\footnotesize $f_{\bar{\mu}}$} (9.5, 5.5);
    % Flecha curva que representa la función g_{\bar{\mu}}, en la parte inferior (simétrica)
    \draw[<-, thick] (5.5, -0.5) to[out=-30, in=-150] node[below] {\footnotesize $g_{\bar{\mu}}$} (9.5, -0.5);

    % Colocar el número 1 en la parte derecha
    \node at (14.4, 4.5) {\footnotesize \textbf{$ \bullet \, 1$}};

    \draw[redun, dotted, thick] (13.5, 5) to[out=270, in=180] (15, 3.8);
    
    % define ran g
    \draw[thick] (-1, 0) -- (-1, 1.60);

    \draw[thick] (-1, 1.97) -- (-1, 3.5);

    \draw[thick] (-1, 0) -- (-.88, 0);

    \draw[thick] (-1, 3.5) -- (-0.88, 3.5);
    
    % Dividir la línea a la izquierda del cuadrado izquierdo
    \node at (-1, 1.75) {\footnotesize \textbf{$\mathrm{ran} g_{\bar{\mu}}$}}; % Texto en la línea izquierda

    % Definiendo dom g
    \draw[thick] (16, 0) -- (16, 1.60);
    \draw[thick] (16, 1.97) -- (16, 3.5);
    \node at (16, 1.75) {\footnotesize \textbf{$\mathrm{dom} g_{\bar{\mu}}$}};

    \draw[thick] (16, 0) -- (15.88, 0);
    \draw[thick] (16, 3.5) -- (15.88, 3.5);
    
    % Etiqueta de de f_{\bar{\mu}}[Q] 
    
    \node at (0.9, 4.25) {\textbf{ \footnotesize $f_{\bar{\mu}}^{-1}[Q_{\bar{\mu}}]$}};

    % Etiqueta de Q 
    
    \node at (10.5, 4.25) {\textbf{\footnotesize $Q_{\bar{\mu}}$}};

    % Etiqueta de G

    % Delimitar la inversa de Q 
    
    \draw[thick] (0, 3.5) -- (5, 3.5);

    %\draw (3.3, 5) -- (3.3, 3.5);

     \draw (3.3, 5) to[out=270, in=180] (4, 3.5);

    % Delimitar Q

    \draw[thick] (10, 3.5) -- (15, 3.5);

    % Etiqueta de dom f

    \draw[thick,greenun] (6, 2.83) -- (6, 5);
    \draw[thick,greenun] (6, 0) -- (6, 2.18);
    \node[greenun] at (6, 2.5) {\textbf{\footnotesize $V^*_{\bar\mu}$}};

    \draw[thick,greenun] (6, 5) -- (5.88, 5);
    
    \draw[thick,greenun] (6, 0) -- (5.88, 0);
    % Etiqueta de ran f

    \draw[thick, redun] (9, 2.7) -- (9, 5);
    \draw[thick, redun] (9, 0) -- (9, 2.3);
    \node[redun] at (9, 2.5) {\footnotesize \textbf{$\ran f_{\bar{\mu}}$}};

    \draw[thick, redun] (9, 5) -- (9.12, 5);
    \draw[thick, redun] (9, 0) -- (9.12, 0);

    % definiendo ran f 

    \draw[thick, redun] (10, 0) rectangle (15, 5);

    \draw[thick, black] (10, 0) rectangle (15, 3.5);
\end{tikzpicture}
    \caption{The situation in~\autoref{b090}, that is, when $\lambda^{\bar{\mu}}$ is free. In this case, $ V^*_{\bar\mu} = [T]$, hence $G_{\bar{\mu}} = \emptyset$. The shaded regions are homeomorphic via $f_{\bar{\mu}} \rest \ran g_{\bar{\mu}}$ and its inverse is $g_{\bar{\mu}}$, and the dotted curve indicates that $1$ may or may not belong to $\ran f_{\bar{\mu}}$}
    \label{f52}
\end{figure}
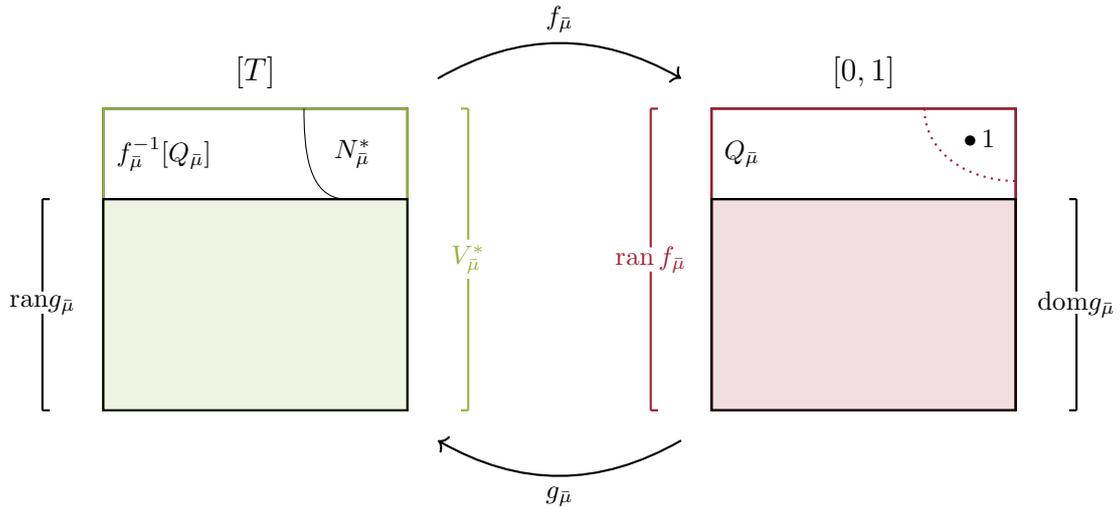

When $\lambda^{\bar{\mu}}$ is free, some properties listed in~\autoref{b087} and~\autoref{b0100} can be simplified.  For example, $ V^*_{\bar\mu} = [T]$ and hence $G_{\bar{\mu}}^{\ast} = \emptyset$ (see~\autoref{f52}). As a consequence: 

\begin{theorem}\label{b090}
    Assume that $\lambda^{\bar{\mu}}$ is free.
    \begin{enumerate}[label=\normalfont(\alph*)]
    
    \item\label{b090b} $f_{\bar\mu}(g_{\bar\mu}(y)) = y$ for $y\in [0,1]\menos Q_{\bar\mu}$, i.e.\ $g_{\bar\mu}$ is one-to-one.
    
    \item\label{b090d} If $x\in [T]$ and $f_{\bar\mu}(x)\notin Q_{\bar\mu}$ then $g_{\bar\mu}(f_{\bar\mu}(x))=x$.
    
    \item\label{b090dd} $\ran g_{\bar\mu} = [T]\menos f^{-1}_{\bar\mu}[Q_{\bar\mu}]$.
    
    \item\label{b090c} $f_{\bar\mu}$ is continuous and $[0,1)\subseteq\ran f_{\bar\mu}$.
    
    \item\label{b090f} $f_{\bar\mu}\frestr \ran g_{\bar\mu}$ is an homeomorphism from $\ran g_{\bar\mu}$ onto $[0,1]\menos Q_{\bar\mu}$ with inverse $g_{\bar\mu}$.

    \item\label{b090g} $\max T\subseteq N^*_{\bar\mu} \subseteq f^{-1}_{\bar\mu}[Q_{\bar\mu}]$. 
    \end{enumerate}
\end{theorem}

\begin{PROOF}{\ref{b090}}

By~\autoref{b087} it just remains to prove that $[0, 1) \subseteq \ran f_{\bar{\mu}}$. For this, by~\ref{b090b}, $[0,1]\menos Q_{\bar\mu} \subseteq \ran f_{\bar\mu}$. Now, if $y\in Q_{\bar\mu}\menos\{1\}$ then, by \autoref{b047}, $y=a_t$ for some $t\in T$. Thus, $y=a_t=f_{\bar\mu}(\concat{t}{\la 0,0,0,\ldots\ra})$.\qedhere %This shows that $[0,1)\subseteq\ran f_{\bar\mu}$. 

\end{PROOF}

From~\autoref{b080} and~\autoref{b090}~\ref{b090c} it follows like in~\autoref{b042}~\ref{b042-1} and~\ref{b042-2} that $Q_{\bar{\mu}}$ is dense when $\lambda^{\bar{\mu}}$ is free. This is actually a characterization of freeness.
%Moreover, we can characterize when $\lambda^{\bar{\mu}}$ is free using $Q_{\bar{\mu}}$. 

\begin{corollary}\label{b095.001}
    $\lambda^{\bar{\mu}}$ is free iff $Q_{\bar{\mu}}$ is dense in $[0, 1]$. 
\end{corollary}

\begin{PROOF}{\ref{b095.001}}
    Assume that $\lambda^{\bar{\mu}}$ is free. Let $a, b \in [0, 1]$ such that $a < b$. Pick some $y \in (a, b)$. By~\autoref{b090}~\ref{b090c}, there is a $x \in [T]$ such that $y = f_{\bar{\mu}}(x)$, hence $\lim_{n \to \infty} a_{x \rest n} = y$. Therefore, we can find an $N < \omega$ such that $a_{x \rest N} \leq y$ and $ y - a_{x \rest N}  < y - a$, that is, $a_{t} \in (a, b)$, where $t \coloneqq X \rest N$. Thus, $Q_{\bar{\mu}}$ is dense. To prove the converse, assume that $\lambda^{\bar{\mu}}$ is not free, that is, there exists some $x \in [T]$ such that $\lambda^{\bar{\mu}}(\{ x \}) > 0$. Let us show that $(f_{\bar{\mu}}^{-}(x), f_{\bar{\mu}}^{+}(x)) \cap Q_{\bar{\mu}} = \emptyset$. Towards contradiction, assume that there exists a $t \in T$ such that $a_{t} \in (f_{\bar{\mu}}^{-}(x), f_{\bar{\mu}}^{+}(x))$. Consider $s \coloneqq x \rest \vert t \vert$, hence $a_{s} \leq f_{\bar{\mu}}^{-}(x) < a_{t} < f_{\bar{\mu}}^{+}(x) \leq b_{s}$, which is not possible. Thus $Q_{\bar{\mu}}$ is not dense in $[0, 1]$. 
\end{PROOF}

As mentioned in~\autoref{b0103}, to establish equivalences in~\autoref{b0100} while preserving the measure there, we face two potential obstacles: on the one hand, $G_{\bar{\mu}}^{\ast}$ may be non-empty, and on the other hand, $N_{\bar{\mu}}^{\ast}$ may include non-measurable subsets. The first issue is solved when $\lambda^{\bar{\mu}}$ is free, and the second by completing the measure, since $N_{\bar{\mu}}^{\ast}$ is null in $[T]$. Thus, under these conditions,~\autoref{b0100} and~\autoref{b101} get simplified.

\begin{theorem}\label{b095}
    Assume that $\lambda^{\bar{\mu}}$ is free, $A \subseteq [T]$ and $B \subseteq [0, 1]$. Then:
    \begin{enumerate}[label=\normalfont(\alph*)]
        \item\label{b095a} $A\in\Bwf_T$ iff $f_{\bar\mu}[A]\in\Bwf([0,1])$ and $A\cap N^*_{\bar\mu}\in\Bwf_T$. In this case, $\Lb(f_{\bar\mu}[A]) = \lambda^{\bar{\mu}}(A)$.
        
        \item\label{b095c} $B\in\Bwf([0,1])$ iff $f^{-1}_{\bar\mu}[B]\in\Bwf_T$. In this case $\lambda^{\bar{\mu}}(f^{-1}_{\bar\mu}[B]) = \Lb(B)$.
        
        \item\label{b105a} 
        $A\in\calL_{\bar\mu}$ iff $f_{\bar\mu}[A]\in\calL(\R)$. In this case, $\Lb(f_{\bar\mu}[A]) = \lambda^{\bar{\mu}}(A)$.
        
        \item\label{b105c} $B\in\calL(\R)$ iff $f^{-1}_{\bar\mu}[B]\in\calL_{\bar\mu}$. In this case $\lambda^{\bar{\mu}}(f^{-1}_{\bar\mu}[B]) = \Lb(B)$.

        \item $B \in \calB([0, 1])$ iff $g_{\bar{\mu}}[B] \in \calB_{T}$. In this case, $\Lb(B) = \lambda^{\bar{\mu}}(g_{\bar{\mu}}[B])$.

        \item $A \in \calB_{T}$ iff $g_{\bar{\mu}}^{-1}[A] \in \calB([0, 1])$ and $A\cap N^*_{\bar\mu}\in\Bwf_T$. In this case,  $\Lb(g_{\bar{\mu}}^{-1}[A]) = \lambda^{\bar{\mu}}(A)$. 

        \item\label{b105a.1} 
        $A\in\calL_{\bar\mu}$ iff $g_{\bar\mu}^{-1}[A] \in \calL(\R)$. In this case, $\Lb(g_{\bar\mu}^{-1}[A]) = \lambda^{\bar{\mu}}(A)$.
        
        \item\label{b105c.1} $B\in\calL(\R)$ iff $g_{\bar\mu}[B]\in\calL_{\bar\mu}$. In this case $\lambda^{\bar{\mu}}(g_{\bar\mu}[B]) = \Lb(B)$.
    \end{enumerate}
\end{theorem}

In addition, we can show a connection between $\la T,\bar\mu\ra$ and the binary probability tree.

\begin{definition}\label{b125}
    For $\alpha\leq\omega$, denote by $\bar 1^\alpha$ the constant sequence of $1$'s of length $\alpha$. Given a representation of $T$, define $h \coloneqq h^T\colon T\to {}^{<\omega}2$ by recursion as follows:
    \begin{itemize}
        \item $h(\la\,\ra) \coloneqq \la\,\ra$;
        \item if $\alpha_t=1$, define $h(\concat{t}{\la 0\ra})\coloneqq h(t)$;
        \item if $1<\alpha_t<\omega$, define $$h(\concat{t}{\la k\ra})\coloneqq 
          \begin{cases}
              \concat{\concat{h(t)}{\bar 1^k}}{\la0\ra} & \text{if $k<\alpha_t-1$,}\\
              \concat{h(t)}{\bar 1^{\alpha_t-1}} & \text{if $k=\alpha_t-1$;}
          \end{cases}$$
        \item if $\alpha_t=\omega$, define $h(\concat{t}{\la k\ra})\coloneqq \concat{\concat{h(t)}{\bar 1^k}}{\la0\ra}$ for all $k<\omega$.
    \end{itemize}
    Let $S \coloneqq S^T$ be the smallest subtree of ${}^{<\omega}2$ containing $\ran h^T$, i.e.\ $s\in S$ iff $s$ is below some node in $\ran h^T$.
\end{definition}

We list below some basic properties of the function $h^T$.

\begin{lemma}\label{b130}
    Let $s,t\in T$. 
    \begin{enumerate}[label = \normalfont (\alph*)]
        \item $s\subseteq t$ implies $h(s)\subseteq h(t)$.
        \item $s\perp t$ iff $h(s)\perp h(t)$.
        \item $h(s)\subseteq h(t)$ iff $[t]_T\subseteq [s]_T$.\footnote{See \autoref{fct:clp}~\ref{clp-c}.}
        \item If $\alpha_t<\omega$ and $t\notin\max(T)$ then $$ [h(t)]_{{}^{<\omega}2} = \bigcup_{t'\in\suc_T(t)} [h(t')]_{{}^{<\omega}2} \text{ (disjoint union).}$$ 
        \item If $\alpha_t=\omega$ and $t\notin\max(T)$ then 
        $$ [h(t)]_{{}^{<\omega}2} = \{\concat{h(t)}{\bar 1^\omega}\}\cup\bigcup_{t'\in\suc_T(t)} [h(t')]_{{}^{<\omega}2} \text{ (disjoint union).}$$
        \item If $t\in\max T$ then $h(t)\in\max S$. 
        \item $h(t)\in\max S$ iff there are no splitting nodes above $t$ in $T$ (including it).
        \item $h(t)\in\max S$ iff $h(t)$ is not a splitting node of $S$. 
        \item Any node in $S$ is either maximal or splitting, and $\max S\subseteq \ran h$.
        \item $T$ is a perfect tree iff $S={}^{<\omega}2$.
        %\item If $x\in [T]$ then there is a unique $y\in[S]$ extending $h(x\frestr n)$ for all $n<\omega$. In particular, if $x\in \max T$ then $y\in\max S$. 
    \end{enumerate}
\end{lemma}

\begin{definition}\label{b131}
    Thanks to \autoref{b130}, we can define a map $e_T\colon [T]\to [S]$ that sends $x\in[T]$ to some $y\in[S]$ extending $h(x\frestr n)$ for all $n<\omega$.
\end{definition}

\begin{lemma}\label{b134}\ 

    \begin{enumerate}[label = \normalfont (\alph*)]
        \item The function $e_T$ is well-defined, i.e.\ the $y$ in \autoref{b131} exists and is unique.
        \item $e_T$ is a topological embedding.
        \item $[S]\menos \ran e_T = \set{\concat{h(t)}{\bar 1^\omega}}{t\in T\menos\max T,\ \alpha_t = \omega}$, which is countable.
        \item $e_T$ is onto iff $T$ is finitely splitting.
    \end{enumerate}
\end{lemma}

\begin{definition}\label{b132}
    Using $\bar\mu$, we define a measure $\Xi^*$ on $\pts(S^T)$ as follows: if $s=h(t)$ for some $t\in T$, set $\Xi^*(\{s\})\coloneqq \Xi^{\bar\mu}(\{t\})$ (this value does not depend on $t$ because $h(t)=h(t')$ implies $[t]_T=[t']_T$), otherwise, if $s\in S\menos\ran h$, pick the largest node $t\in T$ such that $h(t)\subseteq s$, and set $\Xi^*(\{s\})\coloneqq 
    \Xi^{\bar\mu}(\set{t'\in\suc_T(t)}{s\subseteq h(t')}$.
\end{definition}

\begin{theorem}\label{b133}
    $\Xi^*\in\IP$. Moreover, for any $\bar\nu\in\TP$ satisfying $\Xi^*=\Xi^{\bar\nu}$, $A\subseteq[T]$ and $B\subseteq[S]$, 
    \begin{enumerate}[label = \normalfont (\alph*)]
        \item\label{b133a} $Q_{\bar\mu} = Q_{\bar\nu}$ and $g_{\bar\mu} = e_T^{-1}\circ g_{\bar\nu}$.
        \item\label{b133b} $V^*_{\bar\nu} = e_T[V^*_{\bar\mu}]\cup ([S]\menos \ran e_T)$ and $f_{\bar\mu} = f_{\bar\nu}\circ e_T$.
        \item\label{b133c}  $A\in\calL_{\bar\mu}$ iff $e_T[A]\in\calL_{\bar\nu}$, in which case $\lambda^{\bar\nu}(e_T[A]) = \lambda^{\bar\mu}(A)$. 
        \item\label{b133d} $B\in\calL_{\bar\nu}$ iff $e^{-1}_T[B]\in\calL_{\bar\mu}$, in which case $\lambda^{\bar\mu}(e^{-1}_T[B]) = \lambda^{\bar\nu}(B)$. 
        \item\label{b133e} $e_T[N^*_{\bar\mu}] = N^*_{\bar\nu}\cap \ran e_T$.
    \end{enumerate}
\end{theorem}
\begin{PROOF}{\ref{b133}}
    Clearly, $\Xi^*(\la\, \ra)=\Xi^{\bar\mu}(\la\, \ra)=1$. Now let $s\in S\menos\max S$. In the case that $s=h(t)$ for some $t\in T$, $t\notin \max T$, even more, this $t$ can be found as a splitting node of $T$. Then $h(\concat{t}{\la0\ra})= \concat{s}{\la0\ra}$, so $\Xi^*(\{\concat{s}{\la 0\ra}\}) = \Xi^{\bar\mu}(\{\concat{t}{\la 0\ra}\})$, and $\Xi^*(\{h(\concat{t}{\la1\ra})\}) = \sum_{0<k<\alpha_t}\Xi^{\bar\mu}(\{\concat{t}{k}\})$. Thus, $\Xi^*(\suc_S(s)) = \sum_{k<\alpha_t}\Xi^{\bar\mu}(\{\concat{t}{k}\}) = \Xi^{\bar\mu}(\{t\}) = \Xi^*(\{s\})$. 
    
    One can show by recursion that, for $t\in T$, $I^{\bar\nu}_{h(t)} = I^{\bar\mu}_t$. This implies that $Q_{\bar\nu} = Q_{\bar\mu}$, i.e.\ $\dom g_{\bar\nu} = \dom g_{\bar\mu}$. On the other hand, whenever $y\in [S]\menos \ran e_T$, i.e. $y=\concat{h(t)}{\bar 1^\omega}$ for some $t\in T$ with $\alpha_t = \omega$, we have that  $\lambda^{\bar\nu}(\{y\}) = \lim_{n\to\infty} \sum_{k\geq n}\Xi^{\bar\mu}(\{\concat{t}{\la k\ra}\}) = 0$, so $y\in V^*_{\bar\nu}$ and $f_{\bar\nu}(y)=b_{t}$. This shows that $[S]\menos \ran e_T\subseteq f^{-1}_{\bar\nu}[Q_{\bar\nu}]$, i.e.\ $\ran g_{\bar\nu}\subseteq \ran e_T$, which implies that $\dom (e^{-1}_T\circ g_{\bar\nu}) =\dom g_{\bar\nu} = \dom g_{\bar\mu}$. Now, if $z\in\ran g_{\bar \nu}$ then $y\coloneqq g_{\bar\mu}(z)\in \ran e_T$, so $x\coloneqq e_T^{-1}(y)$ is defined, i.e.\ $y$ is the unique extension of $\seq{h(x\frestr n)}{n<\omega}$. Thus, $z\in\bigcap_{n<\omega}I^{\bar\nu}_{h(x\frestr n)} = \bigcap_{n<\omega} I^{\bar\mu}_{x\frestr n}$, which shows that $g_{\bar\mu}(z) = x = e^{-1}_T(g_{\bar\nu}(z))$. This proves~\ref{b133a}. 

    For $x\in[T]$, $\lambda^{\bar\nu}(\{e_T(x)\}) = \lim_{n\to\infty}\Xi^*(\{h(x\frestr n)\}) = \lim_{n\to\infty}\Xi^{\bar\mu}(\{x\frestr n\}) = \lambda^{\bar\mu}(\{x\})$. This proves that $V^*_{\bar\nu} = e_T[V^*_{\bar\mu}]\cup([S]\menos \ran e_T)$, so $\dom (f_{\bar\nu}\circ e_T) = V^*_{\bar\mu} = \dom f_{\bar\mu}$. Now, for $x\in V^*_{\bar\mu}$, $f_{\bar\mu}(x)$ is the unique point in $\bigcap_{n<\omega} I^{\bar\mu}_{x\frestr n} = \bigcap_{n<\omega} I^{\bar\nu}_{h(x\frestr n)} = \{f_{\bar\nu}(e_T(x))\}$. This concludes~\ref{b133b}. 

    \ref{b133c}: Notice that 
    $g^{-1}_{\bar\mu}[A] = g^{-1}_{\bar\nu}[e_T[A]]$. Then, by \autoref{b101}, $A\in\calL_{\bar\mu}$ iff $e_T[A]\in\calL_{\bar\nu}$, in which case $\lambda^{\bar\mu}(A) = \Lb(g^{-1}_{\bar\mu}[A]) = \Lb(g^{-1}_{\bar\nu}[e_T[A]]) = \lambda^{\bar\nu}(e_T[A])$.

    \ref{b133d}: It follows from~\ref{b133c}, also because $[S]\menos \ran e_T \subseteq f^{-1}_{\bar\nu}[Q_{\bar\nu}]$ has measure zero.
\end{PROOF}

\section{Probability trees and the null ideal}\label{7}

In this section, we aim to prove the \emph{invariance} of the \emph{cardinal invariants} associated with the null ideal, mostly via Tukey connections. Although the most general result (\autoref{a243}) is already known, we offer an elementary and direct proof using probability trees. 
Our starting point is to show that, whenever $\lambda \in\BP$ is free, $\calN(\lambda)$ and $\calN([0, 1])$ are Tukey equivalent (see~\autoref{c46} and~\autoref{a241}).

Based on~\cite{vojtas}, we first introduce key concepts related to \emph{relational systems} and \emph{Tukey connections}, which will provide the necessary framework to formalize our results.

\begin{definition}
    We say that $R=\la X, Y, \lhd \ra$ is a \textit{relational system} if $X$, $Y$ are non-empty sets and $\lhd$ is a relation. 
    \begin{enumerate}[label=(\arabic*)]
        \item A set $F\subseteq X$ is \emph{$R$-bounded} if $\exists y\in Y \forall x\in F\ (x \lhd y)$. 
        
        \item A set $E\subseteq Y$ is \emph{$R$-dominating} if $\forall x\in X \exists y \in E\ ( x \lhd y)$. 
    \end{enumerate}
    
    We can associate two cardinal invariants with  relational systems: 
    \begin{itemize}
        \item[{}] $\bfrak(R) \coloneqq \min\{|F| \colon  F\subseteq X  \text{ is }R\text{-unbounded}\}$, the \emph{unbounding number of $R$}, and
        
        \item[{}] $\dfrak(R) \coloneqq \min\{|D| \colon  D\subseteq Y \text{ is } R\text{-dominating}\}$, the \emph{dominating number of $R$}.
    \end{itemize}
\end{definition}

We now present some examples of relational systems and their associated cardinal invariants.  

\begin{example}\label{c37}
    \index{relational system! $\calI$}
    Let $\calI$ be an ideal on a set $X$ containing all its singletons.  

    \begin{enumerate}
        \item  $\calI = \langle \calI, \calI, \subseteq \rangle$ is a relational system, $\gb(\calI) = \add(\calI)$ and $\gd(\calI) = \cof(\calI)$.

        \item $\Cv_{\calI} \coloneqq \langle X, \calI, \in \rangle$ is a relational system, $\gd(\Cv_{\calI}) = \cov(\calI)$ and $\gb(\Cv_{\calI}) = \non(\calI)$.  
    \end{enumerate}
\end{example}

Now we introduce the  notion of \emph{Tukey connections}, which can be thought of as homomorphisms between relational systems: 

\begin{definition}\label{c46}
    Let $R = \la X, Y, \lhd\ra$ and $ S = (Z, W, \sqsubset)$ be relational systems. We say that $(\psi_{-}, \psi_{+}) \colon R \to S$ is a \emph{Tukey connection} from $R$ into $S$ if $\psi_{-} \colon X \to Z$ and $\psi_{+} \colon W \to Y$ are functions such that $$ \forall x \in X  \forall w \in W\ [ \psi_{-}(x) \sqsubset w \Rightarrow x \lhd \psi_{+}(w)]. $$
    
    In this case, we write $R \leqT S$ and we say that \emph{$R$ is Tukey-below  $S$}. When $R \leqT S$ and $S \leqT R$, we say that $R$ and $S$ are \emph{Tukey equivalent} and we denote it by $R \eqT S.$  
\end{definition}

Tukey equivalences determined inequalities between the associated cardinal invariants. 

\begin{lemma}\label{c47}
    Let $R = \langle X, Y, \lhd \rangle$ and $S = \langle Z, W, \sqsubset \rangle$ be two relational systems.
    \begin{enumerate}[label = \normalfont (\alph*)]
        \item If $R \leqT S$ then $\gb(R) \geq \gb(S)$ and $\gd(R) \leq \gd(S)$. 
        \item If $R \eqT S$ then $\gb(R) = \gb(S)$ and $\gd(R) = \gd(S)$. 
    \end{enumerate}
\end{lemma}

To prove the results of this section, we introduce several maps.  Let $\bar{\mu} \in \TP$ and $T\coloneqq T_{\bar\mu}$ (with some representation), and consider $f_{\bar{\mu}}$ as in~\autoref{m45}. We define the maps
    \begin{align*}
        H_{\bar{\mu}}^{-} &\colon \pts([T])\to \pts([0,1]) && A\mapsto f_{\bar{\mu}}[A]\cup \{1\};\\
        H_{\bar{\mu}}^{+}&\colon \pts([0,1])\to \pts([T]) && B\mapsto f_{\bar{\mu}}^{-1}[B].
    \end{align*}
For this section, we extend $g_{\bar\mu}$ to $[0,1]$ in such a way that, for $y\in Q_{\bar\mu}$, $g_{\bar\mu}(y)$ is some $f_{\bar\mu}$-preimage of $y$, if it exists, otherwise $g_{\bar\mu}(y)$ is sent to some arbitrary point in $[T_{\bar\mu}]$.

\begin{lemma}\label{a240-}
    Let $\bar{\mu} \in \TP$, $T\coloneqq T_{\bar\mu}$, $A\subseteq [T]$, $B\subseteq[0,1]$, $x\in[T]$ and $y\in[0,1]$. Then:
    \begin{enumerate}[label= \normalfont (\alph*)]
        \item $H^-_{\bar\mu}(A)\subseteq B$ implies $A\subseteq H^+_{\bar\mu}(B)$.
        \item If $\lambda^{\bar\mu}$ is free, then $H^+_{\bar\mu}(B)\subseteq A$ implies $B\subseteq H^-_{\bar\mu}(A)$.
        \item $f_{\bar\mu}(x)\in B$ implies $x\in H^+_{\bar\mu}(B)$. 
        \item If $y\in\ran f_{\bar\mu}\cup\{1\}$ then $g_{\bar\mu}(y)\in A$ implies $y\in H^-_{\bar\mu}(A)$.
        \item If $\lambda^{\bar\mu}$ is free then $g_{\bar\mu}(y)\in A$ implies $y\in H^-_{\bar\mu}(A)$.
    \end{enumerate}
\end{lemma}
\begin{PROOF}{\ref{a240-}}
    Easy to check, also using that, by~\autoref{b090}~\ref{b090c}, $[0,1)\subseteq \ran f_{\bar\mu}$ whenever $\lambda^{\bar\mu}$ is free. 
\end{PROOF}

As a consequence, for the null ideal we obtain: 
    
\begin{theorem}\label{a241}
    Let $\lambda\in \BP$. If $\lambda$ is free then $\Ncal(\lambda) \eqT \Ncal([0,1])$ and  $\Cv_{\Ncal(\lambda)}\eqT \Cv_{\Ncal([0,1])}.$
\end{theorem}

\begin{PROOF}{\ref{a241}}
    %We prove both items simultaneously. 
    By \autoref{pp129}~\ref{pp129.3}, there is some $\bar\mu\in\TP$ such that $\lambda = \lambda^{\bar\mu}$. By \autoref{b101}, $H^-_{\bar\mu}$ and $H^+_{\bar\mu}$ send measure zero sets into measure zero sets. Thus, by \autoref{a240-}, it is clear that $\Ncal(\lambda) \eqT \Ncal([0,1])$ is witnessed by the pairs $(H^-_{\bar\mu}\frestr\Ncal(\lambda),H^+_{\bar\mu}\frestr\Ncal([0,1]))$ and $(H^+_{\bar\mu}\Ncal([0,1]),H^-_{\bar\mu}\Ncal(\lambda))$, and $\Cv_{\Ncal(\lambda)}\eqT \Cv_{\Ncal([0,1])}$ is witnessed by $(f_{\bar\mu},H^+_{\bar\mu}\frestr\Ncal([0,1]))$ and $(g_{\bar\mu},H^-_{\bar\mu}\frestr\Ncal(\lambda))$.
\end{PROOF}

% As a consequence of~\autoref{c47} and~\autoref{a241}, the cardinal invariants associated with the null ideal are invariant in the following sense. 

% \begin{corollary}
%     If $\lambda\in\BP$ is free then the cardinal invariants associated with $\Ncal(\lambda)$ are the same as $\Ncal([0,1])$ and $\Ncal(\R)$, i.e.\  
%     \begin{align*}
%         \add(\calN(\lambda)) & = \add(\calN([0, 1])) = \add(\calN(\R)), & 
%         \cof(\calN(\lambda)) & = \cof(\calN([0, 1])) = \cof(\calN(\R)),\\
%         \non(\calN(\lambda)) &= \non(\calN([0, 1])) = \non(\calN(\R)), &
%         \cov(\calN(\lambda)) & = \cov(\calN([0, 1])) = \cov(\calN(\R)). 
%     \end{align*}
% \end{corollary}

Although the following generalization of \autoref{a241} is already known, we offer an alternative proof using our results. Recall that a \emph{Polish space} is a completely metrizable separable space, and that a \emph{Borel isomorphism} between topological spaces is a bijection that sends Borel sets into Borel sets in both directions. 

\begin{theorem}[{\cite[Thm.~17.41]{Kechris}}]\label{a243}
    If $X$ is a Polish space and $\mu\colon \calB(X)\to[0,1]$ is a free probability measure, then there is some Borel isomorphism $f\colon X\to [0,1]$ such that, for any Borel $B\subseteq[0,1]$, $\mu(f^{-1}[B]) = \Lb(B)$. In particular, $\Ncal(\mu) \eqT \Ncal([0,1])$ and $\Cv_{\Ncal(\mu)}\eqT \Cv_{\Ncal([0,1])}$.
\end{theorem}
\begin{PROOF}{\ref{a243}}
    Since $\mu$ is a free probability measure on $\Bwf(X)$, $X$ must be uncountable, so there exists a Borel isomorphism $g\colon X\to\baire$. Let $\lambda$ be the probability measure on $\baire$ induced by $X$, i.e.\ $\lambda(A)\coloneqq \mu(g^{-1}[A])$ for $A\in\Bwf(\baire)$. It is clear that $\lambda$ is free, so it is enough to prove the theorem for $\la\baire,\lambda\ra$ instead of $\la X,\mu\ra$. 
    
    Pick some $\bar\mu\in\TP$ such that $\lambda^{\bar\mu} = \lambda$, and pick some infinite countable $W\subseteq \baire\menos f_{\bar\mu}^{-1}[Q_{\bar\mu}]$ and let $W'\coloneqq f_{\bar\mu}[W]$. 
    We proceed by cases. First assume that $N^*_{\bar\mu}$ is countable. Hence, we can construct a function $f\colon \baire\to[0,1]$ such that $f$ equals $f_{\bar\mu}$ at $\baire\menos (f^{-1}_{\bar\mu}[Q_{\bar\mu}]\cup W)$ and $f_{\bar\mu}[f^{-1}_{\bar\mu}[Q_{\bar\mu}]\cup W] = Q_{\bar\mu}\cup W'$. This is indeed a Borel bijection (hence isomorphism) and, by \autoref{b095}, it preserves measure as desired.

    Now consider the case when $N^*_{\bar\mu}$ is uncountable. Let $C_0$ be the ternary Cantor set in $[0,1]$, which is closed and $\Lb(C_0)=0$. Pick a sequence $\seq{J_n}{1\leq n<\omega}$ of pairwise disjoint closed intervals of positive length contained in $(1/3,2/3)$ such that $\max J_n < \min J_{n+1}$ for all $n$, and let $C_n$ be the image of $C_0$ under the canonical homeomorphism from $[0,1]$ onto $J_n$ (i.e\ the line segment from $(0,\min J_n)$ to $(1,\max J_n)$), so $C_n$ resembles the Cantor ternary set in $J_n$, thus it is closed with measure zero. 

    Define the map $h\colon [0,1]\to [0,1]\menos C$ such that it is the identity on $[0,1]\menos \bigcup_{n<\omega} C_n$ and $h\frestr C_n$ is a linear isomorphism onto $C_{n+1}$ for all $n<\omega$. This is a Borel isomorphism that preserves the Lebesgue measure. 
    Finally, define $f\colon \baire\to [0,1]$ such that $f$ coincides with $h\circ f_{\bar\mu}$ at $\baire\smallsetminus (f^{-1}_{\bar\mu}[Q_{\bar\mu}]\cup W))$, $f[(f^{-1}_{\bar\mu}[Q_{\bar\mu}]\menos N^*_{\bar\mu})\cup W] = h[Q_{\bar\mu}\cup W']$, and $f\frestr N^*_{\bar\mu}$ is a Borel isomorphism onto $C_0$. This is the desired Borel isomorphism.
\end{PROOF}

We also have a similar result for the Lebesgue measure on $\R$.

\begin{corollary}\label{a245}
    If $X$ is a Polish space, $\mu\colon \calB(X)\to[0,\infty]$ is a free $\sigma$-finite measure and $\mu(X)=\infty$, then there is some Borel isomorphism $f\colon X\to \R$ such that, for any Borel $B\subseteq\R$, $\mu(f^{-1}[B]) = \Lb(B)$. In particular, $\Ncal(\mu) \eqT \Ncal(\R)$ and $\Cv_{\Ncal(\mu)}\eqT \Cv_{\Ncal(\R)}$.
\end{corollary}
\begin{PROOF}{\ref{a245}}
    Partition $X$ into Borel sets $\seq{B_n}{n\in \Z}$ such that $0<\mu(B_n)<\infty$ for each $n\in\Z$ and $\sum_{n=1}^\infty \mu(B_n) = \sum_{n=-1}^{-\infty} \mu(B_n)=\infty$. Then, we can partition $\R$ into semi-open intervals $\seq{J_n}{n \in\Z}$ such that each $J_n$ has length $\mu(B_n)$. By \autoref{a243}, there is some Borel isomorphism $f_n\colon B_n\to J_n$ preserving measure (this uses that, for any $B\in\Bwf(X)$, there is some finer Polish topology on $X$ such that its Borel $\sigma$-algebra is $\Bwf(X)$ and $B$ is clopen, see e.g.~\cite[Thm.~13.1]{Kechris}). Thus, we obtain the desired $f$ by putting all the $f_n$ together.
\end{PROOF}

\begin{corollary}\label{a246}
    If $X$ is a Polish space, $\mu\colon \calB(X)\to[0,\infty]$ is a free $\sigma$-finite measure and $\mu(X)>0$ then $\Ncal(\mu) \eqT \Ncal(\R)$ and $\Cv_{\Ncal(\mu)}\eqT \Cv_{\Ncal(\R)}$. In particular,
    \begin{align*}
        \add(\calN(\mu)) & = \add(\calN(\R)), & 
        \cof(\calN(\mu)) & =  \cof(\calN(\R)),\\
        \non(\calN(\mu)) & =  \non(\calN(\R)), &
        \cov(\calN(\mu)) & = \cov(\calN(\R)). 
    \end{align*}
\end{corollary}

% Journals of probability 

% \begin{itemize}
%     \item Probability Theory and Related Fields (Q1)  (2.1). 

%     \item Combinatorics, Probability and Computing (Q1) (1.0). 

%     \item Probability Surveys (Q1). (1.05 Favorita.)
    
%     \item Journal of theoretical probability (Q2) (0.8). (Segunda favorita) 

%     \item Theory of Probability and Its Applications (Q3) (0.5). (Tercera favorita)

%     \item Revista Colombiana de Matematicas. Q4, 0.15
%     \end{itemize}

{\small
%\addcontentsline{toc}{section}{References}

%\bibliography{appl}        %activar para modificar bibliografia, luego actualizar el .bbl

\newcommand{\etalchar}[1]{$^{#1}$}

   %bibliografia final, comentar para modificar

\bibliographystyle{alpha}
}

\end{document}